\documentclass[twoside,11pt]{article}

\usepackage[preprint]{jmlr2e}
\hypersetup{hidelinks}

\usepackage[T1]{fontenc}
\usepackage[utf8]{inputenc}
\usepackage[dvipsnames,table]{xcolor}
\usepackage{amsmath,mleftright,mathtools}
\DeclareMathAlphabet\mathbfcal{OMS}{cmsy}{b}{n}
\usepackage{bm}
\usepackage{dsfont}
\newcommand{\mathbbm}[1]{\mathds{#1}}
\usepackage{mathrsfs}
\usepackage{relsize}
\usepackage{textcomp}
\usepackage[caption=false]{subfig}
\usepackage[ruled,vlined,linesnumbered,hangingcomment]{algorithm2e}
\usepackage{comment}
\usepackage{tikz}
\usepackage{bbding}
\usepackage{pifont}
\usepackage{hhline}
\usepackage{listings}
\usepackage{verbatim}
\usepackage{makecell}
\usepackage{tabu}
\usepackage{enumitem}
\usepackage{booktabs}
\usepackage{tabularx}
\usepackage{array}
\usepackage{threeparttable}
\usepackage{lastpage}

\setlist[itemize]{leftmargin=*}
\setlist[enumerate]{leftmargin=*}
\newcolumntype{Y}{>{\raggedright\arraybackslash}X}
\numberwithin{equation}{section}
\numberwithin{theorem}{section}
\renewenvironment{proof}[1][Proof]{\par\noindent{\bfseries #1\ }}{\hfill\BlackBox\\[2mm]}
\renewcommand{\paragraph}[1]{\par\medskip\noindent\textit{#1}\quad}

\newtheorem{assumption}[theorem]{Assumption}
\newtheorem{pro}[theorem]{Proposition}

\newtheorem{coro}[theorem]{Corollary}
\newtheorem{lem}[theorem]{Lemma}

\newtheorem{rem}[theorem]{Remark}

\allowdisplaybreaks[3]

\newcommand{\nn}{\nonumber}

\newcommand{\calC}{\mathcal{C}}

\newcommand{\calH}{\mathcal{H}}

\newcommand{\calK}{\mathcal{K}}
\newcommand{\calL}{\mathcal{L}}
\newcommand{\calM}{\mathcal{M}}
\newcommand{\calN}{\mathcal{N}}

\newcommand{\calP}{\mathcal{P}}

\newcommand{\calS}{\mathcal{S}}

\newcommand{\calU}{\mathcal{U}}
\newcommand{\calV}{\mathcal{V}}
\newcommand{\calW}{\mathcal{W}}

\newcommand{\ba}{\mathbf{a}}
\newcommand{\bA}{\mathbf{A}}

\newcommand{\bB}{\mathbf{B}}

\newcommand{\bC}{\mathbf{C}}
\newcommand{\bd}{\mathbf{d}}
\newcommand{\bD}{\mathbf{D}}

\newcommand{\bE}{\mathbf{E}}

\newcommand{\bg}{\mathbf{g}}
\newcommand{\bG}{\mathbf{G}}
\newcommand{\bh}{\mathbf{h}}

\newcommand{\bI}{\mathbf{I}}

\newcommand{\bP}{\mathbf{P}}

\newcommand{\br}{\mathbf{r}}
\newcommand{\bR}{\mathbf{R}}
\newcommand{\bs}{\mathbf{s}}

\newcommand{\bu}{\mathbf{u}}
\newcommand{\bU}{\mathbf{U}}
\newcommand{\bv}{\mathbf{v}}
\newcommand{\bV}{\mathbf{V}}
\newcommand{\bw}{\mathbf{w}}
\newcommand{\bW}{\mathbf{W}}
\newcommand{\bx}{\mathbf{x}}
\newcommand{\bX}{\mathbf{X}}

\newcommand{\bY}{\mathbf{Y}}

\newcommand{\bSigma}{\bm{\Sigma}}

\DeclareMathOperator{\rank}{rank}
\DeclareMathOperator{\tucker}{Tucrank}

\newcommand{\qednew}{\nobreak \ifvmode \relax \else
      \ifdim\lastskip<1.5em \hskip-\lastskip
      \hskip1.5em plus0em minus0.5em \fi \nobreak
      \vrule height0.75em width0.5em depth0.25em\fi}

\newcommand{\scrU}{\mathscr{U}}

\newcommand{\tcalA}{\mathbfcal{A}}
\newcommand{\tcalB}{\mathbfcal{B}}
\newcommand{\tcalC}{\mathbfcal{C}}
\newcommand{\tcalD}{\mathbfcal{D}}

\newcommand{\tcalH}{\mathbfcal{H}}

\newcommand{\tcalS}{\mathbfcal{S}}
\newcommand{\tcalT}{\mathbfcal{T}}
\newcommand{\tcalU}{\mathbfcal{U}}
\newcommand{\tcalV}{\mathbfcal{V}}

\newcommand{\tcalX}{\mathbfcal{X}}

\newcommand{\tcalZ}{\mathbfcal{Z}}

\DeclareMathOperator{\sign}{sign}
 
\DeclareMathOperator{\rad}{rad}

\DeclareMathOperator{\dist}{dist}

\DeclareMathOperator{\relu}{ReLU}

\ShortHeadings{}{}
\firstpageno{1}
\jmlrheading{}{2026}{1-\pageref{LastPage}}{}{}{}{Junren Chen, Lijun Ding, Dong Xia, and Ming Yuan}

\begin{document}

\title{Low-Rank Tensor Estimation from Nonlinear Observations: \\ A Unified Framework}

\author{\name Junren Chen$^1$, \name Lijun Ding$^2$, \name Dong Xia$^3$, \name Ming Yuan$^1$\\ \addr $^1$Department of Statistics, Columbia University\\ $^2$Department of Mathematics, University of California, San Diego\\ $^3$Department of Mathematics, Hong Kong University of Science and Technology}

\maketitle

\begin{abstract}%
We consider the estimation of a $d_1\times d_2\times d_3$ tensor $\tcalX^\star$ of Tucker rank $(r_1,r_2,r_3)$ from the nonlinear observations $\{y_i=f_i(\langle\tcalA_i,\tcalX^\star\rangle)\}_{i=1}^n$. We develop a unified approach that first constructs a gradient map from the data and then establishes the {\it tensor restricted approximate invertibility condition} (T-RAIC), a condition that quantifies how well the gradient map aligns with the ideal descent step under a low-rank tensor dual norm. We show that T-RAIC yields local linear convergence guarantees for  a Riemannian gradient descent (RGD) algorithm, which may incorporate a normalization step if $\|\tcalX^\star\|_{\rm F}$ is known a priori. Under $O(r_1r_2r_3+\sum_{1\le i\le 3}r_id_i)$ Gaussian measurements, we establish T-RAICs for single-index models, logistic regression (with known $\|\tcalX^\star\|_{\rm F}$), phase retrieval, ReLU regression, and one-bit compressed sensing. The RAICs imply the local linear convergence of RGD to the underlying tensor, exactly in phase retrieval and ReLU regression, and  up to near-optimal estimation errors in the remaining models. We further show that, in all these models except for phase retrieval, a simple spectral initialization yields the desired initialization from $O(d^{3/2})$ measurements under $d_1=d_2=d_3=d$ (ignoring dependence on the Tucker rank and condition number of $\tcalX^\star$). This is also the best known sample complexity for polynomial-time and end-to-end  algorithms in tensor linear regression and tensor completion. Numerical simulations are provided to corroborate our theoretical findings. 
\end{abstract}

\section{Introduction} 
Tensors, or multiway arrays, provide a natural representation of data with multiple modes and have found widespread applications in statistics, machine learning, signal processing, and scientific computing \citep{kolda2009tensor}. Representative examples include hyperspectral and biomedical imaging \citep{li2010tensor}, recommender systems \citep{bi2018multilayer}, network analysis \citep{papalexakis2016tensors}, and longitudinal data analysis \citep{hoff2015multilinear}. In many such applications,   the tensors have very large ambient dimension but exhibit low-rank structures, thus achieving dimension reduction and allowing for estimation from a small number of observations.

In this paper, we consider  the problem of \emph{low Tucker rank} tensor estimation from \emph{nonlinear observations}. The goal is to learn a $d_1\times d_2\times d_3$  tensor $\tcalX^\star$ of low Tucker rank $(r_1,r_2,r_3)$ from   the covariates $\tcalA_1,\cdots,\tcalA_n$ and the corresponding nonlinear responses 
\begin{align}\label{tnlrmodel}
    y_i = f_i(\langle\tcalA_i,\tcalX^\star\rangle),\qquad i\in [n]:= \{1,2,\cdots,n\},
\end{align}
where $f_1,\cdots,f_n$ are independent copies of a possibly random function $f$, and are independent of the covariates. This model is general and encompasses unknown, nonsmooth or even discontinuous $f$. In particular, it covers not only the linear models such as the widely studied tensor PCA and tensor linear regression \citep{luo2023low,luo2024tensor,zhang2018tensor}, but also the following nonlinear problems: 
\begin{itemize}
    \item \textbf{Tensor single-index models.}  The observations are as in (\ref{tnlrmodel}) and the link function $f$ is unknown.
    Since the scale of $\tcalX^\star$ can be absorbed into the unknown link function, we assume $\tcalX^\star$ has unit Frobenius norm to ensure identifiability. Note that single-index models provide a flexible semiparametric extension of linear regression in statistics (see, e.g., \citealp{li1991sliced}). 
    
    \item \textbf{Tensor logistic regression.}  Logistic regression is one of the most basic models  in statistics and machine learning (e.g., \citealp{mccullagh2019generalized}). Under a tensor parameter $\tcalX^\star$, conditioning on the covariates $\tcalA_i$'s, the responses are independent and follow $ y_i \mid \tcalA_i
        \sim \operatorname{Bernoulli}
        \bigl(s(\langle \tcalA_i,\tcalX^\star\rangle)\bigr)$, where $ s(a):= (1+\exp(-a))^{-1}$. 
    This model fits the above formulation by taking
    $f_i(a)=\mathbbm{1}\{u_i\leq s(a)\}$ for independent $u_i\sim\operatorname{Unif}(0,1)$.
    \item \textbf{Tensor phase retrieval.} Set $f_i=|\cdot|$ so that the responses are the magnitudes of the linear measurements: \(\{y_i=\bigl|\langle \tcalA_i,\tcalX^\star\rangle\bigr|\}_{i=1}^n.\)
    Since $\tcalX^\star$ and $-\tcalX^\star$ generate identical observations, the underlying tensor is identifiable only up to a global sign.   Phase retrieval  is motivated by   many real applications in image and signal processing  where the sensors cannot capture the phase information (e.g., \citealp{candes2015phase,shechtman2015phase}). 
    
    \item \textbf{Tensor ReLU regression.} The responses are generated through the rectified linear unit (ReLU) defined as $\relu(u)=\max\{u,0\}=\frac{1}{2}(u+|u|)$, i.e., 
    \(\{y_i=\relu
        \bigl(\langle \tcalA_i,\tcalX^\star\rangle\bigr)\}_{i=1}^n.
    \)
    The ReLU link is  nonsmooth at the origin and censors all negative linear measurements.  ReLU regression can be viewed as the problem of learning a single neuron and has served as a tractable model for understanding the optimization dynamics of shallow neural networks \citep{soltanolkotabi2017learning,yehudai2020learning}. 
    
    \item \textbf{One-bit tensor sensing.}   Only the signs of the linear measurements are observed:
    \(\{y_i=\sign\bigl(\langle \tcalA_i,\tcalX^\star\rangle\bigr)\}_{i=1}^n. 
    \)
    The discontinuous sign link discards all magnitude information, so only the direction of $\tcalX^\star$ is identifiable. We assume $\tcalX^\star$ has unit Frobenius norm to ensure identifiability.  Quantization plays an important role in modern machine learning and signal processing, with the one-bit sensing problem being widely studied as an extreme  case \citep{jacques2013robust,plan2012robust}. 
\end{itemize}

These examples bring together several fundamental problems arising in different areas, under a low Tucker rank tensor parameter. This paper asks the following question: 
\begin{center}
\emph{Can we develop a unified framework that yields computationally efficient  and statistically optimal procedures for these nonlinear low Tucker rank tensor regression problems?}
\end{center}
To the best of our knowledge, computationally efficient  and statistically optimal algorithms are currently unavailable for any of the models above, as we will discuss in Section \ref{sec:relatedwork}.

\subsection{Our Contributions}
This paper answers the above question affirmatively by developing a unified framework for low Tucker rank tensor estimation from nonlinear observations. Its key ingredient is the \emph{tensor restricted approximate invertibility condition} (T-RAIC), which requires a model-specific empirical gradient map  to approximate the ideal descent step under the dual norm with respect to low Tucker rank tensors; see Definition \ref{def:traic}. Note that the T-RAIC is formulated independently of any particular link function or loss function and, in turn, does not require either to be continuous, smooth or convex. Indeed, we shall see that the framework accommodates nonsmooth, discontinuous, and even unknown link functions.

In particular, the T-RAIC is highly effective for the low Tucker rank tensor estimation problem in light of 
its implication on the convergence of   \emph{Riemannian gradient descent} (RGD):  
under a suitable T-RAIC, RGD with a warm initialization linearly converges to a small neighborhood of $\tcalX^\star$ (Theorems \ref{thm:prgdconver}--\ref{thm:nrgdconver}). As such, given a concrete model, 
the core step of our framework is to show the T-RAIC of some gradient map, which is typically chosen as the (sub)gradient  of a loss function.

\paragraph{T-RAICs for concrete models.} We consider Gaussian covariates $\tcalA_i$'s with i.i.d. $N(0,1)$ entries, which was also adopted in the best known prior results on tensor linear regression (e.g., \citealp{zhang2020islet,tong2022scaling,luo2023low}). We establish T-RAICs for single-index models, logistic regression (with   $\|\tcalX^\star\|_{\rm F}=1$), phase retrieval, ReLU regression, and one-bit sensing under a number of observations at the order of \[T_{\rm df}:=r_1r_2r_3+r_1d_1+r_2d_2+r_3d_3.\] Hence, the ``local complexity'' of these problems is dictated by the degrees of freedom of the low Tucker rank tensor cone, reminiscent of some prior works
\citep{tong2022scaling,soltanolkotabi2019structured}. In particular, the T-RAICs imply that RGD locally converges to $\tcalX^\star$, exactly in phase retrieval and ReLU regression, while up to near-optimal errors in single-index models, logistic regression and one-bit compressed sensing. Moreover, in the proofs of the T-RAICs in Appendix \ref{sec:generalraic}, we indeed establish more general RAICs; see Definition \ref{def:raic}. The general RAICs directly imply the T-RAICs and are also of independent interest; see Remark \ref{rem:technical}.

\paragraph{Spectral initialization.} The convergence under T-RAIC is local and therefore requires a warm initialization. We show that a simple first-moment spectral procedure---which forms $\tcalX_0^{(f)}=n^{-1}\sum_{i=1}^n y_i\tcalA_i$ and applies T-HOSVD---suffices for a broad class of nonlinear observations. Interestingly, this is the same initializer used for   tensor linear regression \citep{tong2022scaling}. Its extension to the nonlinear setting is nontrivial and will be discussed in Remark \ref{rem:iniintuition}.   Our method finds accurate initializations under   $O(\sqrt{d_1d_2d_3})$ observations for the above-considered nonlinear models   except for phase retrieval. 


\paragraph{End-to-end algorithms.}
Combining the spectral initialization with the local convergence results yields end-to-end algorithms with complete theoretical guarantees for single-index models, logistic regression, ReLU regression and one-bit compressed sensing with a low Tucker rank tensor parameter. To highlight the leading dimensional dependence, consider the balanced setting $d_1=d_2=d_3=d$ with  Tucker rank and condition number being $O(1)$.  Some of our results for the end-to-end algorithms are stated informally as follows:

\begin{itemize}
    \item \textbf{Tensor single-index models (informal version of Theorem \ref{thm:e2etsim}).} Suppose $f$ satisfies some conditions. 
    For some constant $C_f$ depending on the link function $f$ only,  
    if $n\ge C _f d^{3/2}$, then with high probability (w.h.p.),
    \[
        \|\tcalX_t-\mu_f\tcalX^\star\|_{\rm F}
        \le C_f \big(2^{-t} 
        +
         \sqrt{ T_{\mathrm{df}}/n}\big),\qquad \forall t\ge 0
    \] 
    holds for some $\mu_f>0$ depending on $f.$

    \item \textbf{Tensor ReLU regression (informal version of Theorem \ref{thm:e2erelu}).} If $n\gtrsim d^{3/2}$, then w.h.p., 
    \[
        \|\tcalX_t-\tcalX^\star\|_{\rm F}
        \lesssim
        2^{-t},\qquad \forall t\ge 0. 
    \]  

    \item \textbf{One-bit tensor sensing (informal version of Theorem \ref{thm:e2e1b}).} 
    If $n\gtrsim d^{3/2}$, then w.h.p.,
    \[
        \|\tcalX_t-\tcalX^\star\|_{\rm F}
        \lesssim
        2^{-t}\|\tcalX^\star\|_{\rm F} 
        +
        \frac{T_{\mathrm{df}}}{n}
        \log^{3/2}\Big(\frac{n}{T_{\mathrm{df}}}\Big),\qquad \forall t\ge 0. 
    \] 
\end{itemize}
 Notice that the $O(d^{3/2})$ sample complexity is statistically suboptimal and is incurred in the warm initialization. Nonetheless, it
 agrees with the best-known polynomial-time guarantees for  tensor linear regression and tensor completion \citep{tong2022scaling,luo2023low,xia2019polynomial,xia2021statistically}. 
 We further emphasize that the optimal error rates hinge on the nonlinear observations. For instance, while exact recovery is achieved from ReLU observations, this is not possible from one-bit observations or in single-index models with general unknown link functions.

\subsection{Related Works}\label{sec:relatedwork}
There has been a vast literature for low Tucker rank tensor estimation in different statistical models. The most relevant results can be roughly categorized into the three classes of (i) \emph{tensor linear regression}, (ii) \emph{tensor completion}, and (iii) \emph{tensor PCA with possibly nonlinear observations}, which we separately discuss in the following. We also compare the present paper with the most relevant works in Table \ref{tab:relatedworks}. \emph{For simplicity, we restrict the discussion to  $d_1=d_2=d_3=d,~r_1=r_2=r_3=r$.}   

\paragraph{Tensor linear regression.}
The theoretical literature on low Tucker rank tensor estimation has focused predominantly on linear observation models. Raskutti, Yuan, and Chen \citeyearpar{raskutti2019convex} developed a general convex regularization framework for linear multi-response tensor regression under various low-dimensional structures. Chen, Raskutti, and Yuan \citeyearpar{chen2019non} subsequently studied nonconvex projected gradient descent for generalized low-rank tensor regression, including generalized linear models. However, the two frameworks do not yield statistically optimal results for low Tucker rank tensors. In particular, the results of \citet{raskutti2019convex,chen2019non} only achieved the suboptimal error rates $O(\sqrt{rd^2/n})$, in contrast to the optimal ones that are proportional to $T_{\rm df}=r^3+rd$.

\noindent 
A number of subsequent works developed statistically optimal and computationally efficient methods for tensor regression with  Gaussian linear measurements. Zhang et al. \citeyearpar{zhang2020islet} proposed the ISLET estimator and established minimax optimal guarantees. Tong et al. \citeyearpar{tong2022scaling} established the local convergence of ScaledGD to the true tensor in noiseless tensor regression, and showed that a spectral method gives the desired initialization under $n=O_r(d^{3/2})$. Luo and Zhang \citeyearpar{luo2023low} developed Riemannian Gauss--Newton to achieve local quadratic convergence under noisy linear measurements. More recently, Luo and Zhang \citeyearpar{luo2024tensor} studied tensor-on-tensor regression under rank over-parameterization and the statistical-computational gap. These results  are largely restricted to \emph{(noisy) linear observations}, in contrast to the present paper that focuses on nonlinear observations. 

\paragraph{Tensor completion.} Tensor completion is another extensively studied linear observation model, in which each sample reveals (a possibly noisy version of) a single entry of the underlying tensor. For low Tucker rank tensors, statistically and computationally efficient procedures have been developed under suitable incoherence conditions \citep{xia2019polynomial,xia2021statistically,tong2022scaling}. The technical challenges are complementary to ours: tensor completion must handle the non-Gaussian covariates and the incoherence  issue (which already occurs in matrix completion; e.g., \citealp{davenport2016overview}), whereas we study nonlinear responses under dense Gaussian covariates. As such, the prior tensor completion results are not directly comparable to ours. While some nonlinear counterparts of  tensor completion have also been considered \citep{ghadermarzy2018learning,lee2020tensor}, to the best of our knowledge, a statistically optimal and computationally efficient theory for low Tucker rank tensors has not yet been developed.

\paragraph{(Nonlinear) tensor PCA.} Most of the existing theoretical results for low Tucker rank tensor estimation from nonlinear observations are in tensor PCA, where the entire tensor is observed, i.e., $n = d^{3}$. (The only exception is \citet[Theorem 2]{chen2019non} that encompasses generalized linear models, but as mentioned, the guarantee is suboptimal.) Specifically,  Han, Willett, and Zhang \citeyearpar{han2022optimal} developed a general projected gradient descent framework that yields minimax optimal guarantees for tensor PCA with Poisson or binomial observations.  Cai, Li, and Xia \citeyearpar{cai2023generalized} developed a Riemannian optimization framework for generalized low Tucker rank plus sparse tensor estimation, which applies to (robust) tensor PCA with Poisson or Bernoulli observations. The main distinction from this paper is that these results require $n=d^3$ and do not address the high-dimensional setting with $n\ll d^{3}$. In contrast, for our results, roughly $r^{3}+rd$ observations are sufficient to guarantee local convergence, while (except for phase retrieval) roughly $d^{3/2}$  observations allow for an accurate initialization.

\begin{table}[ht!]
\centering
\begin{threeparttable}
\caption{Comparison with representative works on low Tucker rank tensor estimation under $d_1=d_2=d_3=d$ and $r_1=r_2=r_3=r$. ``Tensor regression'' refers to tensor regression with Gaussian linear measurements.}
\label{tab:relatedworks}
\renewcommand{\arraystretch}{1.18}
\setlength{\tabcolsep}{4.5pt}
\footnotesize
\begin{tabularx}{\textwidth}{
    >{\raggedright\arraybackslash}p{0.22\textwidth}
    >{\raggedright\arraybackslash}p{0.32\textwidth}
    Y
}
\toprule
\textbf{Work}
&
\textbf{Concrete models}
&
\textbf{Main distinction}
\\
\midrule

Raskutti, Yuan, and Chen
\citeyearpar{raskutti2019convex}
&
Linear multi-response tensor regression
&
Convex regularization; suboptimal guarantee for low Tucker rank tensor  
\\
\addlinespace[3pt]

Chen, Raskutti, and Yuan
\citeyearpar{chen2019non}
&
(Generalized) linear models
&
Projected gradient descent; suboptimal guarantee for low Tucker rank tensor   
\\
\addlinespace[3pt]

Zhang et al.
\citeyearpar{zhang2020islet}
&
 Tensor regression  
&
ISLET; minimax optimal guarantee 
\\
\addlinespace[3pt]

Han, Willett, and Zhang
\citeyearpar{han2022optimal}
&
Tensor regression; Sub-Gaussian/Poisson/binomial tensor PCA
&
Projected gradient descent; minimax optimal guarantees; nonlinear models are PCA with $n=d^{3}$   
\\
\addlinespace[3pt]

Tong et al.
\citeyearpar{tong2022scaling}
&
Tensor regression; tensor completion
&
 ScaledGD; optimal guarantees under $n=O(d^{3/2})$ 
\\
\addlinespace[3pt]

Luo and Zhang
\citeyearpar{luo2023low}
&
Tensor regression; sub-Gaussian tensor PCA
&
Riemannian Gauss--Newton;  quadratic convergence; minimax optimal guarantees
\\
\addlinespace[3pt]

Cai, Li, and Xia
\citeyearpar{cai2023generalized}
&Robust 
sub-Gaussian/heavy-tailed/Poisson/Bernoulli   tensor PCA   
&
    Riemannian optimization; low-rank plus sparse models; worked models are PCA with $n=d^{3}$   
\\
\addlinespace[3pt]

Luo and Zhang
\citeyearpar{luo2024tensor}
&
Tensor-on-tensor regression
&
Riemannian optimization; linear observations 
\\

\midrule
\textbf{This work}
&
  Single-index models, logistic regression, phase retrieval, ReLU regression,  one-bit sensing 
&
Riemannian gradient descent; local convergence with $n=O(T_{\mathrm{df}})$; initialization with $n=O(d^{3/2})$   (except for phase retrieval)
\\
\bottomrule
\end{tabularx}
\end{threeparttable}
\end{table}

\paragraph{Nonconvex optimization via RAIC.}  Departing momentarily from tensor, our work is also technically relevant to a line of recent works that use the RAIC to analyze nonconvex optimization algorithms under nonlinear observations.  
This notion was first introduced in Friedlander et al. \citeyearpar{friedlander2021nbiht} to establish the $O(n^{-1})$ error decay rate of the normalized binary iterative hard thresholding (NBIHT) for one-bit compressed sensing. This was   improved to near-optimal error rate by Matsumoto and Mazumdar \citeyearpar{matsumoto2024binary} who established a sharper RAIC, and was then extended to one-bit matrix sensing and more general quantization models in \citet{chen2024optimal}.

\noindent 
More recently, the RAIC has   proven   effective in other nonlinear observation models, including sparse binary generalized linear models \citep{matsumoto2025learning}, logistic regression \citep{chen2026finite} and one-bit phase retrieval \citep{chen2024one}. These works treat the recovery of sparse vectors or low-rank matrices, and do not touch on the estimation of a low Tucker rank tensor which involves additional statistical and computational challenges. Moreover, the RAICs developed in these works are typically tailored to specific models, in contrast to the unified framework herein. 

\subsection{Notation and Organization}
 We use $C,C_i,c,c_i$ to denote positive universal constants. We write $T_1\lesssim T_2$ or $T_1=O(T_2)$ if $T_1\le CT_2$ for some $C$, and write $T_1\gtrsim T_2$ or $T_1=\Omega(T_2)$ if $T_1\ge cT_2$ for some $c$. We use $\tilde{O}(\cdot)$ and $\tilde{\Omega}(\cdot)$ to hide the logarithmic factors of $(n,d_1,d_2,d_3)$ in $O(\cdot)$ and $\Omega(\cdot)$, respectively. In $\mathbb{R}^d$, we denote the unit $\ell_2$ ball by $\mathbb{B}_2^d$ and the unit $\ell_2$ sphere by $\mathbb{S}^{d-1}$. In matrix and tensor spaces, we instead denote the  unit Frobenius sphere by $\mathbb{S}_{\rm F}$, and denote the unit Frobenius ball by $\mathbb{B}_{\rm F}$, with the dimensions being self-evident from the context. We also write $\mathbb{B}_2^d(\bu;R) = \bu+R\mathbb{B}_2^d$ and $\mathbb{B}_{\rm F}(\tcalX^\star;R)=\tcalX^\star+R\mathbb{B}_{\rm F}$.   For random variable $X$,  define \(
\|X\|_{\psi_2}
:=
\inf\big\{
t>0:
\mathbb{E}\exp(\frac{X^2}{t^2})\le 2
\big\}
\) 
and 
\(
\|X\|_{\psi_1}
:=
\inf\big\{
t>0:
\mathbb{E}\exp(\frac{|X|}{t})\le 2
\big\}.
\) For $a\in \mathbb{R}$, we let $\sign(a)=1$ for $a\ge 0$ and $\sign(a)=-1$ for $a<0$. Further notation will be introduced when appropriate.

The rest of this paper is organized as follows. In Section \ref{sec2}, we give the preliminaries and define the T-RAIC. In Section \ref{sec3:genecon}, we show that T-RAIC implies the local linear convergence of RGD and its normalized counterpart. In Section \ref{sec4:examples}, we show that the T-RAICs hold in single-index models, logistic regression, phase retrieval, ReLU regression, and one-bit compressed sensing with Gaussian covariates, and present their implications on the local convergence of the corresponding RGD algorithms. In Sections \ref{sec:speinitial}--\ref{sec:end2end}, we propose an accurate spectral initialization and combine it with the local convergence results to obtain end-to-end algorithms with statistical guarantees. We validate our theoretical findings by numerical simulations in Section \ref{sec:simulation}, and close the paper with concluding remarks in Section \ref{sec:conclusion}. Most of the proofs are deferred to the appendices.

\section{Preliminaries}\label{sec2}

In this section, we provide necessary preliminaries on   tensors (Section \ref{sec:tensor}) and Riemannian optimization (Section \ref{sec:riemannian}),   then formally define the T-RAIC (Section \ref{sec:RAIC}).

\subsection{Tensor and Tucker Rank}\label{sec:tensor}
Suppose $j\in\{1,2,3\}$. The $j$-th matricization $\calM_j(\cdot)$ unfolds the  tensor $\mathbfcal{A}\in \mathbb{R}^{d_1\times d_2\times d_3}$ along mode $j$ and yields $\calM_j(\mathbfcal{A})\in \mathbb{R}^{d_j\times d_j^-}$ where $d_j^- := d_1d_2d_3/d_j$. Specifically, $\calM_1(\mathbfcal{A})\in\mathbb{R}^{d_1\times d_2d_3}$ is formally defined as  
\[[\calM_1(\mathbfcal{A})]_{i_1,(i_2-1)d_3+i_3}=[\mathbfcal{A}]_{i_1,i_2,i_3},\qquad \forall (i_1,i_2,i_3)\in [d_1]\times[d_2]\times[d_3].\]    
The definitions of $\calM_2(\tcalA)\in\mathbb{R}^{d_2\times d_1d_3},\,\calM_3(\tcalA)\in \mathbb{R}^{d_3\times d_1d_2}$ are similar and omitted.

The mode-$1$ multilinear product of $\mathbfcal{A}\in \mathbb{R}^{d_1\times d_2\times d_3}$ with a matrix $\bB\in \mathbb{R}^{r_1\times d_1}$ is a tensor in $\mathbb{R}^{r_1 \times d_2\times d_3}$. It is denoted by $\mathbfcal{A}\times_1\bB$ and formally defined as  
\[[\mathbfcal{A}\times_1\bB]_{i_1,i_2,i_3}=\sum_{k=1}^{d_1}[\mathbfcal{A}]_{k,i_2,i_3}[\bB]_{i_1,k},\qquad \forall (i_1,i_2,i_3)\in [r_1]\times [d_2]\times [d_3].\]
The definition of the mode-$j$ multilinear product $\mathbfcal{A}\times_j\bB$ when $j\in\{2,3\}$ is similar.  

In this paper, we focus on the Tucker rank defined as \[\tucker(\tcalA):=\big(\!\rank(\calM_1(\mathbfcal{A})),\rank(\calM_2(\mathbfcal{A})),\rank(\calM_3(\mathbfcal{A}))\big)^\top.\]    
We say $\mathbfcal{A}$ has Tucker rank at most $\br:=(r_1,r_2,r_3)$ if $\rank(\calM_j(\mathbfcal{A}))\le r_j$ holds for $j\in[3]$, abbreviated as $\tucker(\tcalA)\preceq \br$. Letting $\bd:=(d_1,d_2,d_3)$, 
we shall denote by $T^{\bd}_{\br}$ the set of $d_1\times d_2\times d_3$ tensors of Tucker rank at most $(r_1,r_2,r_3)$.

We say $\tcalA$ is of Tucker-rank $(r_1,r_2,r_3)$ if $\tucker(\tcalA)=(r_1,r_2,r_3)$. Such $\bA$ satisfies a Tucker decomposition 
\[\mathbfcal{A} = \mathbfcal{C}\times_1 \bU_1 \times_2\bU_2\times_3 \bU_3=\mathbfcal{C}\times_{j=1}^3 \bU_j\] 
    for some core tensor $\mathbfcal{C}\in \mathbb{R}^{r_1\times r_2\times r_3}$ and  orthonormal matrices $\bU_j\in \mathbb{O}^{d_j \times r_j}:=\{\bU\in \mathbb{R}^{d_j\times r_j}:\bU^\top\bU=\bI_{r_j}\}$.   We define $T_{\rm df}:=r_1r_2r_3+\sum_{j=1}^3r_jd_j$ to capture the degrees of freedom   of $T^{\bd}_\br$, in light of the Tucker decomposition. 
    Let $\lambda_{\min}(\bA)$ denote the minimal (positive) singular value of a matrix $\bA$, we also define the minimal singular value, maximal singular value, and the condition number of a tensor $\tcalA $ respectively as 
    \[\Lambda_{\min}(\mathbfcal{A}) := \min_{j \in [3]} \lambda_{\min}(\calM_j(\mathbfcal{A})),\quad\Lambda_{\max}(\mathbfcal{A}) := \max_{j \in [3]}\|\calM_j(\mathbfcal{A})\|_{\rm op},\quad  \kappa(\tcalA) := \frac{\Lambda_{\max}(\tcalA)}{\Lambda_{\min}(\tcalA)}.\]

Notice that it is in general NP-hard to compute the exact projection onto $T^{\bd}_{\br}$ \citep{hillar2013most}. We shall frequently utilize the    \emph{Truncated HOSVD} (T-HOSVD) $H^{\bd}_{\br}(\cdot)$ that approximately finds the best Tucker-rank-$(r_1,r_2,r_3)$ approximation, in the sense that \citet[Chapter 10]{hackbusch2012tensor} 
\begin{align*}
    \|H^{\bd}_{\br}(\mathbfcal{A})-\mathbfcal{A}\|_{\rm F}\le \sqrt{3}\|\calP_{T^{\bd}_\br}(\mathbfcal{A})-\mathbfcal{A}\|_{\rm F},\quad \forall\mathbfcal{A}\in \mathbb{R}^{d_1\times d_2\times d_3}.
\end{align*} 
The detailed procedure is given in  Algorithm \ref{alg:thosvd}, where we write ${\rm SVD}_{r}(\bU)$ to denote an orthonormal matrix with the columns being the $r$ leading left singular vectors of $\bU$. We also let $\overline{d}=\max_j d_j$, $\underline{d}=\min_j d_j$ and $\overline{r}=\max_j r_j$. 
  For a more comprehensive introduction to tensors, see  \citep{kolda2009tensor,boumal2023introduction,absil2008optimization}.
 
\begin{algorithm}[ht!]
	\caption{Truncated Higher-order Singular Value Decomposition (T-HOSVD) \label{alg:thosvd}}
	\textbf{Input}: $\tcalT\in\mathbb{R}^{d_1\times d_2\times d_3}$, Tucker Rank $\br=(r_1,r_2,r_3)$

	Compute $\bU_{j,0}={\rm SVD}_{r_j}(\calM_j(\tcalT))$ for $j=1,2,3$.

        \textbf{Output:}  $\widehat{\tcalT}= \tcalT\times_{j=1}^3 (\bU_{j,0}\bU_{j,0}^\top)$

\end{algorithm}

\subsection{Riemannian Optimization}\label{sec:riemannian}
Since Riemannian optimization algorithms have been widely applied to low-rank tensor estimation (e.g., \citealp{luo2023low,luo2024tensor,cai2023generalized}), we only introduce some basic notions that are necessary for the subsequent developments. 
Suppose  $\mathbfcal{A}\in \mathbb{R}^{d_1\times d_2\times d_3}$ of Tucker rank $(r_1,r_2,r_3)$ has Tucker decomposition $\mathbfcal{C}\times_{j=1}^3 \bU_j$, then we denote the tangent space of the {smooth manifold} $\{\tcalA\in\mathbb{R}^{d_1\times d_2\times d_3}:\tucker(\tcalA)=(r_1,r_2,r_3)\}$ at $\tcalA$ by $T(\mathbfcal{A})$. It  can be expressed as \citet{koch2010dynamical}
\begin{align*}
    T(\mathbfcal{A}) = \bigg\{\tilde{\mathbfcal{C}}\times_{j=1}^3 \bU_j
    + \sum_{j=1}^3 \mathbfcal{C}\times_j \tilde{\bU}_j\times _{k\ne j}\bU_k\bigg| \tilde{\mathbfcal{C}}\in \mathbb{R}^{r_1\times r_2\times r_3},\ 
    \tilde{\bU}_j\in \mathbb{R}^{d_j\times r_j},\ 
    \tilde{\bU}_j^\top \bU_j=0,\ j\in[3]\bigg\}. 
\end{align*}
Note that $T(\mathbfcal{A})$ is a linear subspace contained in $T^{\bd}_{2\br}$. Suppose that   $\mathbfcal{H}(\tcalU)$ is the gradient at $\tcalU$, we   define the Riemannian gradient at $\tcalU$
as its projection onto $T(\tcalU)$:
\[
    \mathbfcal{H}^{\rm (r)}(\mathbfcal{U})=\calP_{T(\mathbfcal{U})}( \mathbfcal{H}(\mathbfcal{U})).
\]  
The Riemannian gradient can be efficiently computed since the projection onto $T(\mathbfcal{A})$ has a closed-form expression \citep{koch2010dynamical,luo2023low}.

A Riemannian optimization algorithm typically consists of two steps: optimization over the tangent space and retraction to the smooth manifold. To be concrete, we use the   T-HOSVD in Algorithm \ref{alg:thosvd} as the retraction.

\subsection{Tensor Restricted Approximate Invertibility Condition}\label{sec:RAIC}
We now define the key component of our framework, the \emph{tensor restricted approximate invertibility condition} (T-RAIC) for the gradient map $\tcalH:\mathbb{R}^{d_1\times d_2\times d_3}\to \mathbb{R}^{d_1\times d_2\times d_3}$. In short, the T-RAIC requires that $\tcalH(\tcalU)$ uniformly approximates the ideal descent step $\tcalU-\tcalX^\star$ over some constraint set $\scrU$,\footnote{To estimate $\tcalX^\star$, $\tcalU-\tcalX^\star$ is the ideal descent step  at $\tcalU$ in view of $\tcalU-(\tcalU-\tcalX^\star)=\tcalX^\star$.} under the low Tucker rank tensor dual norm defined as 
\[\big\|\tcalU\big\|_{(T^{\bd}_{2\br}\cap \mathbb{B}_{\rm F})^\circ} := \sup_{\tcalV\in T^{\bd}_{2\br}\cap \mathbb{B}_{\rm F}}\langle\tcalU,\tcalV\rangle,\qquad \forall\tcalU \in \mathbb{R}^{d_1\times d_2\times d_3}.  
\]


\begin{definition}[Tensor RAIC] \label{def:traic}We say that the gradient map $\tcalH:\mathbb{R}^{d_1\times d_2\times d_3}\to \mathbb{R}^{d_1\times d_2\times d_3}$ satisfies the tensor RAIC with respect to Tucker rank $(r_1,r_2,r_3)$ over a constraint set $\scrU\subset \mathbb{R}^{d_1\times d_2\times d_3}$ for a target tensor $\tcalX^\star\in T^{\bd}_{\br}$ up to an error function $R:\mathbb{R}^{d_1\times d_2\times d_3}\to \mathbb{R}$, if it holds that
    \[\big\|\tcalU-\tcalX^\star-\tcalH(\tcalU)\big\|_{(T^{\bd}_{2\br}\cap \mathbb{B}_{\rm F})^\circ}\le R(\tcalU),\qquad\forall \tcalU\in\scrU.\]
    We denote this by $\tcalH(\tcalU)\sim {\rm RAIC}\big(T^{\bd}_{2\br}\cap \mathbb{B}_{\rm F};\scrU,\tcalX^\star,R(\tcalU)\big)$.
\end{definition}

    We shall see that T-RAIC is a special case of the general RAIC in Definition \ref{def:raic}. 

\section{Riemannian Gradient Descent}\label{sec3:genecon}
Suppose that the gradient  map $\tcalH(\tcalU)$ has been chosen for the problem at hand, then the  Riemannian counterpart is given by 
   $\tcalH^{\rm (r)}(\tcalU)=\calP_{T(\tcalU)}(\tcalH(\tcalU))$. Also, 
   due to the computational infeasibility of finding the best low Tucker rank approximation to a tensor, we instead leverage the T-HOSVD in Algorithm \ref{alg:thosvd} to enforce the low Tucker rank structure.  Taken collectively, we reach the following Riemannian gradient descent  algorithm. 

\begin{algorithm}[ht!]
	\caption{Riemannian Gradient Descent (RGD) \label{alg:prgd}}
	\textbf{Input}: Gradient map $\tcalH:\mathbb{R}^{d_1\times d_2\times d_3}\to\mathbb{R}^{d_1\times d_2\times d_3}$, step size $\eta$, initialization $\tcalX_0$, Tucker rank $(r_1,r_2,r_3)$

	\textbf{For}
	$t = 0, 1, 2,\cdots$  \textbf{do}:
     \[
      \mathbfcal{X}_{t+1}=H_{\br}^{\bd}\big(\mathbfcal{X}_t - \eta\cdot \mathbfcal{H}^{\rm (r)}(\mathbfcal{X}_t)\big);  
\]

        \textbf{Output:}  $\{\tcalX_t\}_{t\ge 0}$

\end{algorithm}

In the following,
we give a general local convergence guarantee of Algorithm \ref{alg:prgd} under the T-RAIC satisfying the conditions that (i) the constraint set is large enough to contain the intersection of a small Frobenius norm ball centered at $\tcalX^\star$ and $T^{\bd}_{\br}$ (see (\ref{RGDconstraintset}) below), and (ii) the error function takes the form $\frac{\|\tcalU-\tcalX^\star\|_{\rm F}}{4(\sqrt{3}+1)}+\mu_2$ for some sufficiently small $\mu_2$.

\begin{theorem}
    [Convergence of Algorithm \ref{alg:prgd}] \label{thm:prgdconver}
   Assume 
   \[\eta\tcalH(\tcalU)\sim {\rm RAIC}\bigg(T^{\bd}_{2\br}\cap \mathbb{B}_{\rm F};\scrU,\tcalX^\star,\frac{\|\tcalU-\tcalX^\star\|_{\rm F}}{4(\sqrt{3}+1)}+\mu_2\bigg)\]
 where $\tcalX^\star$ obeys $\tucker(\tcalX^\star)=(r_1,r_2,r_3)$ and  $\Lambda_{\min}:=\Lambda_{\min}(\tcalX^\star)>0$, and for some $R_{\rm loc}\in(0, \infty]$
 \begin{align}
     \scrU\supset T^{\bd}_\br \cap \mathbb{B}_{\rm F}(\tcalX^\star;R_{\rm loc}).\label{RGDconstraintset}
 \end{align}
     Suppose that $\{\tcalX_t\}_{t=0}^\infty$ is obtained by running Algorithm \ref{alg:prgd} with 
     \begin{align}
         \tcalX_0\in T^{\bd}_\br \cap \mathbb{B}_{\rm F}\bigg(\tcalX^\star; \min\bigg\{R_{\rm loc},\frac{\Lambda_{\min}}{12(\sqrt{3}+1)}\bigg\}\bigg).\label{rgdinicon}
     \end{align} There exist universal constants $C>c>0$ such that, if $\mu_2\le c\min\{R_{\rm loc},\Lambda_{\min}\}$, then 
     \[
         \|\tcalX_t-\tcalX^\star\|_{\rm F} \le\frac{\Lambda_{\min}}{2^t}+C\mu_2,\qquad \forall t\ge 0.
     \]
\end{theorem}
\begin{proof}
    See Appendix \ref{app:provergdconverge}.  
\end{proof}

\subsection{Known Frobenius Norm}
If $\|\tcalX^\star\|_{\rm F}=\beta$ is known  a priori, then we can rescale the norms to ensure that the iterates remain at $\beta \mathbb{S}_{\rm F}$. 
This leads to Algorithm \ref{alg:nprgd} which we refer to as  
 normalized RGD (NRGD).

\begin{algorithm}[ht!]
	\caption{Normalized Riemannian Gradient Descent (NRGD) \label{alg:nprgd}}
	\textbf{Input}: Gradient map $\tcalH:\mathbb{R}^{d_1\times d_2\times d_3}\to\mathbb{R}^{d_1\times d_2\times d_3}$, step size $\eta$, initialization $\tcalX_0$, Tucker rank $(r_1,r_2,r_3)$,  Frobenius norm $\beta=\|\tcalX^\star\|_{\rm F}$

	\textbf{For}
	$t = 0, 1, 2,\cdots$  \textbf{do}:
     \[
      \mathbfcal{X}_{t+1}=\beta\cdot \frac{H_{\br}^{\bd}\big(\mathbfcal{X}_t - \eta\cdot \mathbfcal{H}^{\rm (r)}(\mathbfcal{X}_t)\big)}{\big\|H_{\br}^{\bd}\big(\mathbfcal{X}_t - \eta\cdot \mathbfcal{H}^{\rm (r)}(\mathbfcal{X}_t)\big) \big\|_{\rm F}}; 
\]

        \textbf{Output:}  $\{\tcalX_t\}_{t\ge 0}$

\end{algorithm}

 The merit of incorporating the Frobenius norm rescaling step is that NRGD converges under an RAIC over a constraint set being contained in $\beta\mathbb{S}_{\rm F}$ (see (\ref{NRGDconstraintset}) below), which is weaker than (\ref{RGDconstraintset}).

\begin{theorem}
    [Convergence of Algorithm \ref{alg:nprgd}] \label{thm:nrgdconver} 
  Assume 
  \[\eta\tcalH(\tcalU)\sim{\rm RAIC}\bigg(T^{\bd}_{2\br}\cap \mathbb{B}_{\rm F};\scrU,\tcalX^\star,\frac{\|\tcalU-\tcalX^\star\|_{\rm F}}{8(\sqrt{3}+1)}+\mu_2\bigg) \]
   where $\tcalX^\star$ obeys $\tucker(\tcalX^\star)=(r_1,r_2,r_3)$,  $\Lambda_{\min}:=\Lambda_{\min}(\tcalX^\star)>0$, and $\|\tcalX^\star\|_{\rm F}=\beta$, and for some $R_{\rm loc}\in(0, \infty]$, 
   \begin{align}
       \label{NRGDconstraintset}
       \scrU\supset T^{\bd}_\br \cap \mathbb{B}_{\rm F}(\tcalX^\star;R_{\rm loc})\cap \beta \mathbb{S}_{\rm F}.
   \end{align}
     Suppose that $\{\tcalX_t\}_{t=0}^\infty$ is obtained by running Algorithm \ref{alg:nprgd} with \begin{align}
         \tcalX_0\in T^{\bd}_\br \cap \mathbb{B}_{\rm F}\bigg(\tcalX^\star; \min\bigg\{R_{\rm loc},\frac{\Lambda_{\min}}{24(\sqrt{3}+1)}\bigg\}\bigg)\cap \beta \mathbb{S}_{\rm F}. \label{nrgdinicon}
     \end{align} There exist universal constants $C>c>0$ such that, if $\mu_2\le c\min\{R_{\rm loc},\Lambda_{\min}\}$, then 
     \begin{align}
         \|\tcalX_t-\tcalX^\star\|_{\rm F} \le\frac{\Lambda_{\min}}{2^t}+C\mu_2,\qquad \forall t\ge 0.
     \end{align}
\end{theorem}
\begin{proof}
    See Appendix \ref{app:proof32}.
\end{proof}

\section{Applications to Concrete Models}\label{sec4:examples}
In this section, we establish the T-RAICs in concrete tensor estimation models with nonlinear observations. We treat Gaussian covariates as a theoretical benchmark.
\begin{assumption}
    \label{gaussianAi}
    $\tcalA_1,\tcalA_2,\cdots,\tcalA_n$ have i.i.d. $N(0,1)$ entries.  
\end{assumption}

\subsection{Tensor Single-Index Models}\label{sec:sim}
  We consider single-index models with observations
\begin{align}\label{tsim}
    y_i=f_i(\langle\tcalA_i,\tcalX^\star\rangle),\qquad i\in[n],
\end{align}
under possibly unknown link functions $f_i$'s.
Note that we assume $\|\tcalX^\star\|_{\rm F}=1$ since the norm of $\tcalX^\star$ is absorbed into the possibly unknown $f_i$'s.
 \begin{assumption}\label{ass:unitnorm}
$\|\tcalX^\star\|_{\rm F}=1$ and $\tucker(\tcalX^\star)=(r_1,r_2,r_3)$. 
\end{assumption}

Moreover, we assume the following on the nonlinear transforms $f_i$'s.
\begin{assumption} \label{assump:fi}
 $f_1,\cdots,f_n$ are independent of $\{\tcalA_i\}_{i=1}^n$ and are i.i.d. copies of some (possibly) random function $f$. 
\end{assumption}

Our treatment follows the $\ell_2$ loss minimization approach to single-index models 
\citep{plan2016generalized,oymak2017fast,thrampoulidis2015lasso,genzel2023unified}. As such, we   make the following assumptions on the function $f$, where   (\ref{positivemu}) is   inherited from \citet{plan2016generalized,oymak2017fast}, and (\ref{psi2bound}) is     assumed in \citet{genzel2023unified} to allow for sub-Gaussian concentration. 


\begin{assumption}\label{assumptionmu} Let  $f$ be the (possibly) random function in Assumption \ref{assump:fi} that is independent of $g\sim N(0,1)$. There exist some $\mu>0$ and $\sigma\ge 0$ such that
   \begin{gather}
     \mathbbm{E}\big(gf(g)\big)= \mu,\label{positivemu}\\
     \|f(g)-\mu g\|_{\psi_2}\le \sigma. \label{psi2bound}
\end{gather}  
\end{assumption}

Under the above assumptions, we follow \citet{plan2016generalized,oymak2017fast} and adopt the $\ell_2$ loss
\[\calL(\tcalU)=\frac{1}{2n}\sum_{i=1}^n\big(\langle\tcalA_i,\tcalU\rangle-y_i\big)^2,\] 
whose gradient is 
\begin{align}\label{hsim}
   \tcalH(\tcalU):= \nabla \calL(\tcalU)=\frac{1}{n}\sum_{i=1}^n(\langle\tcalA_i,\tcalU\rangle-y_i)\tcalA_i.
\end{align}  
\begin{rem}
    In this remark, we show that minimizing the $\ell_2$ loss estimates $\mu\tcalX^\star$ rather than $\tcalX^\star$. At the population level,
    \(
        \nabla\mathbb{E}\calL(\tcalU) = \mathbb{E}\nabla\calL(\tcalU) = \mathbb{E}\big(\langle\tcalA_i,\tcalU\rangle - f_i(\langle\tcalA_i,\tcalX^\star\rangle)\big)\tcalA_i= \tcalU - \mathbb{E}[f_i(\langle\tcalA_i,\tcalX^\star\rangle)\tcalA_i].
    \)  
    Moreover, by rotational invariance of $\tcalA_i$, 
    \begin{align}\label{zeromeansim}
        \mathbb{E}[f_i(\langle\tcalA_i,\tcalX^\star\rangle)\tcalA_i] = \mathbb{E}\big[f_i(\langle\tcalA_i,\tcalX^\star\rangle)\langle\tcalA_i,\tcalX^\star\rangle \big]\tcalX^\star\stackrel{(\ref{positivemu})}{=}  \mu\tcalX^\star.  
    \end{align}
    Hence, $\nabla\mathbb{E}\calL(\tcalU)=\tcalU-\mu\tcalX^\star$. Combining with the convexity, $\mu\tcalX^\star$ is the unique minimizer of the population loss. 
\end{rem}

In the following, we show that $\tcalH (\tcalU)$ in (\ref{hsim}) satisfies the  T-RAIC. 

\begin{theorem}[T-RAIC for Tensor Single-Index Models]\label{pro:traicsim} Under Assumptions \ref{gaussianAi}--\ref{assumptionmu}, 
there exist universal constants $C>c>0$ such that  $\tcalH (\tcalU)$ in (\ref{hsim}) satisfies the following. If $n\ge CT_{\rm df}$, then with probability at least $1-C\exp(-cT_{\rm df})$,  
     \begin{align}
         \tcalH (\tcalU)\sim {\rm RAIC}\bigg(T^{\bd}_{2\br}\cap \mathbb{B}_{\rm F};T^{\bd}_{\br},\mu\tcalX^\star,\frac{\|\tcalU-\mu\tcalX^\star\|_{\rm F}}{20}+ C\sigma \sqrt{\frac{T_{\rm df}}{n}}\bigg). \label{traicforsim}
     \end{align}
\end{theorem}
\begin{proof}
    See Appendix \ref{app:raicforsim}.  
\end{proof}

Here and hereafter, we simply write $\Lambda_{\min}=\Lambda_{\min}(\tcalX^\star)$ and $\kappa=\kappa(\tcalX^\star)$ when we consider the recovery of a fixed $\tcalX^\star$. 

\begin{theorem}
    [Tensor single-index models: local convergence] \label{thm:tsim} Under Assumptions \ref{gaussianAi}--\ref{assumptionmu}, let $\tcalH(\tcalU)$ be as in (\ref{hsim}), we obtain $\{\tcalX_t\}_{t=0}^\infty$ by running
    \[ \mathbfcal{X}_{t+1}= H^\bd_{\br}\Big(\mathbfcal{X}_t - \calP_{T(\mathbfcal{X}_t)}\big( \tcalH(\tcalX_t)\big)\Big),\qquad t\ge0.\]
There exist constants $C_{\mu,\sigma}>c_{\mu,\sigma}>0$ that depend on $\mu,\sigma$ only and some universal constants $C>c>0$, such that the following holds. If $\tcalX_0\in T^{\bd}_{\br}\cap \mathbb{B}_{\rm F}\Big(\mu\tcalX^\star;\frac{\mu\Lambda_{\min}}{12(\sqrt{3}+1)}\Big)$ and 
\(n\ge \frac{C_{\mu,\sigma}T_{\rm df}}{\Lambda_{\min}^2},\)
then with probability at least $1-C\exp(-cT_{\rm df})$, 
\begin{align}\label{tsimconvergence}
    \|\tcalX_t-\mu\tcalX^\star\|_{\rm F} \le \frac{\mu\Lambda_{\min}}{2^t} + C\sigma \sqrt{\frac{T_{\rm df}}{n}},\qquad \forall t\ge 0. 
\end{align}  
\end{theorem}
    \begin{proof}
The minimal singular value of $\mu\tcalX^\star$ is $\mu\Lambda_{\min}$. By Theorem \ref{pro:traicsim}, (\ref{traicforsim}) holds with probability at least $1-C\exp(-cT_{\rm df})$. Since $1/20\le [4(\sqrt{3}+1)]^{-1}$, this gives the RAIC required in Theorem \ref{thm:prgdconver}, with $R_{\rm loc}=\infty$ and
$\mu_2=C\sigma\sqrt{T_{\rm df}/n}$. The sample size condition in Theorem \ref{thm:tsim} ensures
$\mu_2\le c\mu\Lambda_{\min}$ after adjusting $C_{\mu,\sigma}$. Applying Theorem \ref{thm:prgdconver} gives \eqref{tsimconvergence}. 
\end{proof}
\begin{rem}
    The error rate in Theorem  \ref{thm:tsim} is in general minimax optimal. In fact, let $f_i$'s be i.i.d. copies of a random function $f:a\mapsto a+\varepsilon$ where $\epsilon\sim N(0,\hat{\sigma}^2)$, then 
    the single-index models reduce to tensor linear regression $y_i=\langle\tcalA_i,\tcalX^\star\rangle+\epsilon_i,~i\in[n]$, for which (\ref{positivemu}) holds with $\mu=1$ and (\ref{psi2bound}) holds with $\sigma\asymp \hat{\sigma}$. Our Theorem \ref{thm:tsim} implies the squared Frobenius norm error rate $O(\hat{\sigma}^2{T_{\rm df}}/{n})$, matching the minimax lower bound in \citet[Theorem 5]{zhang2020islet} for tensor linear regression. 
\end{rem}

\subsection{Tensor Logistic Regression}\label{sec:tenlogi}
We consider tensor logistic regression with the responses $y_i$'s as in Assumption \ref{ass:logi}. 

\begin{assumption}
    \label{ass:logi}
      Conditioning on $\tcalA_i$'s, the responses $y_1,\cdots,y_n$ are independent and follow $y_i\mid \tcalA_i \sim {\rm Bernoulli}(s(\langle\tcalA_i,\tcalX^\star\rangle))$, where $s(a)=(1+\exp(-a))^{-1}$ is the sigmoid function.
\end{assumption}

The goal of this section is to give an initial treatment of tensor logistic regression and show that our framework yields nontrivial convergence and statistical guarantees. To this end, we treat a simplified setting where
  $\|\tcalX^\star\|_{\rm F}=1$ is known a priori so that  Algorithm \ref{alg:nprgd} applies. We note that the current techniques directly extend to known  $\|\tcalX^\star\|_{\rm F}$ that need not equal $1$, but 
  are not sufficient for  unknown $\|\tcalX^\star\|_{\rm F}$, which is left for future work. It is natural to consider the gradient of the logistic loss 
\begin{align}\label{grainglm}
    \tcalH(\tcalU) = \frac{1}{n}\sum_{i=1}^n \big(s(\langle\tcalA_i,\tcalU\rangle)-y_i\big)\tcalA_i. 
\end{align}  This gradient map satisfies the following RAIC.

\begin{theorem}[T-RAIC in tensor logistic regression] 
    \label{traicglm}
    Under Assumptions \ref{gaussianAi}, \ref{ass:unitnorm} and \ref{ass:logi}, consider $\tcalH(\tcalU)$ in (\ref{grainglm}) and $\eta:=[\mathbbm{E}_{g\sim N(0,1)}(s'(g))]^{-1}$ as in Theorem \ref{thm:raicglm}, there exist universal constants $C>c>0$ such that the following holds. If $n\ge CT_{\rm df}$, then with probability at least $1-C\exp(-cT_{\rm df})$,    
    \begin{align}
        \eta \tcalH(\tcalU)\sim {\rm RAIC} \bigg(T^{\bd}_{2\br}\cap\mathbb{B}_{\rm F};T^{\bd}_{\br}\cap \mathbb{S}_{\rm F},\tcalX^\star, C\sqrt{\frac{T_{\rm df}}{n}}\bigg). \label{raicforlogi}
    \end{align}  
\end{theorem}
\begin{proof}
    See Appendix \ref{app:tenlogi}. 
\end{proof}
\begin{theorem}[Tensor Logistic Regression: local convergence] \label{thm:tlogistic}
In the setting of Theorem \ref{traicglm}, suppose that we run 
\[
    \tcalX_{t+1}= \frac{H^{\bd}_{\br}\Big(\tcalX_t - \eta \cdot \calP_{T(\tcalX_t)}\Big(\tcalH(\tcalX_t)\Big)\Big)}{\Big\|H^{\bd}_{\br}\Big(\tcalX_t - \eta \cdot \calP_{T(\tcalX_t)}\Big(\tcalH(\tcalX_t)\Big)\Big)\Big\|_{\rm F}},\qquad t\ge 0
\]
to obtain $\{\tcalX_t\}_{t= 0}^\infty$. There exist some universal constants $C>c>0$ such that the following holds. If $n\ge \frac{C}{\Lambda_{\min}^2}T_{\rm df}$ and $\tcalX_0\in T^{\bd}_{\br}\cap \mathbb{B}_{\rm F}(\tcalX^\star; \frac{\Lambda_{\min}}{24(\sqrt{3}+1)})\cap \mathbb{S}_{\rm F}$, then with probability at least $1-C\exp(-cT_{\rm df})$, 
    \begin{align} \label{tlogiconvergence}
        \|\tcalX_t-\tcalX^\star\|_{\rm F}\le \frac{\Lambda_{\min}}{2^t} + C\sqrt{\frac{T_{\rm df}}{n}},\qquad \forall t\ge 0. 
    \end{align} 
\end{theorem}
 \begin{proof}
  Theorem \ref{traicglm} gives the RAIC in (\ref{raicforlogi})
with probability at least $1-C\exp(-cT_{\rm df})$. Thus the hypothesis of Theorem \ref{thm:nrgdconver} holds with $\beta=1$, $R_{\rm loc}=\infty$, $\mu_1=0$ and
$\mu_2=C\sqrt{T_{\rm df}/n}$. The sample size condition in Theorem \ref{thm:tlogistic} ensures $\mu_2\le c\Lambda_{\min}$. Thus,  \eqref{tlogiconvergence} follows from Theorem \ref{thm:nrgdconver}. 
 \end{proof}
\begin{rem}[Optimality] The upper bound $O(\sqrt{T_{\rm df}/n})$ is minimax optimal in light of existing matching lower bound \citep[Corollary 1]{taki2023structured}.
\end{rem}

\subsection{Tensor Phase Retrieval} \label{sec:NPR}
We consider the tensor phase retrieval problem, which concerns the recovery of   $\pm\tcalX^\star$ from the phaseless observations \begin{align}
     y_i = |\langle\tcalA_i,\tcalX^\star\rangle|,\qquad i\in[n].\label{tpryi}
\end{align}
  In general, we shall measure the estimation error by $\dist(\tcalU,\tcalX^\star)=\min\{\|\tcalU-\tcalX^\star\|_{\rm F},\|\tcalU+\tcalX^\star\|_{\rm F}\}$, while in subsequent development we  assume $\|\tcalX_0-\tcalX^\star\|_{\rm F}\le c\|\tcalX^\star\|_{\rm F}$ (without loss of generality) so that we can simply work with the regular $\ell_2$ distance.

We shall consider the amplitude-based loss function  \citep{wang2017solving,zhang2017nonconvex} \[\calL_{\tcalX^\star}(\tcalU) = \frac{1}{2n}\sum_{i=1}^n \big(|\langle\tcalA_i,\tcalU\rangle|-y_i\big)^2 
\] 
with a specific subgradient
\begin{align}
    \tcalH_{\tcalX^\star}(\tcalU):= \frac{1}{n}\sum_{i=1}^n \big(|\langle\tcalA_i,\tcalU\rangle|-y_i\big)\sign(\langle\tcalA_i ,\tcalU\rangle)\tcalA_i 
    . \label{tprgradient}
\end{align}

\begin{theorem}
    [T-RAIC for phase retrieval]\label{traicpr} Under Assumption \ref{gaussianAi} and (\ref{tpryi}), there exist universal constants $C>c>0$ such that, if $n\ge CT_{\rm df}$, then with probability at least $1-C\exp(-cn)$ the following holds: uniformly for all $\tcalX^\star\in T^{\bd}_{\br}\setminus\{0\}$, $\tcalH_{\tcalX^\star}(\tcalU)$ in (\ref{tprgradient}) satisfies \begin{align}
        \tcalH_{\tcalX^\star}(\tcalU)\sim {\rm RAIC} \bigg(T^{\bd}_{2\br}\cap \mathbb{B}_{\rm F};T^{\bd}_{\br}\cap \mathbb{B}_{\rm F}\big(\tcalX^\star;c\|\tcalX^\star\|_{\rm F}\big),\tcalX^\star,\frac{1}{20}\|\tcalU-\tcalX^\star\|_{\rm F}\bigg).\label{raicforpr}
    \end{align} 
\end{theorem}
\begin{proof}
    See Appendix \ref{app:raicpr}. 
\end{proof}

Using the RAIC and the arguments for proving Theorem \ref{thm:prgdconver}, we arrive at the local contraction guarantee for RGD using the gradients in (\ref{tprgradient}). Note that the guarantee holds \emph{uniformly for all} $\tcalX^\star \in T^{\bd}_{\br}\setminus\{0\}$, due to the uniformity of the T-RAIC in Theorem \ref{traicpr}.

\begin{theorem}[Tensor phase retrieval: local convergence] \label{thm:TPR}
In the setting of Theorem \ref{traicpr}, we recover $\tcalX^\star\in T^{\bd}_{\br}$ from $\{\tcalA_i\}_{i=1}^n$ and $\{y_i=|\langle\tcalA_i,\tcalX^\star\rangle|\}_{i=1}^n$ by running 
\begin{align}
     \tcalX_{t+1}=H^{\bd}_{\br}\Big(\tcalX_t-\calP_{T(\tcalX_t)}\big(\tcalH_{\tcalX^\star}(\tcalX_t)\big)\Big),\qquad \forall t\ge 0,\label{TPRRGDiterate}
\end{align}
which yields $\{\tcalX_t\}_{t\ge 0}$. There exist some universal constants $C>c>0$ such that, if $n\ge CT_{\rm df}$, then with probability at least $1-C\exp(-cn)$, the following holds uniformly for all $\tcalX^\star\in T^{\bd}_{\br}$: if $\tucker(\tcalX^\star)=(r_1,r_2,r_3)$ and $\tcalX_0\in T^{\bd}_{\br}\cap \mathbb{B}_{\rm F}(\tcalX^\star;c\Lambda_{\min}(\tcalX^\star))$, then the corresponding $\{\tcalX_t\}_{t\ge 0}$ satisfies 
\[
    \|\tcalX_t - \tcalX^\star\|_{\rm F} \le \frac{\Lambda_{\min}(\tcalX^\star)}{2^t},\qquad \forall t\ge 0.\]
\end{theorem}
\begin{proof} By Theorem \ref{traicpr}, with probability at least $1-C\exp(-cn)$, (\ref{raicforpr}) holds uniformly for all $\tcalX^\star\in T^{\bd}_{\br}\setminus\{0\}$. Since $\frac{1}{20}<\frac{1}{4(\sqrt{3}+1)}$, this gives the RAIC required by Theorem \ref{thm:prgdconver} with $\mu_2=0$. For each $\tcalX^\star$ of Tucker rank $(r_1,r_2,r_3)$, we invoke Theorem \ref{thm:prgdconver} to yield the claim. In view of (\ref{raicforpr}), we set $R_{\rm loc}=c_1\|\tcalX^\star\|_{\rm F}$ for some $c_1>0$ in Theorem \ref{thm:prgdconver}.
Choose the constant $c$ in Theorem \ref{thm:TPR} sufficiently small. Since 
 $\Lambda_{\min}(\tcalX^\star)\le\|\tcalX^\star\|_{\rm F}$, $\tcalX_0\in T^{\bd}_{\br}\cap \mathbb{B}_{\rm F}(\tcalX^\star;c\Lambda_{\min}(\tcalX^\star))$ in Theorem \ref{thm:TPR} suffices to ensure (\ref{rgdinicon}) in Theorem \ref{thm:prgdconver}. 
\end{proof} 

\subsection{Tensor ReLU Regression}
We consider the recovery of $\tcalX^\star\in T^{\bd}_{\br}$ from i.i.d. Gaussian $\{\tcalA_i\}_{i=1}^n$ and \begin{align}
    y_i=\relu(\langle\tcalA_i,\tcalX^\star\rangle),\qquad i\in[n]. \label{treluyi}
\end{align} We adopt the squared $\ell_2$ loss on the ReLU observations  
\[\calL_{\tcalX^\star}(\tcalU) = \frac{1}{n}\sum_{i=1}^n\big(\!\relu(\langle\tcalA_i,\tcalU\rangle)-y_i\big)^2,\]
for which a generalized gradient is given by
\begin{align}
    \tcalH_{\tcalX^\star}(\tcalU) = \frac{1}{n}\sum_{i=1}^n\big[\relu(\langle\tcalA_i,\tcalU\rangle)-y_i\big]\big[1+\sign(\langle\tcalA_i,\tcalU\rangle)\big]\tcalA_i 
    . \label{relutengra}
\end{align}

\begin{theorem}[T-RAIC for ReLU regression] \label{thm:traicrelu} Under Assumption \ref{gaussianAi} and (\ref{treluyi}), there exist universal constants $C>c>0$ such that the following holds. If $n\ge CT_{\rm df}$, then with probability at least $1-C\exp(-cn)$, for all $\tcalX^\star\in T^{\bd}_{\br}\setminus\{0\}$, $\tcalH_{\tcalX^\star}(\tcalU)$ in (\ref{relutengra}) satisfies \begin{align}
    \tcalH_{\tcalX^\star}(\tcalU)\sim {\rm RAIC} \bigg(T^{\bd}_{2\br}\cap\mathbb{B}_{\rm F};T^{\bd}_{\br}\cap \mathbb{B}_{\rm F}(\tcalX^\star;c\|\tcalX^\star\|_{\rm F}),\tcalX^\star,\frac{\|\tcalU-\tcalX^\star\|_{\rm F}}{20}\bigg).\label{raicforrelu}
\end{align}
\end{theorem}
\begin{proof}
    See Appendix \ref{app:raicrelu}. 
\end{proof}

In the following, we present our local convergence guarantee for tensor ReLU regression.  
As with Theorem \ref{thm:TPR}, the local convergence guarantee holds uniformly for all $\tcalX^\star\in T^{\bd}_{\br}$ of Tucker rank $(r_1,r_2,r_3)$, due to the uniformity of the T-RAIC in Theorem \ref{thm:traicrelu} over $\tcalX^\star\in T^{\bd}_{\br}\setminus\{0\}$. 

\begin{theorem}[Tensor ReLU regression: local convergence] \label{thm:trelulocal}
In the setting of Theorem \ref{thm:traicrelu}, we recover $\tcalX^\star\in T^{\bd}_{\br}$   from $\{\tcalA_i\}_{i=1}^n$ and $\{y_i=\relu(\langle\tcalA_i,\tcalX^\star\rangle)\}_{i=1}^n$ by running
\begin{align}
       \tcalX_{t+1}=H^{\bd}_{\br}\left(\tcalX_t-\calP_{T(\tcalX_t)}\big(\tcalH_{\tcalX^\star}(\tcalX_t)\big)\right),\qquad \forall t\ge 0,\label{reluRGD}
\end{align}
which yields $\{\tcalX_t\}_{t\ge 0}$. There exist universal constants $C>c>0$ such that, if $n\ge CT_{\rm df}$, then with probability at least $1-C\exp(-cn)$, the following holds uniformly for all $\tcalX^\star\in T^{\bd}_{\br}$: if $\tucker(\tcalX^\star)=(r_1,r_2,r_3)$ and $\tcalX_0 \in T^{\bd}_{\br} \cap \mathbb{B}_{\rm F}(\tcalX^\star;c\Lambda_{\min}(\tcalX^\star))$, then the corresponding $\{\tcalX_t\}_{t\ge 0}$ satisfies 
    \begin{align}
        \label{treluconver}
        \|\tcalX_t-\tcalX^\star\|_{\rm F}\le \frac{\Lambda_{\min}(\tcalX^\star)}{2^t},\qquad \forall t\ge 0.
    \end{align} 
\end{theorem}
\begin{proof}
By Theorem \ref{thm:traicrelu},  we   assume that (\ref{raicforrelu}) holds for all $\tcalX^\star\in T^{\bd}_{\br}\setminus\{0\}$ with the promised probability. Since $\frac{1}{20}<\frac{1}{4(\sqrt{3}+1)}$, for each $\tcalX^\star \in T^{\bd}_{\br}$ of Tucker rank $(r_1,r_2,r_3)$,   the RAIC required in 
  Theorem \ref{thm:prgdconver} holds with $\mu_2=0$ and $R_{\rm loc}=c_1\|\tcalX^\star\|_F$ for some universal constant $c_1>0$. Choose sufficiently small $c$ in Theorem \ref{thm:trelulocal}, then $\tcalX_0 \in T^{\bd}_{\br} \cap \mathbb{B}_{\rm F}(\tcalX^\star;c\Lambda_{\min}(\tcalX^\star))$ ensures (\ref{rgdinicon}) in Theorem \ref{thm:prgdconver}. Since $\mu_2\le c\min\{R_{\rm loc},\Lambda_{\min}(\tcalX^\star)\}$ holds trivially, Theorem \ref{thm:prgdconver} applies and yields the claim. 
 \end{proof}

\subsection{One-Bit Tensor Sensing} \label{sec:1btsensing}
We consider the  recovery of $\tcalX^\star\in T^{\bd}_{\br} \cap \mathbb{S}_{\rm F}$ from 
\begin{align}\label{t1bcsyi}
    y_i = \sign(\langle\tcalA_i,\tcalX^\star\rangle),\qquad i\in[n].
\end{align}
This is referred to as the one-bit tensor sensing problem, and is the low Tucker rank tensor counterpart of the extensively studied problem of one-bit compressed sensing \citep{jacques2013robust,plan2012robust,boufounos20081,matsumoto2024binary} and halfspace learning \citep{long1995sample,zhang2020efficient}. Here, $\|\tcalX^\star\|_{\rm F}=1$ is assumed to render the well-posedness of the problem, as the tensors $\{\lambda\tcalX^\star:\lambda>0\}$ are all indistinguishable from the one-bit measurements.

Following the near-optimal algorithm for 1-bit compressed sensing (of sparse vector) \citet{matsumoto2024binary,chen2024optimal}, we adopt the ReLU loss (a.k.a. one-sided $\ell_1$ loss) \[\calL_{\tcalX^\star}(\tcalU) = \frac{1}{n}\sum_{i=1}^n \big(|\langle\tcalA_i,\tcalU\rangle| - y_i\langle\tcalA_i,\tcalU\rangle\big).\]  
Substituting $y_i=\sign(\langle\tcalA_i,\tcalX^\star\rangle)$, a specific subgradient of $\calL_{\tcalX^\star}(\tcalU)$ is given by 
\begin{align}\label{1bcsgrad}
    \tcalH_{\tcalX^\star}(\tcalU):= \frac{1}{n}\sum_{i=1}^n \big[\sign(\langle\tcalA_i,\tcalU\rangle)-y_i\big]\tcalA_i. 
\end{align}

\begin{theorem}
    [T-RAIC for one-bit sensing] \label{traic1bcs} 
 Under Assumption \ref{gaussianAi} and (\ref{t1bcsyi}), there exist universal constants $C>c>0$ such that, 
    if $n\ge CT_{\rm df}$, then with probability at least $1-C\exp(-cT_{\rm df})$, uniformly for all $\tcalX^\star \in T^{\bd}_{\br}\cap \mathbb{S}_{\rm F}$,    $\tcalH_{\tcalX^\star}(\tcalU)$ in (\ref{1bcsgrad}) satisfies \begin{align}
        \sqrt{\frac{\pi}{2}}\tcalH_{\tcalX^\star}(\tcalU)\sim {\rm RAIC} \bigg(T^{\bd}_{2\br}\cap\mathbb{B}_{\rm F};T^{\bd}_{\br}\cap \mathbb{S}_{\rm F},\tcalX^\star,\frac{\|\tcalU-\tcalX^\star\|_{\rm F}}{40}+\frac{CT_{\rm df}}{n}\log^{3/2}\Big(\frac{n}{T_{\rm df}}\Big)\bigg).\label{tensorraic1bcs}
    \end{align}
\end{theorem}
\begin{proof}
    See Appendix \ref{app:raic1bcs}. 
\end{proof}
 
We are ready to present our local convergence guarantee for one-bit tensor sensing, which holds uniformly for all $\tcalX^\star$ of Tucker rank $(r_1,r_2,r_3)$. 

\begin{theorem}
    [One-bit tensor sensing] 
    \label{thm:t1bcs} 
   Under Assumption \ref{gaussianAi}, we recover  $\tcalX^\star\in T^{\bd}_{\br}\cap \mathbb{S}_{\rm F}$ from $\tcalA_i$'s and $\{y_i=\sign(\langle\tcalA_i,\tcalX^\star\rangle)\}_{i=1}^n$ by running 
    \begin{align}\label{1btsnrgd}
    \mathbfcal{X}_{t+1}= \frac{H^{\bd}_{\br}\big(\mathbfcal{X}_t - \calP_{T(\mathbfcal{X}_t)}\big(\sqrt{\frac{\pi}{2}}\tcalH_{\tcalX^\star}(\tcalX_t)\big)\big)}{\big\|H^{\bd}_{\br}\big(\mathbfcal{X}_t - \calP_{T(\mathbfcal{X}_t)}\big(\sqrt{\frac{\pi}{2}}\tcalH_{\tcalX^\star}(\tcalX_t)\big)\big)\big\|_{\rm F}},\qquad t\ge 0,
\end{align}
which yield $\{\tcalX_t\}_{t\ge 0}$. There exist some universal constants $C>c>0$ such that, if $n\ge CT_{\rm df}$,  then with probability at least $1-C\exp(-cT_{\rm df})$, for any $\tcalX^\star\in \mathbb{S}_{\rm F}$ of Tucker rank $(r_1,r_2,r_3)$, if $\tcalX_0\in T^{\bd}_{\br}\cap \mathbb{S}_{\rm F}\cap \mathbb{B}_{\rm F}\big(\tcalX^\star;\frac{\Lambda_{\min}(\tcalX^\star)}{32(\sqrt{3}+1)}\big)$ and $\Lambda_{\min}(\tcalX^\star)\ge \frac{CT_{\rm df}}{n}\log^{3/2}(\frac{n}{T_{\rm df}})$,
    then the corresponding $\{\tcalX_t\}_{t\ge 0}$ satisfies
    \begin{align}\label{1btsconvergence}
        \|\tcalX_t-\tcalX^\star\|_{\rm F}\le \frac{\Lambda_{\min}(\tcalX^\star)}{2^t} + \frac{CT_{\rm df}}{n}\log^{3/2}\Big(\frac{n}{T_{\rm df}}\Big),\qquad \forall t\ge 0.
    \end{align}  
\end{theorem}
\begin{proof}
By Theorem \ref{traic1bcs}, we assume that (\ref{tensorraic1bcs}) holds for all $\tcalX^\star \in T^{\bd}_{\br}\cap \mathbb{S}_F$ with the promised probability. For each $\tcalX^\star \in T^{\bd}_{\br}\cap \mathbb{S}_F$ of Tucker rank $(r_1,r_2,r_3)$ and obeying $\Lambda_{\min}(\tcalX^\star) \ge \frac{CT_{df}}{n}\log^{3/2}\frac{n}{T_{df}}$, we invoke Theorem \ref{thm:nrgdconver} with $\beta =1$. Since $\frac{1}{40}<\frac{1}{8(\sqrt{3}+1)}$,   the T-RAIC required in Theorem \ref{thm:nrgdconver} holds with $\mu_2=\frac{CT_{df}}{n}\log^{3/2}\frac{n}{T_{df}}$. Note that we can set $R_{\rm loc}=\infty$,   hence  $\tcalX_0\in T^{\bd}_{\br}\cap \mathbb{S}_{\rm F}\cap \mathbb{B}_{\rm F}\big(\tcalX^\star;\frac{\Lambda_{\min}(\tcalX^\star)}{32(\sqrt{3}+1)}\big)$ implies (\ref{nrgdinicon}), and $\mu_2\le c\min\{R_{\rm loc},\Lambda_{\min}\}$ also holds. Therefore, Theorem \ref{thm:nrgdconver} applies and yields the claim for each   $\tcalX^\star$ under consideration. 
\end{proof}
\begin{rem}[Optimality]
Since $T^{\bd}_{\br}$ contains $r_1r_2r_3$, $r_id_i$-dimensional linear subspaces ($1\le i\le 3$), \citet[Theorem 1]{chen2024optimal} implies a lower bound $\Omega(T_{\rm df}/n)$ on the uniform recovery error in one-bit tensor sensing. The error rate in  Theorem \ref{thm:t1bcs} matches this lower bound, up to logarithmic factors. 
\end{rem}

We close this section by discussing our technical contributions. 

\begin{rem}[Technical contributions]
    \label{rem:technical} 
  The T-RAICs in this section are direct outcomes of the corresponding general RAICs   established in Appendix \ref{sec:generalraic}. Technically, while \citet{soltanolkotabi2019structured,soltanolkotabi2017learning} treated ReLU regression and phase retrieval by implicitly establishing the RAICs, our Theorems \ref{thm:prraic}, \ref{thm:raicrelusparse} are established through different arguments, involving an in-depth hyperplane tessellation result (Lemma \ref{lem:globinary}). Consequently, our RAICs hold uniformly for all signals of interest, which is stronger than those implied by the techniques in  \citet{soltanolkotabi2019structured,soltanolkotabi2017learning}, which only hold for a fixed underlying signal. This improvement leads to the uniform local convergence guarantees in Theorems \ref{thm:TPR}, \ref{thm:trelulocal}. Also, in subsequent work, \citet{chen2026instance} build upon Theorems \ref{thm:prraic}, \ref{thm:raicrelusparse} to establish instance optimal sparse recovery guarantees. We believe the general RAICs are of independent interest to some other high-dimensional estimation problems. 
\end{rem}  

\section{Spectral Initialization} \label{sec:speinitial}  
This section develops a   spectral initialization that is accurate under a general class of nonlinear observations. For some $\rho>0$, we assume the underlying tensor parameter satisfies 
\begin{align}
    \tcalX^\star\in T^{\bd}_{\br}\cap \rho \mathbb{S}_{\rm F}.
\end{align}
For $j\in\{1,2,3\}$, write
\(
 \lambda_j:=\lambda_{\min}\bigl(\calM_j(\tcalX^\star)\bigr),
 ~
 s_j:=\|\calM_j(\tcalX^\star)\|_{\rm op}.\)
We abbreviate
\(
 \Lambda_{\min}:=\Lambda_{\min}(\tcalX^\star),~
 \Lambda_{\max}:=\Lambda_{\max}(\tcalX^\star),~
 \kappa:=\kappa(\tcalX^\star).\)

As in Assumption \ref{assump:fi}, we consider observations \(\{y_i=f_i\bigl(\langle \tcalA_i,\tcalX^\star\rangle\bigr)\}_{i=1}^n\), where 
the functions $f_1,\ldots,f_n$ are independent copies of a possibly random function $f$, and are independent of the Gaussian covariates $\tcalA_1,\ldots,\tcalA_n$.  A deterministic link is included as the degenerate case (in which $f$ is deterministic). We write $D:=d_1d_2d_3,~d_{-j}:=D/d_j$ and make the following Assumption. Note that (\ref{Cpoint2}) is an adaptation of Assumption \ref{assumptionmu} under $\|\tcalX^\star\|_{\rm F}=\rho$, with the sub-Gaussianity being directly imposed on $f(\rho g)$. 

\begin{assumption} \label{assump:initial}
   The true tensor satisfies $\tucker(\tcalX^\star)=(r_1,r_2,r_3)$ and $\|\tcalX^\star\|_{\rm F}=\rho$ for some $\rho>0$. Let $g\sim N(0,1)$ be independent of $f$, for some $\mu,\nu>0$ we have
   \begin{equation}
 \frac{1}{\rho}\,\mathbb E\bigl[gf(\rho g)\bigr]=\mu,
 \qquad
 \bigl\|f(\rho g)\bigr\|_{\psi_2}\leq \nu. \label{Cpoint2}
\end{equation}
\end{assumption}

We now propose a simple HOSVD initialization in Algorithm \ref{alg:ini}.

\begin{algorithm}[ht!]
    \caption{HOSVD Initialization from Nonlinear Observations \label{alg:ini}}
    \textbf{Input}: $\{(\tcalA_i,y_i)\}_{i=1}^n$ and $\br=(r_1,r_2,r_3)$

    Compute  
    \[
        \tcalX_0^{(f)}:=\frac{1}{n}\sum_{i=1}^n y_i\tcalA_i.
    \]

    \textbf{Output}: 
    $\displaystyle \tcalX_0= H^{\bd}_{\br}(\tcalX_0^{(f)}).$
\end{algorithm}

\begin{rem} 
\label{rem:iniintuition}
The same first-moment HOSVD procedure was used in
\citet{tong2022scaling} for tensor linear regression with Gaussian
covariates (see also similar initializer and analysis in \citealp{zhang2020islet}). In fact, our proof builds on a Gaussian
orthogonal decomposition argument that also appears in the proof of \citet[Lemma 11]{tong2022scaling}. 
However, our proof is not a straightforward adaptation due to 
additional difficulties  under nonlinear observations. Specifically, under nonlinear and therefore non-Gaussian responses, the key technical tools in their proof, including tensor RIP and Lemma 23 therein, no longer apply. 
Instead, we control the corresponding random processes via the powerful
concentration inequality in Lemma \ref{lem:product_process} for product process. 
\end{rem}

  In the following, we provide the theoretical guarantees. 
 

\begin{theorem}[Performance of Algorithm \ref{alg:ini}]\label{thm:initsim}
There exist universal constants $C>c>0$ such that the following holds.  Suppose that
$n\geq CT_{\mathrm{df}}$ and
\begin{equation}
 n\geq C\Big(\frac{\nu}{\mu}\Big)^2
 \max_{1\leq j\leq 3}
 \bigg\{
 \frac{d_j+\sqrt D}{\lambda_j^2},
 \frac{s_j^2d_j}{\lambda_j^4}
 \bigg\}.\label{Cpoint6}
\end{equation}
Then, with probability at least $1-C\exp(-c\underline d)$, the output of Algorithm \ref{alg:ini} satisfies
\begin{equation}
 \|\tcalX_0-\mu \tcalX^\star\|_{\rm F}
 \leq C\nu\sqrt{\frac{T_{\mathrm{df}}}{n}}
 +C\rho\sum_{j=1}^3
 \bigg\{
 \frac{\nu^2(d_j+\sqrt D)}{\mu\lambda_j^2n}
 +\frac{\nu s_j}{\lambda_j^2}\sqrt{\frac{d_j}{n}}
 \bigg\}. \label{Cpoint7}
\end{equation}
\end{theorem}

Theorem \ref{thm:initsim} yields the following corollary, which is often more convenient to apply. 
 
\begin{coro}\label{corCpoint1}
For every fixed $c_0\in(0,1)$, there exists a constant $C_{c_0}>0$ such that 
\begin{equation}
 n\geq C_{c_0}\Big(\frac{\nu}{\mu}\Big)^2
 \bigg\{
 \frac{T_{\mathrm{df}}}{\Lambda_{\min}^2}
 +\frac{\rho\sqrt D}{\Lambda_{\min}^3}
 +\frac{\rho^2\kappa^2\overline d}{\Lambda_{\min}^4}
 \bigg\} \label{Cpoint8}
\end{equation}
implies
\begin{equation}
 \|\tcalX_0-\mu \tcalX^\star\|_{\rm F}\leq c_0\mu\Lambda_{\min},\qquad  \left\|\frac{\tcalX_0}{\|\tcalX_0\|_{\rm F}}-\frac{\tcalX^\star}{\rho}\right\|_{\rm F}
 \leq 2c_0\frac{\Lambda_{\min}}{\rho} \label{Cpoint9}
\end{equation}
with probability at least $1-C\exp(-c\underline d)$. 
Moreover, condition \eqref{Cpoint8} is implied by 
\begin{equation}
 n\geq C_{c_0}\Big(\frac{\nu}{\mu\rho}\Big)^2
 \bigg\{
 \kappa^3\overline r^{3/2}\sqrt D
 +\kappa^6\overline r^2\overline d
 +\kappa^2\overline r^4
 \bigg\}. \label{Cpoint11}
\end{equation}
\end{coro}

 
\begin{rem} \label{rem:mpr}
In phase retrieval with $f(a)=|a|$, the first condition in (\ref{Cpoint2}) is not satisfied, hence Algorithm \ref{alg:ini} is not able to provide an accurate initialization. However,
to the best of our knowledge, provable initialization procedures in phase retrieval are currently available only for unstructured or sparse signals (e.g., \citealp{candes2015phase,cai2016optimal}). For phase retrieval of $d_1\times d_2$ low-rank matrices, the initialization problem already becomes substantially more challenging, and no efficient algorithm is known to achieve provable accuracy when $n \leq d_1d_2$ \citep[Section 4.2]{lee2021phase}. As such, we   leave   a provable initialization for tensor phase retrieval for future study, which is even more intricate than the unresolved matrix case.
\end{rem} 

\section{End-to-End Algorithms for Concrete Models}\label{sec:end2end}
In this section, we combine the local convergence guarantee in  Section \ref{sec4:examples} and the spectral initialization in Section \ref{sec:speinitial} to yield  theoretical guarantees for the end-to-end algorithms of spectral initialization followed by (N)RGD. Before proceeding, we note the following discrepancies between the local convergence   and the initialization. 
 \begin{rem}
     While the local convergence guarantees in Section \ref{sec4:examples} hold under $n=O(T_{\rm df})$ (up to dependence on $\Lambda_{\min}(\tcalX^\star)$), and Theorems \ref{thm:trelulocal}, \ref{thm:t1bcs} even hold uniformly over a large subset of $T_{\br}^{\bd}$, our initialization in Section \ref{sec:speinitial} applies to a fixed $\tcalX^\star$ only and requires at least $\sqrt{D}$ observations. As such, the forthcoming  results for the end-to-end algorithms are nonuniform, meaning that they are only valid for a fixed unknown $\tcalX^\star$, and require $O(\sqrt{D})$ sample complexity.
 \end{rem} 

\subsection{Tensor Single-Index Models}
Consider the tensor single-index models treated in Section \ref{sec:sim}. We now use Algorithm \ref{alg:ini} to obtain $\tcalX_0\in T^{\bd}_{\br}\cap \mathbb{B}_{\rm F}\Big(\mu\tcalX^\star;\frac{\mu\Lambda_{\min}}{12(\sqrt{3}+1)}\Big)$ required in Theorem \ref{thm:tsim}, yielding the following result. 

\begin{theorem}[Tensor single-index models] \label{thm:e2etsim}
In the setting of Theorem \ref{thm:tsim}, let $\tcalX_0=\tcalZ_0$, where $\tcalZ_0$ is the output of Algorithm \ref{alg:ini}. If  
\(
    n\ge C_{\mu,\sigma}\Big[\frac{T_{\rm df}}{\Lambda_{\min}^2}+\kappa^6 \big((\bar{r}\bar{d})^{3/2}+\bar{r}^4\big)\Big],\)
then (\ref{tsimconvergence}) holds with probability at least $1-C\exp(-c\underline{d})$.   
\end{theorem}
\begin{proof}
Note that Assumption \ref{assump:initial} holds with  $\rho=1$, $\mu>0$ the same as in Assumption \ref{assumptionmu}, and
$\nu:=\|f(g)\|_{\psi_2}\le \|f(g)-\mu g\|_{\psi_2}+\|\mu g\|_{\psi_2}\lesssim \sigma+\mu$. Choose $c_0=[12(\sqrt{3}+1)]^{-1}$ in Corollary \ref{corCpoint1}. Since $\sqrt{D}\le \bar d^{3/2}$, $\bar r\le \bar d$, and $\kappa\ge 1$, the sample size assumed in Theorem \ref{thm:e2etsim} implies \eqref{Cpoint11}. Hence, with probability at least $1-C\exp(-c\underline d)$, we have 
\(
    \|\tcalZ_0-\mu\tcalX^\star\|_{\rm F}\le \frac{\mu\Lambda_{\min}}{12(\sqrt{3}+1)}.
\)
Note that the   sample size in Theorem \ref{thm:e2etsim} also implies the sample size condition in Theorem \ref{thm:tsim}. Thus, on the intersection of the initialization event and the event in Theorem \ref{thm:tsim}, $\tcalX_0=\tcalZ_0$ satisfies the required warm start condition, and \eqref{tsimconvergence} follows. A union bound on these events completes the proof.
\end{proof}

\subsection{Tensor Logistic Regression}
We continue the treatment of tensor logistic regression in Section \ref{sec:tenlogi}, using Algorithm \ref{alg:ini} to fulfill the initialization condition in Theorem \ref{thm:tlogistic}. This leads to the following. 
\begin{theorem}[Tensor logistic regression]\label{thm:e2elogi}
In the setting of Theorem \ref{thm:tlogistic}, if $\tcalX_0= \tcalZ_0/\|\tcalZ_0\|_{\rm F}$, where $\tcalZ_0$ is the output of Algorithm \ref{alg:ini}, and 
\(
    n \ge C\Big[\frac{T_{\rm df}}{\Lambda_{\min}^2}+ \kappa^6 \big((\bar{r}\bar{d})^{3/2}+\bar{r}^4\big)\Big],
\)
then (\ref{tlogiconvergence}) holds with probability at least $1-C\exp(-c\underline{d})$. 
\end{theorem}
\begin{proof}
We write $y_i=f_i(\langle\tcalA_i,\tcalX^\star\rangle)$ with
$f_i(a)=\mathbbm{1}\{u_i\le s(a)\}$ and independent $u_i\sim{\rm Unif}(0,1)$ to see that our setting is encompassed by Assumption \ref{assump:initial} with $\rho=1$,
\(
    \mu=\mathbb{E}[g s(g)]=\mathbb{E}[s'(g)]>0\) (by Stein's identity), \(\nu=\|f_i(g)\|_{\psi_2}=O(1)\) (due to $|y_i|\le 1$). 
Choose $c_0=[48(\sqrt{3}+1)]^{-1}$ in Corollary \ref{corCpoint1}. Since the   sample size in Theorem \ref{thm:e2elogi} implies \eqref{Cpoint11}, 
\(
    \big\|\frac{\tcalZ_0}{\|\tcalZ_0\|_{\rm F}}-\tcalX^\star\big\|_{\rm F} 
    \le \frac{\Lambda_{\min}}{24(\sqrt{3}+1)}\)
holds with probability at least $1-C\exp(-c\underline d)$. It also implies the sample size condition in Theorem \ref{thm:tlogistic}. The conclusion now follows from Theorem \ref{thm:tlogistic} and a union bound on these events.
\end{proof}

\subsection{Tensor ReLU Regression}
Recall that the local convergence for tensor ReLU regression has been established in Theorem \ref{thm:trelulocal}, under the initialization condition of $\tcalX_0\in T^{\bd}_{\br}\cap \mathbb{B}_{\rm F}(\tcalX^\star;c\Lambda_{\min}(\tcalX^\star))$. We seek to invoke Algorithm \ref{alg:ini} to achieve this, yet with $f=\relu$ we have \(\mu = \rho^{-1}\mathbb{E}_{g\sim N(0,1)}[g\relu (\rho g)]= \mathbb{E}_{g\sim N(0,1)}[g\relu(g)] = 1/2\) in Assumption \ref{assump:initial}, and therefore we use $2\tcalZ_0$ (where $\tcalZ_0$ is the output of Algorithm \ref{alg:ini}) as the initialization. 
\begin{theorem}[Tensor ReLU regression] \label{thm:e2erelu}
Under Assumption \ref{gaussianAi}, consider the recovery of $\tcalX^\star\in T^{\bd}_{\br}$ of Tucker rank $(r_1,r_2,r_3)$ from $\tcalA_i$'s and $\{y_i=\relu(\langle\tcalA_i,\tcalX^\star\rangle)\}_{i=1}^n$. Suppose that  \(\tcalX_0=2\tcalZ_0,\)
where $\tcalZ_0$ is the output of Algorithm \ref{alg:ini}, and we run $\tcalX_{t+1}=H^{\bd}_{\br}\left(\tcalX_t-\calP_{T(\tcalX_t)}\big(\tcalH_{\tcalX^\star}(\tcalX_t)\big)\right),~t\ge 0$ with $\tcalH_{\tcalX^\star}(\tcalU)$ defined in (\ref{relutengra}) to obtain $\{\tcalX_t\}_{t\ge 0}$. If
\(n\ge C\kappa^6\big[(\bar{r}\bar{d})^{3/2}+ \bar{r}^4\big],\) 
then (\ref{treluconver}) holds with probability at least $1-C\exp(-c\underline{d})$.
\end{theorem}
\begin{proof}
For $f(a)=\relu(a)$, Assumption \ref{assump:initial} holds with $\mu = 1/2$, $\nu=O(\rho)$. Choose a sufficiently small universal constant $c_0>0$ in Corollary \ref{corCpoint1}. The   sample size in Theorem \ref{thm:e2erelu} implies \eqref{Cpoint11}, and therefore, with the promised probability, (\ref{Cpoint9}) gives
\(\|\tcalX_0-\tcalX^\star\|_{\rm F}\le c \Lambda_{\min}\) for some small enough $c$, thus satisfying   the warm-start condition of Theorem \ref{thm:trelulocal}. The   sample size in Theorem \ref{thm:e2erelu} also implies $n\ge CT_{\rm df}$, so Theorem \ref{thm:trelulocal} applies and yields the local convergence guarantee. 
\end{proof}

\subsection{One-Bit Tensor Sensing}
Similarly, we use Algorithm \ref{alg:ini} to fulfill the initialization condition in Theorem \ref{thm:t1bcs}, yielding the following. 

\begin{theorem}[One-bit tensor sensing]\label{thm:e2e1b}Under Assumption \ref{gaussianAi}, consider the recovery of $\tcalX^\star\in T^{\bd}_{\br}\cap\mathbb{S}_{\rm F}$ of Tucker rank $(r_1,r_2,r_3)$ from $\tcalA_i$'s and $\{y_i=\sign(\langle\tcalA_i,\tcalX^\star\rangle)\}_{i=1}^n$. 
Suppose that $\tcalX_0=\frac{\tcalZ_0}{\|\tcalZ_0\|_{\rm F}}$, where $\tcalZ_0$ is the output of Algorithm \ref{alg:ini}, and we run (\ref{1btsnrgd}) to obtain $\{\tcalX_t\}_{t\ge 0}$. There exist universal constants $C>c>0$ such that, if \(n \ge C \Big[\frac{\log^{3/2}(e/\Lambda_{\min})T_{\rm df}}{\Lambda_{\min}}+\kappa^6 \big((\bar{r}\bar{d})^{3/2}+\bar{r}^4\big)\Big],\)
then (\ref{1btsconvergence}) holds with probability at least $1-C\exp(-c\underline{d})$. 
\end{theorem}
\begin{proof}
In one-bit tensor sensing with $f(a)=\sign(a)$, Assumption \ref{assump:initial} holds with $\rho=1$,  $\mu=\mathbb{E}|g|=\sqrt{2/\pi}$, and $\nu=O(1)$. Choose $c_0=[64(\sqrt{3}+1)]^{-1}$ in Corollary \ref{corCpoint1}. The   sample size in Theorem \ref{thm:e2e1b} implies \eqref{Cpoint11}, and hence
\(
    \big\|\frac{\tcalZ_0}{\|\tcalZ_0\|_{\rm F}}-\tcalX^\star\big\|_{\rm F}
    \le \frac{\Lambda_{\min}}{32(\sqrt{3}+1)}\)
with probability at least $1-C\exp(-c\underline d)$. The first term in the displayed sample size is precisely a sufficient   condition for $\Lambda_{\min}\gtrsim \frac{T_{\rm df}}{n}\log^{3/2}(\frac{n}{T_{\rm df}})$ required by Theorem \ref{thm:t1bcs}. Therefore, the normalized spectral initializer satisfies the warm-start condition in Theorem \ref{thm:t1bcs}, which then guarantees the local convergence with the promised probability. 
\end{proof}

\begin{rem} \label{rem:subopsam} Our Theorems~\ref{thm:e2etsim}--\ref{thm:e2e1b} establish statistically optimal and computationally efficient  guarantees for tensor single-index models, tensor logistic regression, tensor ReLU regression, and one-bit tensor sensing. In the balanced setting $d_1=d_2=d_3=d$ with fixed Tucker rank and bounded condition number, these guarantees use a sample size of order $d^{3/2}$, which is statistically suboptimal. Nonetheless, the same $d^{3/2}$ sample complexity appears in the best-known polynomial-time guarantees for tensor linear regression and tensor completion  \citep{tong2022scaling,xia2019polynomial,han2022optimal,luo2023low}. Moreover, \citet[Theorem~9, Proposition~5]{luo2024tensor} provides evidence against the existence of a polynomial-time estimator achieving constant relative accuracy from $n=O(d^{3/2-\epsilon})$ observations, where $\epsilon$ is any fixed constant in $(0,1)$, in Gaussian linear regression with rank-one symmetric tensor parameter. Since the nonlinear observations only increase the difficulty of tensor estimation, our $d^{3/2}$ sample complexity may reflect a fundamental computational barrier. 
\end{rem}  

\section{Numerical Simulations}\label{sec:simulation}
Numerical experiments are conducted to validate the theoretical findings of the proposed statistical procedures. All experiments were implemented using MATLAB R2022a on a laptop with an Intel CPU up to $2.5$ GHz and $32$ GB RAM. In each simulation, the covariates have i.i.d.\ $N(0,1)$ entries. We consider $d_1=d_2=d_3=d$ and $r_1=r_2=r_3=r$, and generate the true parameter tensor $\tcalX^\star\in T^{\bd}_{\br}\cap \mathbb{S}_{\rm F}$ as
\(
    \tcalX^\star=\tcalS\times_{j=1}^3\bU_j,
    ~
    \tcalS\sim{\rm Unif}(\mathbb{S}_{\rm F}),
    ~
    \bU_j\sim{\rm Unif}(\mathbb{O}_{d_j\times r_j})\) (i.e., Haar distribution over the Stiefel manifold). We report the   Frobenius error 
after $40$ iterations, averaged over $50$ independent trials. The codes are available online (see the \href{https://github.com/junrenchen58/Nonlinear-tensor-regression}{github repository}).

\subsection{End-to-End Recovery under Nonlinear Observations}
We test the end-to-end algorithms for tensor single-index models (Theorem \ref{thm:e2etsim}), tensor ReLU regression (Theorem \ref{thm:e2erelu}), and one-bit tensor sensing (Theorem \ref{thm:e2e1b}).

For tensor single-index models, we specifically test 
\begin{align}
     y_i=T_1(\langle\tcalA_i,\tcalX^\star\rangle),\qquad \textrm{where}~~
    T_1(a)
    :=
    a\mathbbm{1}(|a|\leq1)
    +
    \sign(a)\mathbbm{1}(|a|>1), \label{T1tsim}
\end{align}
 and $\mu=\mathbbm{E}_{g\sim N(0,1)}[gT_1(g)]\approx 0.683$. Since our algorithm for single-index models estimates $\mu\tcalX^\star$ (see Theorem \ref{thm:e2etsim}), we scale the estimate by $1/\mu$ before the error is evaluated.

  For ReLU regression and one-bit sensing, we directly test the algorithms in Theorems \ref{thm:e2erelu} and \ref{thm:e2e1b}.
  Interestingly, our single-index models algorithm also applies. Hence, we shall compare the model-specific algorithms and the generic single-index models method in these two models.

\begin{figure}[ht!]
    \centering
    \includegraphics[width=\textwidth]{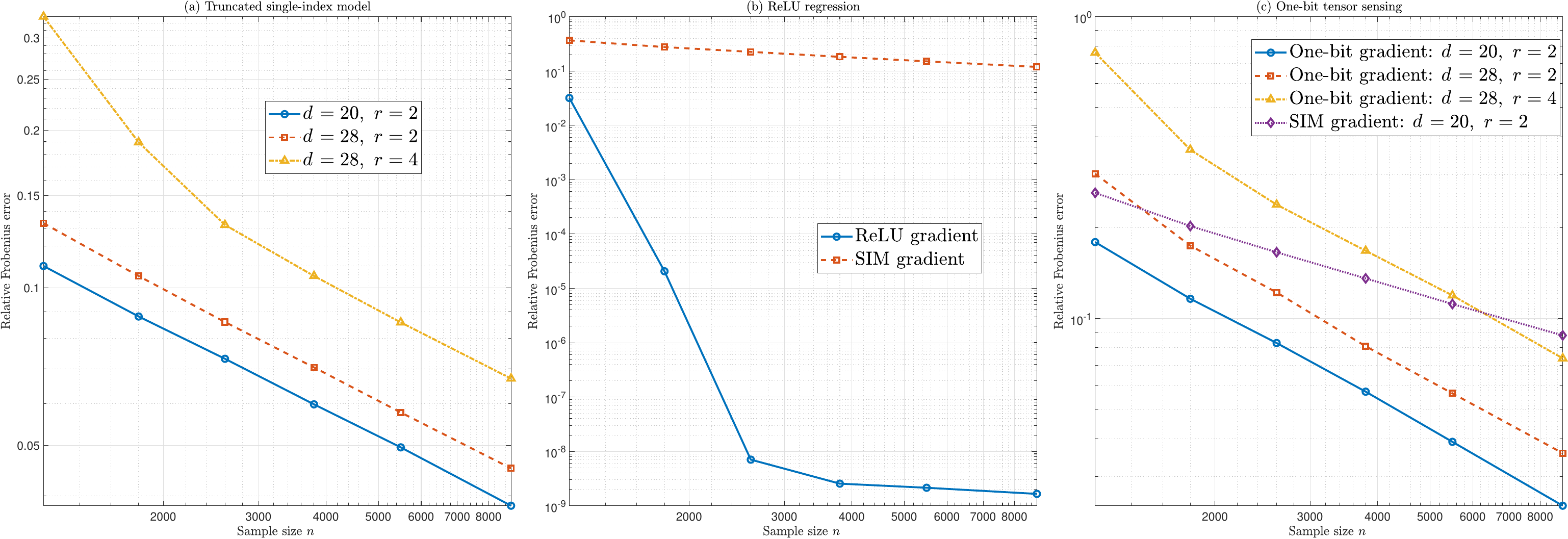}
    \caption{End-to-end estimation procedures for (a) the truncated single-index model with $f_i=T_1$, (b) ReLU regression, and (c) one-bit tensor sensing. The   single-index model curves use   the gradient in (\ref{hsim}).} 
    \label{fig:nonlinear3examples}
\end{figure}

Our results for $(d,r)=(20,2),\,(28,2),\,(28,4)$ under $n=1200,\,1800,\,2600,\,3800,\,5500,\,8800$ are shown in Figure~\ref{fig:nonlinear3examples}.
Figure~\ref{fig:nonlinear3examples}(a) for (\ref{T1tsim}) shows that the estimation error decreases  with the sample size, while increasing the tensor dimension or Tucker rank leads to a larger error, consistently with the dependence on $T_{\mathrm{df}}$.

In ReLU regression, we only test $(d,r)=(20,2)$ since the iterates converge \emph{exactly} to the true parameter in Theorem \ref{thm:e2erelu}. Consistent with our theoretical findings, the model-specific algorithm in Theorem \ref{thm:e2erelu} reaches numerical precision once sufficiently many observations are available, whereas the   single-index models method in Theorem \ref{thm:e2etsim} exhibits a substantially larger statistical error. This confirms the benefit of using a specialized solver under a specific link function. See Figure~\ref{fig:nonlinear3examples}(b).

Similarly, Figure~\ref{fig:nonlinear3examples}(c) shows that the model-specific algorithm in  Theorem \ref{thm:e2e1b} yields a markedly faster error decay than the generic single-index models gradient. This is consistent with the  faster rate $O(n^{-1})$ in Theorem \ref{thm:e2e1b} than the $O(n^{-1/2})$ in Theorem \ref{thm:e2etsim}. Moreover, increasing $d$ or $r$
leads to larger estimation errors of the algorithm in Theorem \ref{thm:e2e1b}, which is consistent with the $\tilde{O}(T_{\rm df}/n)$ theoretical error rate. 
\subsection{Tensor Phase Retrieval with a Hyperspectral Image Dataset}\label{sec:realtpr}

\begin{figure}[ht!]
    \centering
    \subfloat[True tensor $\tcalX^\star$]{%
        \includegraphics[width=0.3\textwidth]{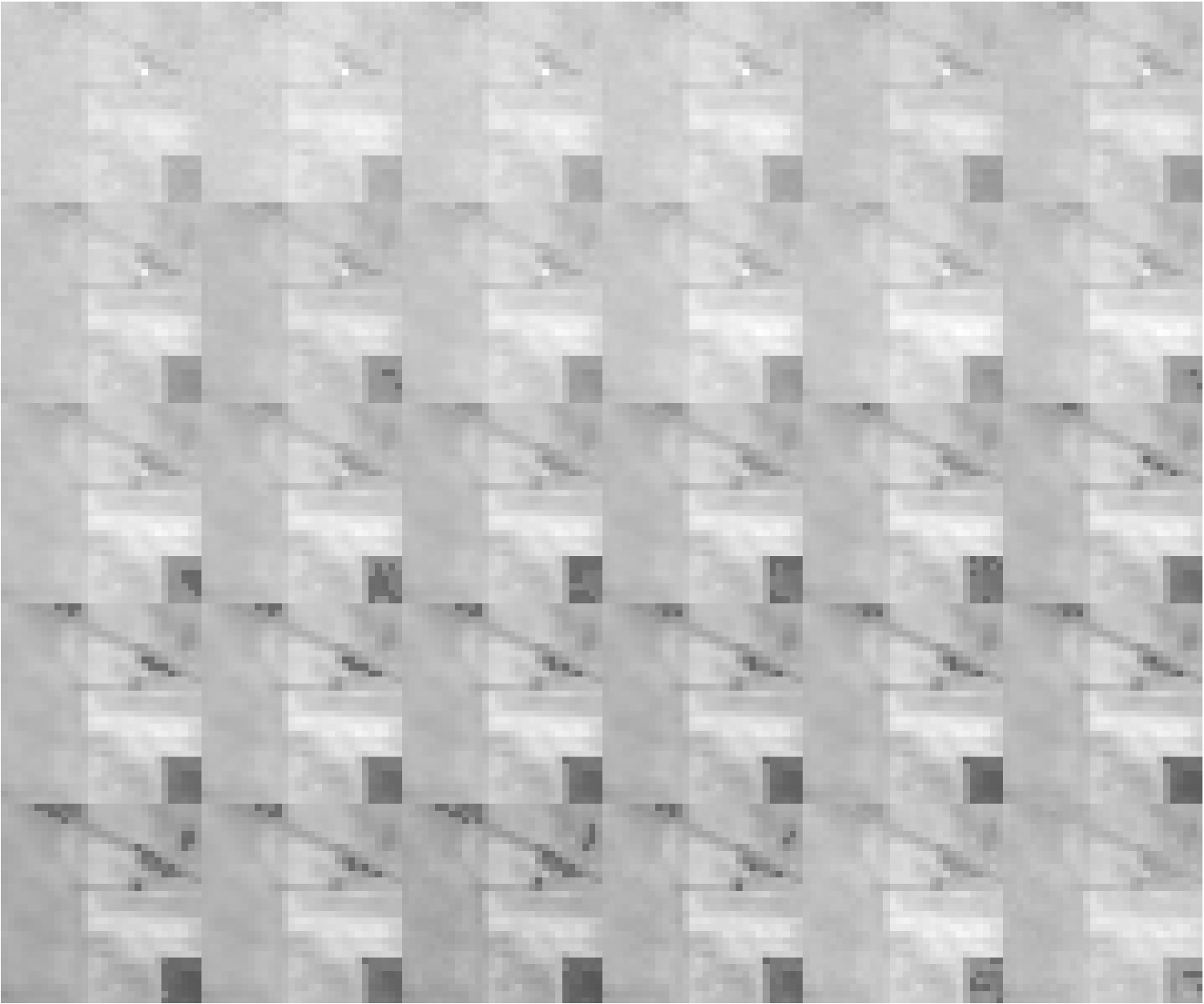}%
        \label{fig:gttt}
    }
    \hspace{0.1\textwidth}
    \subfloat[Initialization $\tcalX_0$]{%
        \includegraphics[width=0.3\textwidth]{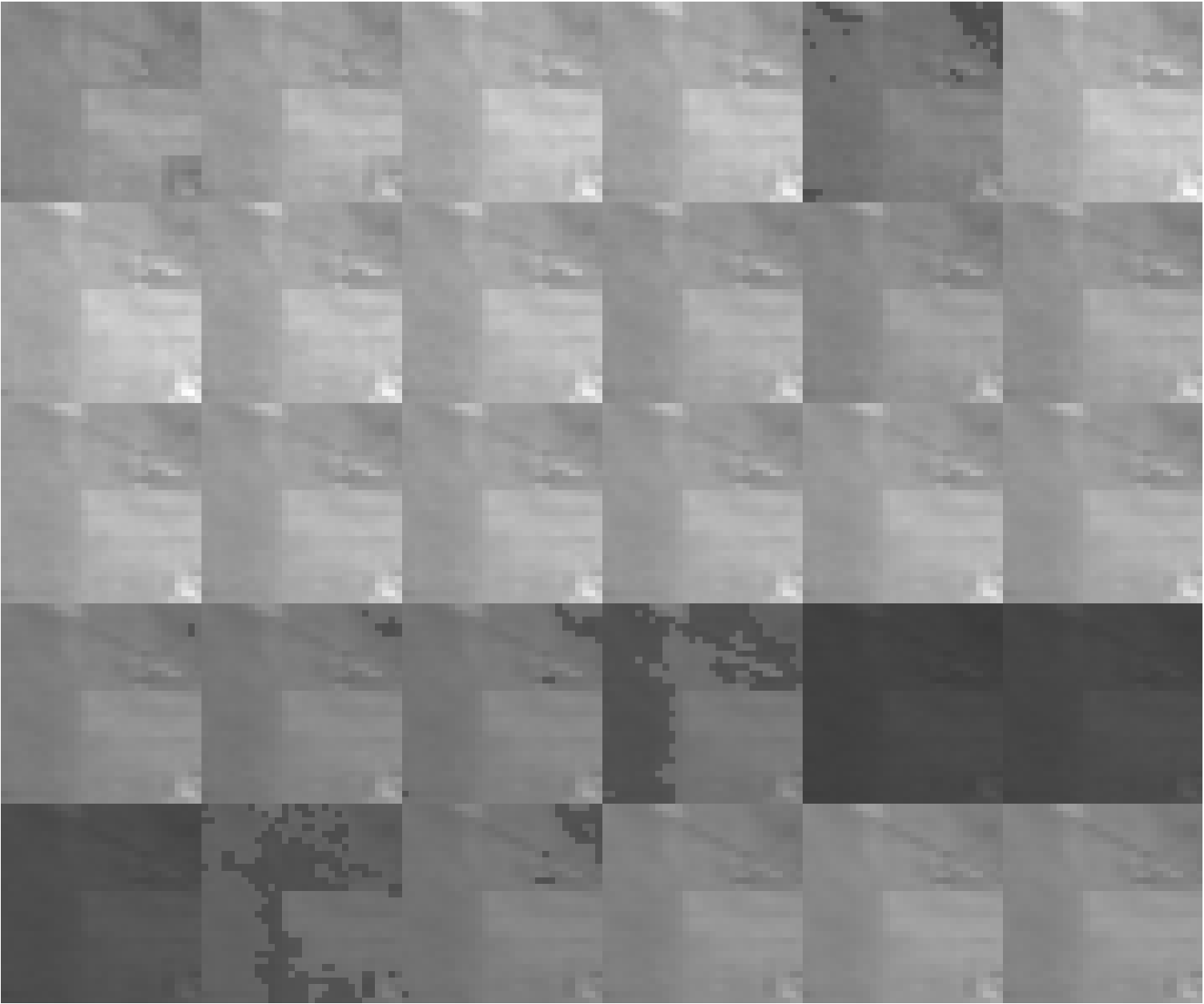}%
        \label{fig:inin}
    }

    \vspace{1ex}

    \subfloat[Tensor method (PSNR $=35.91$)]{%
        \includegraphics[width=0.3\textwidth]{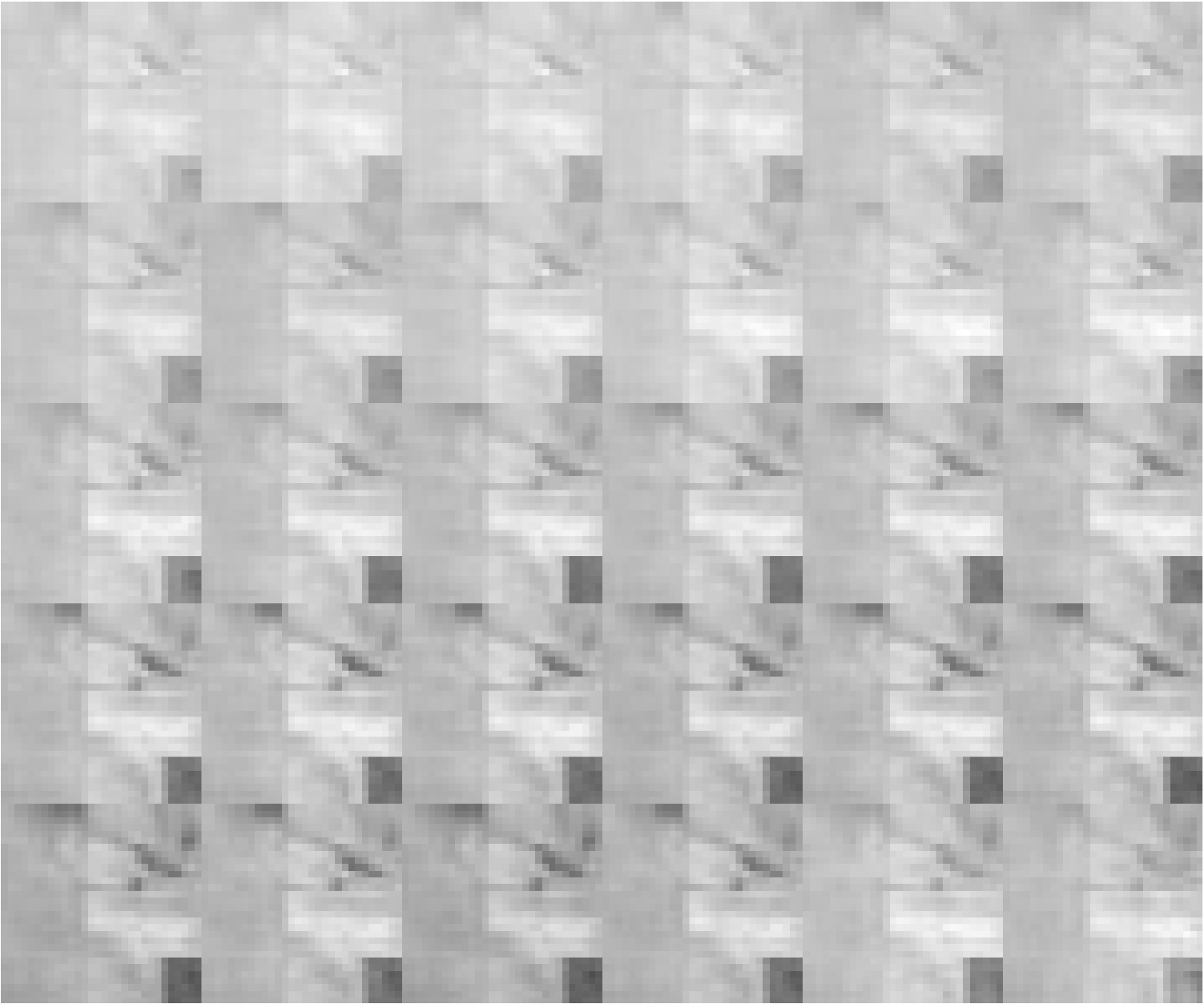}%
        \label{fig:tprpr}
    }\hfill
    \subfloat[Matrix method (PSNR $=20.84$)]{%
        \includegraphics[width=0.3\textwidth]{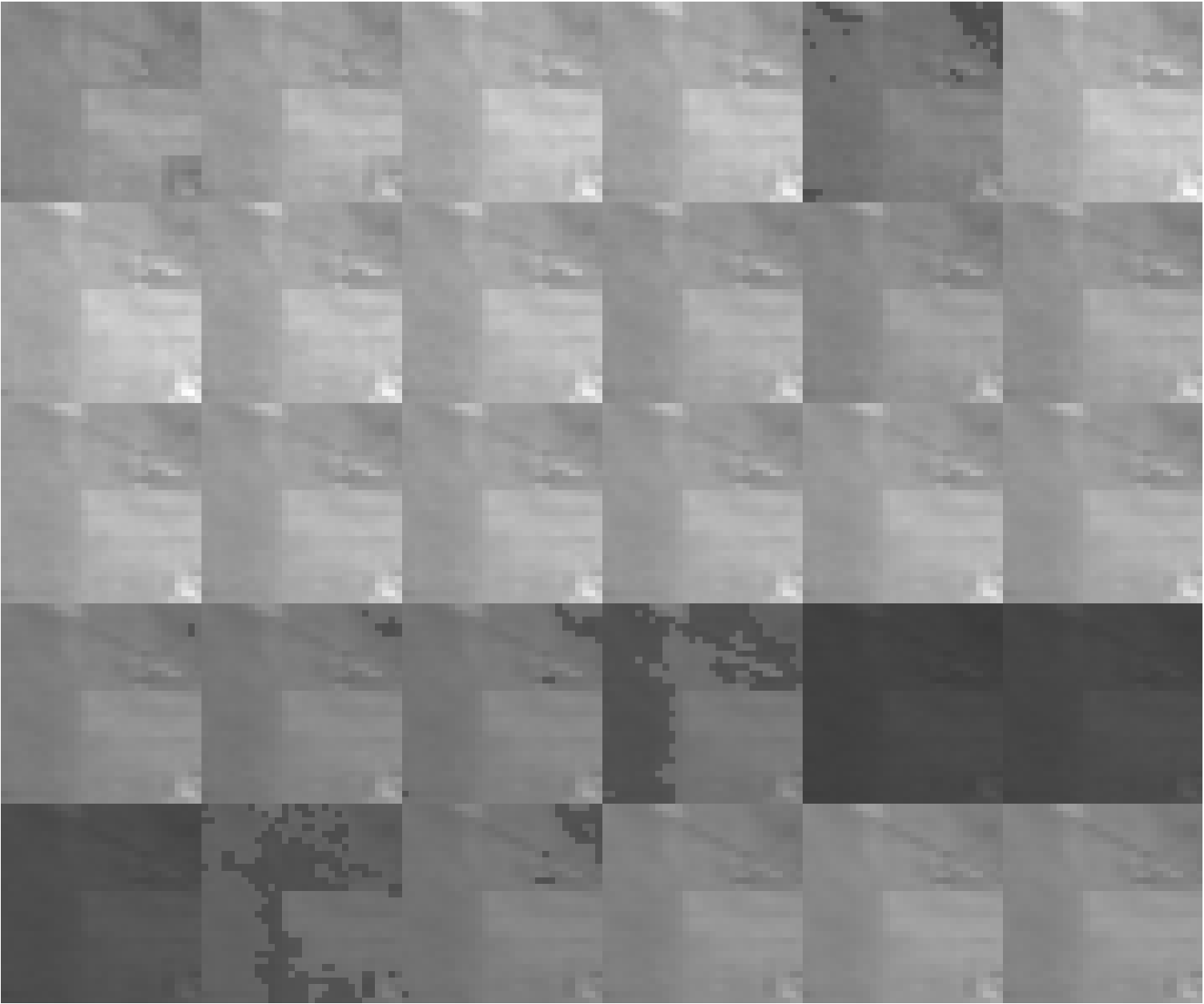}%
        \label{fig:mprpr}
    }\hfill
    \subfloat[Vector method (PSNR $=14.12$)]{%
        \includegraphics[width=0.3\textwidth]{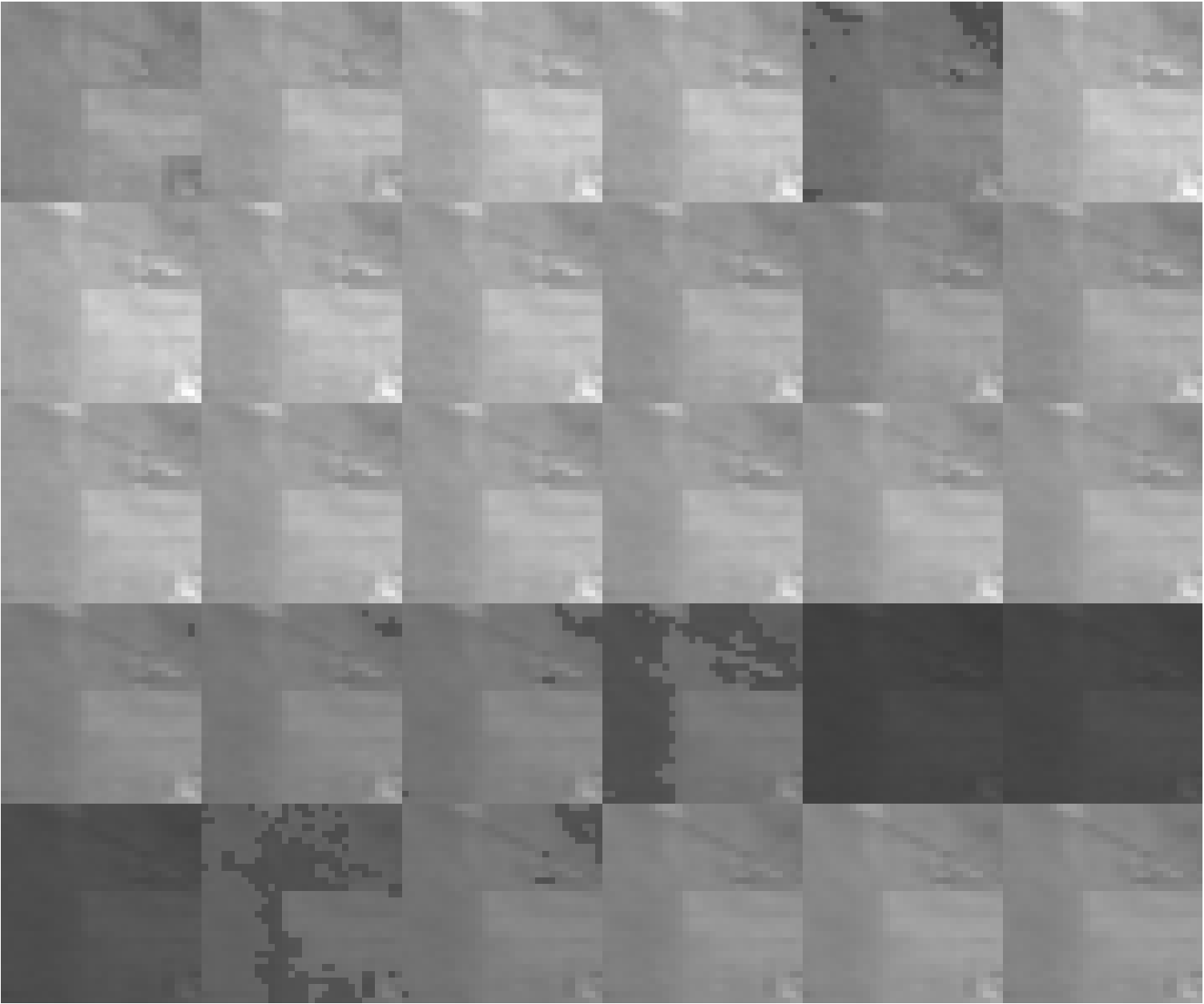}%
        \label{fig:vprpr}
    }

    \caption{Comparison of tensor, matrix, and vector phase retrieval methods on the Indian Pines hyperspectral dataset. The spectral slices are arranged as montages for visualization. The ground-truth tensor in (a) is formed from spectral bands $6$--$35$, and the initialization in (b) is formed from the adjacent bands $36$--$65$ over the same spatial patch. For the tensor method, we test target Tucker ranks $(5,5,5)$ and $(10,10,10)$, and report the better result obtained at rank $(10,10,10)$. For the matrix method, we test target ranks $r=5,10$ after mode-$1$ matricization, and report the better result obtained at rank $r=5$.}
    \label{fig:compare1}
\end{figure}

We further present an experiment based on the Indian Pines hyperspectral dataset. The corrected data cube naturally forms a third-order tensor with two spatial modes and one spectral mode, and is \emph{approximately} low Tucker rank due to the strong correlations among nearby spectral bands. We extract a $30\times30$ spatial patch from rows $61$--$90$ and columns $61$--$90$, and use spectral bands $6$--$35$ to form the underlying tensor $\tcalX^\star\in\mathbb{R}^{30\times30\times30}$ (visualized in Figure~\ref{fig:compare1}(a)). We normalize $\tcalX^\star$ to unit Frobenius norm and generate
\(
    n=6750=30^3/4
\)
Gaussian phaseless measurements $y_i=|\langle\tcalA_i,\tcalX^\star\rangle|$. Since a provably accurate initialization for tensor phase retrieval remains unavailable, we use a data-driven initialization: $\tcalX_0$ is constructed from the adjacent spectral bands $36$--$65$ over the same spatial region and is also normalized to unit Frobenius norm (visualized in Figure~\ref{fig:compare1}(b)).

Starting from the same initialization $\tcalX_0$, we compare three methods, each run for $50$ iterations. The tensor method is the RGD iteration in (\ref{TPRRGDiterate}) with target Tucker rank $(r,r,r)$, where $r=5,10$. The matrix method matricizes the tensor along the first mode and applies low-rank matrix phase retrieval via projected gradient descent \citep{soltanolkotabi2019structured} to $\calM_1(\tcalX^\star)\in\mathbb{R}^{30\times900}$ with target rank $r=5,10$. The vector method runs the reshaped Wirtinger flow \citep{zhang2017nonconvex} with fixed step size one and does not exploit any low-dimensional structure. The results are displayed in Figure~\ref{fig:compare1}. The tensor method with rank $(10,10,10)$ achieves $\|\widehat{\tcalX}_{\rm TPR}-\tcalX^\star\|_{\rm F}=0.021$ and PSNR $=35.91$, whereas the best matrix method, attained at rank $5$, gives $\|\widehat{\tcalX}_{\rm MPR}-\tcalX^\star\|_{\rm F}=0.119$ and PSNR $=20.84$. The vector method gives $\|\widehat{\tcalX}_{\rm VPR}-\tcalX^\star\|_{\rm F}=0.258$ and PSNR $=14.12$. This comparison is consistent with the tensor method's ability to jointly exploit the low-rank structure across all three tensor modes.

\section{Concluding Remarks}\label{sec:conclusion}

This paper develops a unified framework for low Tucker rank tensor estimation from nonlinear observations. Its central ingredient is the tensor restricted approximate invertibility condition (T-RAIC), which quantifies how well a model-specific empirical gradient map approximates the ideal  descent step under the tensor dual norm. We show that T-RAIC yields local linear convergence guarantees for RGD and its normalized counterpart. Under Gaussian covariates, we establish T-RAICs for single-index models, logistic regression, phase retrieval, ReLU regression, and one-bit tensor sensing. These results imply exact local recovery for phase retrieval and ReLU regression, and near-optimal statistical errors for the remaining models. For all models except phase retrieval, we further show that a simple spectral initialization followed by RGD yields computationally efficient end-to-end algorithms.

Several questions remain open. First, it would be important to develop a provable initialization for low-rank tensor phase retrieval, which is open even for the simpler problem of phase retrieval of low-rank matrices (Remark \ref{rem:mpr}). Second, extending the framework to   multi-response tensor regression would accommodate substantially richer models. Third, it is of interest to develop robust variants of the framework under different settings, such as
low-rank-plus-sparse models, heavy-tailed covariates, or adversarially corrupted observations. We believe that some robust counterpart of the T-RAIC would be useful for this investigation.

\bibliography{libr}

  \appendix

\section{Proof of Convergence Guarantees (Theorems \ref{thm:prgdconver}, \ref{thm:nrgdconver})} 
\subsection{Proof of Theorem \ref{thm:prgdconver} (Convergence of Algorithm \ref{alg:prgd})}\label{app:provergdconverge}
\begin{proof} 
By definition, the T-RAIC reads  
\begin{align}\label{traicthm2}
       \big\|\tcalU-\tcalX^\star-\eta\cdot\tcalH (\tcalU)\big\|_{(T^{\bd}_{2\br}\cap \mathbb{B}_{\rm F}) ^\circ}\le\frac{\|\tcalU-\tcalX^\star\|_{\rm F}}{4(\sqrt{3}+1)}+\mu_2,\qquad\forall \tcalU\in\scrU. 
   \end{align}
We define a sequence $\{f_t\}_{t=0}^\infty$ by $f_0= \min\big\{R_{\rm loc},\frac{\Lambda_{\min}}{12(\sqrt{3}+1)}\big\}$ and the recurrence 
\(
    f_{t+1} = \frac{f_t}{2}+(\sqrt{3}+1)\mu_2,~ \forall t\ge 0.\) 
Solving the recurrence gives 
\begin{align}\label{ftestimate31}
    f_t= 2^{-t}\min\bigg\{R_{\rm loc}, \frac{\Lambda_{\min}}{12(\sqrt{3}+1)}\bigg\} + 2(\sqrt{3}+1)\mu_2(1-2^{-t}),\qquad \forall t\ge 0.
\end{align}
In light of $\mu_2\le c\min\{R_{\rm loc},\Lambda_{\min}\}$ for small enough $c$, we have $f_0>2(\sqrt{3}+1)\mu_2$, and hence $\{f_t\}_{t=0}^\infty$ is monotonically decreasing, which implies $f_t\le f_0= \min\big\{R_{\rm loc},\frac{\Lambda_{\min}}{12(\sqrt{3}+1)}\big\}$ for any $t\ge 0$.

We now seek to prove $\|\tcalX_t-\tcalX^\star\|_{\rm F}\le f_t$ for $t\ge 0$. We now use induction to show $$\|\mathbfcal{X}_t-\mathbfcal{X}^\star\|_{\rm F}\le f_t,\quad\text{and}\quad \tcalX_t\in\scrU,\quad\text{and}\quad \tucker(\tcalX_t)=(r_1,r_2,r_3) $$ 
for $t\ge 0$.

    We first observe that this holds   for $t=0$: by $f_0=\min\{R_{\rm loc},\frac{\Lambda_{\min}}{12(\sqrt{3}+1)}\}$, $\tcalX_0\in T^{\bd}_{\br}\cap \mathbb{B}_{\rm F}\big(\tcalX^\star;\min\big\{R_{\rm loc},\frac{\Lambda_{\min}}{12(\sqrt{3}+1)}\big\}\big)$ and $\scrU\supset T^{\bd}_{\br}\cap \mathbb{B}_{\rm F}(\tcalX^\star;R_{\rm loc})$, we have $\|\tcalX_0-\tcalX^\star\|_{\rm F}\le f_0$ and $\tcalX_0\in \scrU$;  moreover, Weyl's inequality (see, e.g., \citealp[Theorem 1]{stewart1998perturbation}), together with $\|\tcalX_0-\tcalX^\star\|_{\rm F}\le \frac{1}{2}\Lambda_{\min}(\tcalX^\star)$ and $\tucker(\tcalX^\star)=(r_1,r_2,r_3)$, ensures that $\lambda_{r_j}(\calM_{j}(\tcalX_0))>0$ for $j\in[3]$, where we suppose that $\lambda_{1}(\calM_{j}(\tcalX_0))\ge \lambda_{2}(\calM_{j}(\tcalX_0))\ge \cdots\ge \lambda_{d_j}(\calM_{j}(\tcalX_0))$ are the singular values of $\calM_{j}(\tcalX_0)$; this yields $\tucker(\tcalX_0)\succeq(r_1,r_2,r_3)$, and by combining with $\tcalX_0\in T^{\bd}_{\br}$, we obtain $\tucker(\tcalX_0)=(r_1,r_2,r_3)$.

    To complete the induction, it remains to show \[\|\mathbfcal{X}_{t+1}-\mathbfcal{X}^\star\|_{\rm F}\le f_{t+1},\quad\text{and}\quad \tcalX_{t+1}\in \scrU,\quad\text{and}\quad\tucker(\tcalX_{t+1})=(r_1,r_2,r_3)\] under the hypothesis of \[\|\mathbfcal{X}_t-\mathbfcal{X}^\star\|_{\rm F}\le f_t,\quad\text{and}\quad \tcalX_t\in\scrU,\quad\text{and}\quad \tucker(\tcalX_t
    )=(r_1,r_2,r_3).\] To this end, we first use Lemma \ref{pertur} to obtain  
    \begin{align} \label{B.37}
        \|\tcalX_{t+1}-\tcalX^\star\|_{\rm F}=\big\|H^{\bd}_{\br}\big(\tcalX_t-\eta\cdot \tcalH^{\rm (r)} (\tcalX_t)\big)-\mathbfcal{X}^\star\big\|_{\rm F}
        \le (\sqrt{3}+1)\big\|\mathbfcal{X}_t - \mathbfcal{X}^\star - \eta\cdot \mathbfcal{H}^{\rm (r)}(\mathbfcal{X}_t)\big\|_{\rm F}.
    \end{align}
    Using $\mathbfcal{H}^{\rm (r)}(\mathbfcal{X}_t) = \calP_{T(\mathbfcal{X}_t)}\big(\mathbfcal{H}(\mathbfcal{X}_t) \big)$ and  the   bound $\|\mathbfcal{X}^\star-\calP_{T(\mathbfcal{X}_t)}(\mathbfcal{X}^\star)\|_{\rm F}\le \frac{3\|\mathbfcal{X}_t-\mathbfcal{X}^\star\|^2_{F}}{\Lambda_{\min}}$ from Lemma \ref{rie} (note that this bound relies on $\tucker(\tcalX_t)=(r_1,r_2,r_3)$), along with $\|\tcalX_t-\tcalX^\star\|_{\rm F} \le f_t \le \frac{\Lambda_{\min}}{12(\sqrt{3}+1)}$, we continue from (\ref{B.37}) to obtain
    \begin{align}\nn
         &\big\|\mathbfcal{X}_{t+1}-\mathbfcal{X}^\star\big\|_{\rm F}\\&\nn\le (\sqrt{3}+1)\|\tcalX^\star-\calP_{T(\tcalX_t)}(\tcalX^\star)\|_{\rm F} + (\sqrt{3}+1)\big\|\mathbfcal{X}_t-\calP_{T(\mathbfcal{X}_t)}(\mathbfcal{X}^\star)-\eta\cdot\calP_{T(\mathbfcal{X}_t)}\big(\mathbfcal{H}(\mathbfcal{X}_t)\big)\big\|_{\rm F}\\
         &\le \frac{3(\sqrt{3}+1)\|\mathbfcal{X}_t-\mathbfcal{X}^\star\|^2_{F}}{\Lambda_{\min}} + (\sqrt{3}+1)\big\|\mathbfcal{X}_t-\calP_{T(\mathbfcal{X}_t)}(\mathbfcal{X}^\star)-\eta\cdot\calP_{T(\mathbfcal{X}_t)}\big(\mathbfcal{H}(\mathbfcal{X}_t)\big)\big\|_{\rm F}\nn\\
         &\le \frac{f_t}{4} +(\sqrt{3}+1)\big\|\mathbfcal{X}_t-\calP_{T(\mathbfcal{X}_t)}(\mathbfcal{X}^\star)-\eta\cdot\calP_{T(\mathbfcal{X}_t)}\big(\mathbfcal{H}(\mathbfcal{X}_t)\big)\big\|_{\rm F}.\label{2ov333}
    \end{align}
  We then further control the second term above: 
   \begin{align}\nn
       &\big\|\mathbfcal{X}_t-\calP_{T(\mathbfcal{X}_t)}(\mathbfcal{X}^\star)-\eta\cdot\calP_{T(\mathbfcal{X}_t)}\big(\mathbfcal{H}(\mathbfcal{X}_t)\big)\big\|_{\rm F} \\\nn
       &\stackrel{(a)}{=} \big\|\calP_{T(\mathbfcal{X}_t)}\big(\mathbfcal{X}_t-\mathbfcal{X}^\star-\eta\cdot \mathbfcal{H}(\mathbfcal{X}_t)\big)\big\|_{\rm F}\\ \nn
       &\stackrel{(b)}{=} \sup_{\mathbfcal{U}\in T(\mathbfcal{X}_t)\cap \mathbb{B}_{\rm F} }\big\langle \mathbfcal{U}, \mathbfcal{X}_t-\mathbfcal{X}^\star-\eta\cdot \mathbfcal{H}(\mathbfcal{X}_t)\big\rangle
       \\  \nn
       &\stackrel{(c)}{\le} \sup_{\mathbfcal{U}\in T^{\bd}_{2\br}\cap \mathbb{B}_{\rm F}}\big\langle \mathbfcal{U}, \mathbfcal{X}_t-\mathbfcal{X}^\star-\eta\cdot \mathbfcal{H}(\mathbfcal{X}_t)  
\big\rangle  \\
\nn
       &\stackrel{(d)}{=} \big\|\mathbfcal{X}_t-\mathbfcal{X}^\star-\eta\cdot \mathbfcal{H}(\mathbfcal{X}_t) \big\|_{(T^{\bd}_{2\br}\cap \mathbb{B}_{\rm F})^\circ}\\\label{B.25}
       &\stackrel{(e)}{\le} \frac{\|\tcalX_t-\tcalX^\star\|_{\rm F}}{4(\sqrt{3}+1)}+\mu_2\le \frac{f_t}{4(\sqrt{3}+1)}+\mu_2,  
\end{align}
where in $(a)$ we notice that $T(\mathbfcal{X}_t)$ is a linear subspace containing $\tcalX_t$, $(b)$ expresses the norm of the projection in a variational form,  $(c)$ holds due to $T(\tcalX_t)\subset T^{\bd}_{2\br}$, $(d)$ is due to the definition of the dual norm, and in $(e)$ we note that  $\tcalX_t \in \scrU$ allows us to apply the RAIC in (\ref{traicthm2}). Substituting (\ref{B.25}) into (\ref{2ov333}), along with $f_{t+1}=\frac{f_t}{2}+(\sqrt{3}+1)\mu_2$, yields 
\begin{align}\label{B.26}
    \|\tcalX_{t+1}-\tcalX^\star\|_{\rm F}\le \frac{f_t}{2}+ (\sqrt{3}+1)\mu_2 = f_{t+1}.
\end{align}
This also implies $\|\tcalX_{t+1}-\tcalX^\star\|_{\rm F}\le f_{t+1}\le f_0\le R_{\rm loc}$. Combining with $\tcalX_{t+1}\in T^{\bd}_{\br}$, we have $\tcalX_{t+1}\in T^{\bd}_{\br}\cap \mathbb{B}_{\rm F}(\tcalX^\star;R_{\rm loc})$, and hence $\tcalX_{t+1}\in \tcalU$. It remains to show $\tucker(\tcalX_{t+1})=(r_1,r_2,r_3)$, and in fact, this follows from $\|\tcalX_{t+1}-\tcalX^\star\|_{\rm F}\le f_{t+1}\le f_0\le \frac{\Lambda_{\min}}{2}$ and Weyl's inequality; we omit these details. As such, we have completed the induction and shown $\|\tcalX_{t}-\tcalX^\star\|_{\rm F}\le f_t$ for all $t\ge 0$. In light of (\ref{ftestimate31}), the proof is complete. 
\end{proof}

\subsection{Proof of Theorem \ref{thm:nrgdconver} (Convergence of Algorithm \ref{alg:nprgd})}\label{app:proof32}
\begin{proof}
    The proof is largely parallel to that of Theorem \ref{thm:prgdconver}, so we will omit some of the details. In fact, the only additional ingredient is a standard error bound for analyzing normalization in Euclidean norm (see, e.g., \citealp[Equation (2.12)]{chen2024robust}): for any $\tcalX_t\ne 0$ and $\|\tcalX^\star\|_{\rm F}=\beta>0$, 
    \begin{align}\label{normalbound}
        \bigg\| \beta\frac{\tcalX_t}{\|\tcalX_t\|_{\rm F}}-\tcalX^\star\bigg\|_{\rm F} \le 2\|\tcalX_t-\tcalX^\star\|_{\rm F}. 
    \end{align}

    We now start the proof. We define a sequence $\{f_t\}_{t=0}^\infty$ by $f_0=\min\big\{R_{\rm loc},\frac{\Lambda_{\min}}{24(\sqrt{3}+1)}\big\}$ and the recurrence 
    \(f_{t+1}=\frac{1}{2}f_t+2(\sqrt{3}+1)\mu_2,~t\ge0.\)
    Solving the recurrence finds that
    \begin{align}
        f_t =  \frac{1}{2^t}\min\Big\{R_{\rm loc},\frac{\Lambda_{\min}}{24(\sqrt{3}+1)}\Big\}+\Big(1-\frac{1}{2^t}\Big)4(\sqrt{3}+1)\mu_2,\quad t\ge 0.\label{linerecurrence}
    \end{align}
    Moreover, under $\mu_2\le c\min\{R_{\rm loc},\Lambda_{\min}\}$ for small enough $c$, we have $4(\sqrt{3}+1)\mu_2<f_0$, and hence $f_t$ is monotonically decreasing. This implies $f_{t+1}<f_t$ for $t\ge 0$, which also implies $f_{t}\le f_0 = \min\big\{R_{\rm loc},\frac{\Lambda_{\min}}{24(\sqrt{3}+1)}\big\}$ for all $t\ge 0$.

    We now seek to prove $\|\tcalX_t-\tcalX^\star\|_{\rm F}\le f_t$ for all $t\ge 0$. To this end, we use induction to prove that $\|\tcalX_t-\tcalX^\star\|_{\rm F}\le f_t$, $\tcalX_t\in \scrU$, and $\tucker(\tcalX_t)=(r_1,r_2,r_3)$ hold for all $t\ge 0.$

    It is easy to see that these hold for $t=0$ by $\tcalX_0\in T^{\bd}_{\br}\cap \mathbb{B}_{\rm F}(\tcalX^\star;\min\{R_{\rm loc},\frac{\Lambda_{\min}}{24(\sqrt{3}+1)}\})\cap \beta \mathbb{S}_{\rm F}\subset \scrU$ and Weyl's inequality. Therefore, it remains to establish
    \[\|\mathbfcal{X}_{t+1}-\mathbfcal{X}^\star\|_{\rm F}\le f_{t+1},\quad\text{and}\quad \tcalX_{t+1}\in \scrU,\quad\text{and}\quad\tucker(\tcalX_{t+1})=(r_1,r_2,r_3)\] under the hypothesis of \[\|\mathbfcal{X}_t-\mathbfcal{X}^\star\|_{\rm F}\le f_t,\quad\text{and}\quad \tcalX_t\in\scrU,\quad\text{and}\quad \tucker(\tcalX_t
    )=(r_1,r_2,r_3).\]
    By revisiting the arguments in (\ref{B.37})--(\ref{B.26}), along with straightforward adaptations, we are able to show that
    \[\big\|H^{\bd}_{\br}\big(\tcalX_t-\eta \cdot \tcalH^{\rm (r)}(\tcalX_t)\big)-\tcalX^\star\big\|_{\rm F}\le \frac{f_t}{4}+(\sqrt{3}+1)\mu_2.\]
    Further using (\ref{normalbound}), the update rule of Algorithm \ref{alg:nprgd}, and $f_{t+1}=\frac{1}{2}f_t+2(\sqrt{3}+1)\mu_2$,
    \[\|\tcalX_{t+1}-\tcalX^\star\|_{\rm F} \le \frac{f_t}{2}+2(\sqrt{3}+1)\mu_2 = f_{t+1}.\]
    This implies $\|\tcalX_{t+1}-\tcalX^\star\|_{\rm F} \le f_0\le R_{\rm loc}$. Combining with $\scrU\supset T^{\bd}_{\br}\cap \mathbb{B}_{\rm F}(\tcalX^\star;R_{\rm loc})\cap \beta\mathbb{S}_{\rm F}$ and $\tcalX_{t+1}\in T^{\bd}_{\br}\cap\beta\mathbb{S}_{\rm F}$ ensured by the algorithm, we reach $\tcalX_{t+1}\in\scrU$. Also, combining $\|\tcalX_{t+1}-\tcalX^\star\|_{\rm F}\le \frac{1}{2}\Lambda_{\min}$ and Weyl's inequality, it follows that $\tucker(\tcalX_{t+1})=(r_1,r_2,r_3)$. Therefore, we have completed the induction and shown $\|\tcalX_{t}-\tcalX^\star\|_{\rm F}\le f_t$ for all $t\ge 0$. The proof is complete in view of (\ref{linerecurrence}).
\end{proof}

\section{Proofs of the RAICs (Theorems \ref{pro:traicsim}, \ref{traicglm}, \ref{traicpr}, \ref{thm:traicrelu}, \ref{traic1bcs})} \label{sec:generalraic}
 For the convenience of future research, we establish RAICs for parameter  in $\mathbb{R}^d$ with a general structure (covering sparsity and low-rankness). They directly imply the T-RAICs in Theorems \ref{pro:traicsim}, \ref{traicglm}, \ref{traicpr}, \ref{thm:traicrelu}, and \ref{traic1bcs}. For this purpose, we need to define a general RAIC and introduce the complexity metrics for a general set in Section \ref{sec:prelimimore}. 

 \subsection{General RAIC, Gaussian Width, and Covering Number}\label{sec:prelimimore}
Let $\bh:\mathbb{R}^d\to \mathbb{R}^d$ be a gradient map, then the RAIC requires that $\bh(\bu)$ well approximates the ideal descent step $\bu-\bx$ under the dual norm with respect to some set $\calK\subset \mathbb{B}_2^d$, uniformly over some constraint set $\calU\subset\mathbb{R}^d$. Here, as with \citet{friedlander2021nbiht}, the dual norm of $\bu$ with respect to a symmetric set $\calW\subset \mathbb{R}^d$ is defined as \(\|\bu\|_{\calW^\circ}=\sup_{\bw\in\calW}\langle\bw,\bu\rangle.\) 

\begin{definition}[General RAIC]
\label{def:raic}
 We say that the gradient map $\bh:\mathbb{R}^d\to \mathbb{R}^d$ satisfies    RAIC with respect to a symmetric set $\calK\subset \mathbb{B}_2^d$ over a constraint set $\calU\subset \mathbb{R}^d$ for a target parameter $\bx$ up to an error function $R:\mathbb{R}^d\to \mathbb{R}$, if it holds that 
 \[
        \big\|\bu-\bx- \bh(\bu)\big\|_{\calK^\circ} \le R(\bu),\qquad \forall \bu\in\calU.
 \] 
    We denote this by \(\bh(\bu)\sim {\rm RAIC}\big(\calK;\calU,\bx,R(\bu)\big).\)
\end{definition}

\begin{rem}
    By identifying $\mathbb{R}^{d_1\times d_2\times d_3}$ with $\mathbb{R}^{d_1d_2d_3}$ via vectorization, the general RAIC in Definition \ref{def:raic} with $\calK = T^{\bd}_{2\br}\cap \mathbb{B}_{\rm F}$, $\calU=\scrU$ recovers the T-RAIC in Definition \ref{def:traic}. 
\end{rem}


We will establish the general RAICs under a number of observations proportional to the complexities of $\calK$ and $\calU$, captured by the Gaussian width and covering number to be introduced shortly. Specifically,   for   $\calW$ in $\mathbb{R}^d$:
\begin{itemize}
    \item The Gaussian width defined as 
\[\omega(\calW) = \mathbb{E}_{\bg\sim N(0,\bI_d)}\sup_{\bw\in\calW}\,\langle\bg,\bw\rangle;\]
\item The radius-$\epsilon$  covering number defined as 
\[\calN(\calW;\epsilon)=\inf\bigg\{|\calS|:\calS\subset \calW,~\calW\subset \bigcup_{\bw\in\calS}\mathbb{B}_2^d(\bw;\epsilon)\bigg\}.\]
\end{itemize}  

\begin{rem}
When $\calK,\calU$ are low-complexity sets (such as the set  of sparse vectors), RAIC holds in high-dimensional settings with $n\ll d$. The main components in Definition \ref{def:raic} that allow for $n\ll d$ are the dual norm $\|\cdot\|_{\calK^\circ}$ and the constraint $\bu\in\calU$.      
\end{rem}

  We have the following bounds on the Gaussian width and covering number of $T^{\bd}_{\br}\cap \mathbb{B}_{\rm F}$. 
\begin{pro}[See, e.g., \citealp{rauhut2017low}] \label{pro:tensorbound}For some universal constant $C$ and any $\epsilon\in (0,1)$,
\[\calN(T^{\bd}_{\br}\cap\mathbb{B}_{\rm F};\epsilon)\le(C\epsilon^{-1})^{T_{\rm df}},\qquad\omega(T^{\bd}_{\br}\cap\mathbb{B}_{\rm F})\le C\sqrt{T_{\rm df}}.\] 
\end{pro}

In the remainder of this appendix, we establish general RAICs under Gaussian design.
\begin{assumption}\label{gaussianaiassump}
    $\ba_1,\cdots,\ba_n$ are i.i.d. $N(0,\bI_d)$ vectors. 
\end{assumption}

 \subsection{Single-Index Models (Proof of Theorem \ref{pro:traicsim})}\label{app:raicforsim}
 
 Note that the vector counterpart of (\ref{hsim}) is given by 
 \begin{align}\label{husimvec}
     \bh(\bu)=\frac{1}{n}\sum_{i=1}^n(\ba_i^\top\bu -f_i(\ba_i^\top\bx))\ba_i. 
 \end{align} 
It satisfies the following   RAIC. Let $\calC$ be a cone in $\mathbb{R}^d$, then we let 
\[\calC_{(1)}:=(\calC-\calC)\cap \mathbb{B}_2^d.\]

\begin{theorem}[RAIC for Single-Index Models]\label{thm:raicsim} Let $\calC$ be a cone in $\mathbb{R}^d$,   fix  $\bx\in \calC\cap \mathbb{S}^{d-1}$, and let $\calK$ be a symmetric set obeying $\calK\subset \mathbb{B}_2^d$. Under Assumption \ref{gaussianaiassump}, 
  assume that  $f_1,\cdots,f_n$ are independent of $\ba_i$'s and are i.i.d. copies of a (possibly) random function $f$  satisfying 
Assumption \ref{assumptionmu}. There exist some universal constants $C>c>0$ such that the following holds. For any $\delta\in(0,1)$, 
if  $n\ge C[\omega^2(\calK)+\omega^2(\calC_{(1)})]$, then with probability at least $1-C\exp(-c\delta^2 n)$, $\bh(\bu)$ in (\ref{husimvec}) satisfies
\[\bh(\bu)\sim {\rm RAIC}\bigg(\calK;\calC,\mu\bx,\frac{\|\bu-\mu\bx\|_{2}}{20}+ C\sigma \Big(\frac{\omega(\calK)}{\sqrt{n}}+\delta\Big)\bigg).\]  
\end{theorem}
\begin{proof}
By definition, we seek to establish  
\begin{align}
    \big\| \bu -\mu\bx - \bh(\bu)\big\|_{\calK^\circ}\le \frac{\|\bu-\mu\bx\|_{2}}{20}+C\sigma\bigg(\frac{\omega(\calK)}{\sqrt{n}}+\delta\bigg),\qquad \forall \bu\in\calC. \label{simraic}
\end{align} 
By substituting (\ref{husimvec}) and triangle inequality, 
    \begin{align} 
        \big\| \bu -\mu\bx - \bh(\bu)\big\|_{\calK^\circ}\label{tristartsim}
        &\le \underbrace{\Big\|\Big(\frac{1}{n}\sum_{i=1}^n\ba_i\ba_i^\top-\bI_d\Big)(\bu-\mu\bx)\Big\| _{\calK^\circ}}_{:=\Xi_1}+ \underbrace{\Big\|\frac{1}{n}\sum_{i=1}^n\big( f_i(\ba_i^\top\bx)-\mu\ba_i^\top\bx\big)\ba_i\Big\|_{\calK^\circ}}_{:=\Xi_2}.
    \end{align}

   \paragraph{Bounding $\Xi_1$.} We can only consider $\bu\ne\mu\bx$ (we are done otherwise). Since $\bu,\mu\bx\in\calC$, we have $\bu-\mu\bx\in \calC-\calC$, and hence  $\frac{\bu-\mu\bx}{\|\bu-\mu\bx\|_{2}}\in   \calC_{(1)}$. Therefore, uniformly for all $\bu\in\calC$, 
    \begin{align*}
        \Xi_1&=\|\bu-\mu\bx\|_2\cdot\sup_{\bw\in\calK}\, \bw^\top \bigg(\frac{1}{n}\sum_{i=1}^n\ba_i\ba_i^\top-\bI_d\bigg) \frac{\bu-\mu\bx}{\|\bu-\mu\bx\|_2}
        \\&\le  \|\bu-\mu\bx\|_{2}\cdot\sup_{\bs\in \calC_{(1)}}\sup_{\bw\in \calK}\,\bigg|\frac{1}{n}\sum_{i=1}^n \bs^\top\ba_i\ba_i^\top\bw - \bs^\top\bw \bigg|:=\Xi_{1,1}\cdot\|\bu-\mu\bx\|_2.
    \end{align*}
    By Lemma \ref{lem:product_process} with ``$\calU=\calC_{(1)},~\calV=\calK,~K_{\calU}=O(1),~r_{\calU}=O(1),~K_{\calV}=O(1),~r_{\calV}=O(1)$,''
    \[\mathbb{P}\bigg(\Xi_{1,1}\lesssim \frac{\omega(\calK)+\omega(\calC_{(1)})+t}{\sqrt{n}}+\frac{(\omega(\calK)+t)(\omega(\calC_{(1)})+t)}{n}\bigg)\ge 1-2\exp(-t^2),\qquad \forall t>0. 
    \] 
    Under $n\gtrsim \omega^2(\calK)+\omega^2(\calC_{(1)})$, we   set $t=c_1\sqrt{n}$ with small enough $c_1$ to establish 
    \[\mathbb{P}\bigg(\Xi_{1,1}\le\frac{1}{20}\bigg)\ge 1-2\exp(-cn).\]
    Hence, $\Xi_1\le \frac{1}{20}\|\bu-\mu\bx\|_2$ holds for all $\bu\in\calC$ with the promised probability.

\paragraph{Bounding $\Xi_2$.} Under (\ref{positivemu}), rotational invariance of $\ba_i$ ensures $\mathbbm{E}\big( f_i(\ba_i^\top\bx)-\mu\ba_i^\top\bx\big)\ba_i=0$; see (\ref{zeromeansim}). Hence \[\Xi_2 = \sup_{\bw\in\calK}\,\frac{1}{n}\sum_{i=1}^n\big(f_i(\ba_i^\top\bx)-\mu\ba_i^\top\bx\big)\ba_i^\top\bw\]
is the supremum of a centered random process.  By (\ref{psi2bound}) we also have \(\| f_i(\ba_i^\top\bx)-\mu\ba_i^\top\bx\|_{\psi_2}\le \sigma .\) Since $\bx\in\calC\cap \mathbb{S}^{d-1}$ is fixed, Lemma \ref{lem:product_process} with ``$\calU=\{\bx\},~K_{\calU}=0,~r_{\calU}=\sigma,~\calV=\calK,~K_{\calV}=O(1),~r_{\calV}=O(1)$'' yields  
    \[
    \mathbb{P}\bigg(\Xi_2 \lesssim \frac{\sigma t(\omega(\calK)+t)}{n}+\frac{\sigma[\omega(\calK)+t]}{\sqrt{n}} \bigg)\ge 1-2\exp(-t^2),\qquad \forall t\ge 1.
    \]
   Setting $t=\delta\sqrt{n}$ yields
  \(\mathbb{P}\Big(\Xi_2 \lesssim \sigma \big(\frac{\omega(\calK)}{\sqrt{n}}+\delta\big)\Big)\ge 1-2\exp(-\delta^2n).\) 
    Substituting the bounds on $\Xi_1,\Xi_2$ into (\ref{tristartsim}) yields the desired RAIC (\ref{simraic}). 
\end{proof}

We are now ready to specialize Theorem \ref{thm:raicsim} to the T-RAIC in Theorem \ref{pro:traicsim}. 

\begin{proof}[Proof of Theorem \ref{pro:traicsim}] By a simple vectorization, we can identify $\mathbb{R}^{d_1\times d_2\times d_3}$ with the vector space $\mathbb{R}^{d_1d_2d_3}$. As such, we   invoke Theorem \ref{thm:raicsim} with $\calC= T^{\bd}_{\br}$, $\calK= T^{\bd}_{2\br}\cap \mathbb{B}_{\rm F}$, $\delta= \sqrt{T_{\rm df}/n}$, along with the observation $(T^{\bd}_{\br}-T^{\bd}_{\br})\cap \mathbb{B}_{\rm F} \subset T^{\bd}_{2\br}\cap \mathbb{B}_{\rm F}$, to obtain the following: if $n\gtrsim T_{\rm df}+\omega^2(T^{\bd}_{2\br}\cap \mathbb{B}_{\rm F})$, then the tensor RAIC claimed in Theorem \ref{pro:traicsim} holds  with probability at least $1-C\exp(-cT_{\rm df})$. The proof is concluded by using $\omega^2(T^{\bd}_{2\br}\cap \mathbb{B}_{\rm F})\lesssim T_{\rm df}$ from Proposition \ref{pro:tensorbound}. 
\end{proof}

\subsection{Logistic Regression (Proof of Theorem \ref{traicglm})}\label{app:tenlogi}

Note that the vector counterpart of  (\ref{grainglm}) is given by 
\begin{align}\label{hlogivec}
    \bh(\bu):= \frac{1}{n}\sum_{i=1}^n\big(s(\ba_i^\top\bu)-y_i\big)\ba_i,
\end{align}
where the observations in logistic regression satisfy the following. 

\begin{assumption}
    \label{assumplogis}
      Conditioning on $\ba_i$'s, the responses $y_1,\cdots,y_n$ are independent and follow $y_i\mid \ba_i \sim {\rm Bernoulli}(s(\langle\ba_i,\bx\rangle))$, where $s(a)=(1+\exp(-a))^{-1}$ is the sigmoid function.
\end{assumption}
 We establish the general RAIC of $\bh(\bu)$ in the following.

\begin{theorem}[RAIC for Logistic Regression] \label{thm:raicglm}  Fix $\bx\in \calC\cap\mathbb{S}^{d-1}$ for a cone $\calC$ and let $\calK$ be a symmetric set obeying $\calK\subset \mathbb{B}_2^d$. Under Assumptions \ref{gaussianaiassump} and \ref{assumplogis}, we define $\bh(\bu)$ as in (\ref{hlogivec}) and let $\eta:=[\mathbbm{E}_{g\sim N(0,1)}(s'(g))]^{-1}$. There exist universal constants $C>c>0$ such that the following holds. For any $\delta\in(0,1)$, if $n\ge C[\omega^2(\calC_{(1)})+\omega^2(\calK)]$, then with probability at least $1-C\exp(-c\delta^2n)$,    \[\eta\bh(\bu) \sim {\rm  RAIC} \bigg(\calK;\calC\cap \mathbb{S}^{d-1},\bx,C\bigg(\frac{\omega(\calK)+\omega(\calC\cap \mathbb{S}^{d-1})}{\sqrt{n}}+\delta\bigg)\bigg).\] 
\end{theorem}
\begin{proof}
   By definition, we only need to establish 
\begin{align}\label{logi-raic}
    \big\|\bu-\bx-\eta\bh(\bu)\big\|_{\calK^\circ} \le  C\bigg(\frac{\omega(\calK)+\omega(\calC\cap\mathbb{S}^{d-1})}{\sqrt{n}}+\delta\bigg),\qquad \forall \bu \in \calC\cap\mathbb{S}^{d-1}. 
\end{align} We start with
\begin{align}\label{logitrian}
    \big\|\bu-\bx-\eta\bh(\bu)\big\|_{\calK^\circ} \le \underbrace{\eta\big\|\bh(\bu)- \mathbb{E}[\bh (\bu)]\big\|_{\calK^\circ}}_{:=\Xi_3}+\underbrace{\big\|\eta \mathbb{E}[\bh(\bu)] - (\bu-\bx)\big\|_{\calK^\circ}}_{:=\Xi_4}. 
\end{align}

\paragraph{Bounding $\Xi_3$.} Since $\eta=O(1)$, for a fixed $\bx$, it remains to bound 
\begin{align*} 
     &\sup_{\bu\in \calC\cap \mathbb{S}^{d-1}}
     \big\|\bh(\bu)-\mathbb{E}[\bh(\bu)]\big\|_{\calK^\circ}\\
     &\quad= \sup_{\bu\in \calC\cap \mathbb{S}^{d-1}}
     \bigg\|\frac{1}{n}\sum_{i=1}^n(s(\ba_i^\top\bu)-y_i)\ba_i
     - \mathbb{E}[(s(\ba_i^\top\bu)-y_i)\ba_i]\bigg\|_{\calK ^\circ}\\
     &\quad= \sup_{\bu\in \calC\cap \mathbb{S}^{d-1}} \sup_{\bv\in\calK}
     \bigg\{\frac{1}{n}\sum_{i=1}^n(s(\ba_i^\top\bu)-y_i)\ba_i^\top\bv- \mathbb{E}[(s(\ba_i^\top\bu)-y_i)\ba_i^\top\bv]\bigg\}:=\Xi_{3,1}. 
\end{align*} 
For the fixed $\bx$ and any $\bu,\bu_1\in\calC\cap\mathbb{S}^{d-1}$, the boundedness of $s(\ba_i^\top\bu)$ and $y_i$ gives $ \|s(\ba_i^\top\bu)-y_i\|_{\psi_2} = O(1)$. Moreover, by $\sup_{a\in \mathbb{R}}s'(a)\le 1/4$, we have   
\begin{align*} 
    &\big\|\big[s(\ba_i^\top \bu)-y_i\big]-\big[s(\ba_i^\top \bu_1)-y_i\big]\big\|_{\psi_2}\\
    &\qquad= \|s(\ba_i^\top\bu)-s(\ba_i^\top\bu_1)\|_{\psi_2}
    \lesssim \|\ba_i^\top(\bu-\bu_1)\|_2 =O(\|\bu-\bu_1\|_2). 
\end{align*}
Obviously, for any $\bv,\bv_1\in\calK $, we have $\|\ba_i^\top\bv\|_{\psi_2}=O(1)$ and $\|\ba_i^\top\bv-\ba_i^\top\bv_1\|_{\psi_2}=O(1)\|\bv-\bv_1\|_{2}$. Therefore, by Lemma \ref{lem:product_process} with ``$\calU=\calC\cap \mathbb{S}^{d-1}$, $K_{\calU}=O(1)$,~$r_{\calU}=O(1)$,~$\calV=\calK$,~$K_{\calV}=O(1)$,~$r_{\calV}=O(1)$'',
\begin{align*}
     \mathbb{P}\bigg(\Xi_{3,1}\lesssim
     \frac{(\omega(\calK )+t)(\omega(\calC\cap \mathbb{S}^{d-1})+t)}{n}+ \frac{\omega(\calK)+\omega(\calC\cap\mathbb{S}^{d-1})+t}{\sqrt{n}} \bigg)
     \ge 1-2\exp(-ct^2),~~\forall t\ge 1.  
\end{align*} 
Under $n\ge\omega^2(\calK)+\omega^2(\calC\cap \mathbb{S}^{d-1})$, setting $t=\delta\sqrt{n}$ yields that $\Xi_{3,1}\lesssim \frac{\omega(\calK)+\omega(\calC\cap \mathbb{S}^{d-1})}{\sqrt{n}}+\delta$ holds with probability at least $1-2\exp(-c\delta^2n)$.

\paragraph{Showing $\Xi_4=0$.} By a direct application of the Stein's identity (see, e.g., \citealp[Appendix A.2]{chen2026finite}),
\[
    \mathbb{E}[\bh(\bu)] = \mathbb{E}_{g\sim N(0,1)}[s'(g)]\big(\bu-\bx\big)= \eta^{-1}(\bu-\bx).
\]
Combining with $\eta = (\mathbb{E}_{g\sim N(0,1)}[s'(g)])^{-1}$, it follows that $\Xi_4=0$. By (\ref{logitrian}) and the bound on $\Xi_3$, the claim follows. 
\end{proof}

We now prove Theorem \ref{traicglm}.
\begin{proof}[Proof of Theorem \ref{traicglm}] 
By identifying $\mathbb{R}^{d_1\times d_2\times d_3}$ with $\mathbb{R}^{d_1d_2d_3}$, we invoke Theorem \ref{thm:raicglm} with $\calC=T^{\bd}_{\br}$, $\calK=T^{\bd}_{2\br}\cap \mathbb{B}_{\rm F}$, and $\delta= \sqrt{T_{\rm df}/n}$ and apply the Gaussian width estimate in Proposition \ref{pro:tensorbound} to conclude the proof. 
\end{proof}

\subsection{Phase Retrieval (Proof of Theorem \ref{traicpr})}\label{app:raicpr}
Recall that we construct the gradient $\tcalH_{\tcalX^\star}(\tcalU)$ in (\ref{tprgradient}) for the recovery of $\tcalX^\star$. As such, the counterpart in $\mathbb{R}^d$ is given by 
\begin{align}
    \label{hxprvec}
    \bh_{\bx}(\bu) = \frac{1}{n}\sum_{i=1}^n \big(|\ba_i^\top\bu|-|\ba_i^\top\bx|\big)\sign(\ba_i^\top\bu)\ba_i. 
\end{align}
In the following, we establish the RAIC of $\bh_{\bx}(\bu)$ uniformly for all $\bx\in\calC\setminus\{0\}$.  

\begin{theorem}[RAIC for Phase Retrieval] \label{thm:prraic}
Assume that $\calC$ is a cone in $\mathbb{R}^d$, $\calK$ is a symmetric set being contained in $\mathbb{B}_2^d$. Under Assumption \ref{gaussianaiassump}, we let $\bh_{\bx}(\bu)$ be as in (\ref{hxprvec}). 
There exist universal constants $C>c >0$ such that the following holds. If   \(
    n\ge C \left(\log \calN(\calC\cap \mathbb{S}^{d-1};c)+\omega^2(\calK)+\omega^2(\calC_{(1)})\right),\) then with probability at least $1-C\exp(-cn)$,  
\[\bh_{\bx}(\bu)\sim {\rm RAIC}\bigg(\calK;\calC \cap\mathbb{B}^d_2(\bx;c\|\bx\|_2),\bx,\frac{1}{20}\|\bu-\bx\|_2\bigg),\qquad \forall \bx\in\calC\setminus\{0\}.\]  
\end{theorem}
\begin{proof} We seek to establish \begin{align}
    \label{prraic}
    \big\|\bu - \bx-  \bh_{\bx}(\bu)\big\|_{\calK^\circ} \le \frac{\|\bu-\bx\|_2}{20},\qquad \forall \bu \in \calC\cap \mathbb{B}^d_2(\bx;c\|\bx\|_2),~\bx\in\calC\setminus\{0\}.  
 \end{align} To this end we restrict our attention to $(\bu,\bx)$ obeying $\bx\ne0$ and \begin{align}
     \|\bu-\bx\|_2\le c_0\|\bx\|_2\label{c0place}
 \end{align} for some sufficiently small universal constant $c_0$. We define the index set \[\bR_{\bu,\bx}:=\big\{i\in [n]:\sign(\ba_i^\top\bu)\ne\sign(\ba_i^\top\bx)\big\}\] and   decompose $\bh_{\bx}(\bu)$ into 
 \begin{align*}
     \bh_{\bx}(\bu)
     &= \frac{1}{n}\sum_{i\notin\bR_{\bu,\bx}} \ba_i\ba_i^\top(\bu-\bx)
     + \frac{1}{n}\sum_{i\in\bR_{\bu,\bx}} \ba_i\ba_i^\top(\bu+\bx)= \underbrace{\frac{1}{n}\sum_{i=1}^n \ba_i\ba_i^\top(\bu-\bx)}_{\bh_1(\bu,\bx)}
     + \underbrace{\frac{2}{n}\sum_{i\in\bR_{\bu,\bx}}\ba_i\ba_i^\top\bx}_{\bh_2(\bu,\bx)}.
 \end{align*} 
By triangle inequality,
    \begin{align}\label{decompr}
        \big\|\bu-\bx-\bh_{\bx}(\bu)\big\|_{\calK ^\circ} \le \big\|\bu-\bx-\bh_1(\bu,\bx)\big\|_{\calK^\circ} + \|\bh_2(\bu,\bx)\|_{\calK^\circ}:=\Xi_5+\Xi_6.
    \end{align}
  \paragraph{Bounding $\Xi_5$.} Note that 
    $
        \Xi_5:= \Big\|\Big(\frac{1}{n}\sum_{i=1}^n\ba_i\ba_i^\top-\bI_d\Big)(\bu-\bx)\Big\|_{\calK^\circ}$ 
    is identical to the term $\Xi_1$ appearing in the proof of Theorem \ref{thm:raicsim}. Under $n\gtrsim \omega^2(\calK)+\omega^2(\calC_{(1)})$, repeating the arguments therein and changing the constants as needed yield  that 
    \begin{align}
        \Xi_5 \le\frac{1}{40}\|\bu-\bx\|_2, \qquad \forall (\bu,\bx)\in \calC\times (\calC\setminus\{0\})\label{t1boundpr}
    \end{align}
    holds 
    with probability at least $1-2\exp(-cn)$.  
\paragraph{Bounding $\Xi_6$.} We use triangle inequality to obtain
     \begin{align}\nn
         &\Xi_6 = \sup_{\bw\in\calK}\frac{2}{n}\sum_{i\in\bR_{\bu,\bx}}\bw^\top\ba_i\ba_i^\top\bx\\&\stackrel{(a)}{\le} \sup_{\bw\in\calK}\frac{2}{n}\sum_{i\in\bR_{\bu,\bx}}|\bw^\top\ba_i\ba_i^\top(\bu-\bx)|\nn \\
         &\stackrel{(b)}{\le} \|\bu-\bx\|_{2}\cdot \sup_{\bw\in\calK}\sup_{\bs\in \calC_{(1)}}\frac{2}{n}\sum_{i\in\bR_{\bu,\bx}}|\bw^\top\ba_i\ba_i^\top\bs|\nn\\ \nn 
         &\stackrel{(c)}{\le} \|\bu-\bx\|_{2}\cdot \bigg(\sup_{\bw\in\calK}\,\frac{1}{n}\sum_{i\in\bR_{\bu,\bx}}(\ba_i^\top\bw)^2+ \sup_{\bs\in\calC_{(1)}}\frac{1}{n}\sum_{i\in\bR_{\bu,\bx}}(\ba_i^\top\bs)^2\bigg)     
         \\:&\stackrel{(d)}{=} \|\bu-\bx\|_{2}\cdot \Big[\Xi_6'+\Xi_6''\Big],\label{116}
     \end{align}
     where in $(a)$ we observe that $\sign(\ba_i^\top\bu)\ne\sign(\ba_i^\top\bx)$ and hence $|\ba_i^\top\bx|\le |\ba_i^\top(\bu-\bx)|$ when $i\in\bR_{\bu,\bx}$, in $(b)$ we take supremum over $\bs:=\frac{\bu-\bx}{\|\bu-\bx\|_{2}}\in \calC_{(1)}$, in $(c)$ we use $2|\bw^\top\ba_i\ba_i^\top\bs|\le |\ba_i^\top\bw|^2+|\ba_i^\top\bs|^2$, and in $(d)$ we define 
     \begin{align}\label{defineT2p}
         \Xi_6':= \sup_{\bw\in\calK}\frac{1}{n}\sum_{i\in\bR_{\bu,\bx}}(\ba_i^\top\bw)^2,\qquad \Xi_6'':=\sup_{\bs\in\calC_{(1)}}\frac{1}{n}\sum_{i\in\bR_{\bu,\bx}}(\ba_i^\top\bs)^2.  
     \end{align}

     \paragraph{Bounding $\Xi_6'$, $\Xi_6''$.} We first bound $|\bR_{\bu,\bx}|$ uniformly by the {\it global binary embedding} results obtained in  \citet{oymak2015near,plan2014dimension}; a version sufficient for our need is restated in our Lemma \ref{lem:globinary}. We notice that for any $\bu,\bx\in\calC\setminus\{0\}$, we have $\bu/\|\bu\|_2,\bx/\|\bx\|_2\in\calC\cap \mathbb{S}^{d-1}$, and $\bR_{\bu,\bx}=\bR_{\bu/\|\bu\|_2,\bx/\|\bx\|_2}$. For the small enough $c_0$ in (\ref{c0place}), Lemma \ref{lem:globinary} has the following implication: for some universal constants $C>c>0$ depending on $c_0$ only, if $n\ge C(\log\calN(\calC \cap \mathbb{S}^{d-1},c)+\omega^2(\calC_{
     (1)
     }))$, then  
     \begin{align}\label{c113}
          &\left|\frac{|\bR_{\bu,\bx}|}{n}- \frac{1}{\pi}\arccos\Big(\Big\langle\frac{\bu}{\|\bu\|_2},\frac{\bx}{\|\bx\|_2}\Big\rangle\Big)\right|\le c_0  ,\qquad\forall(\bu,\bx)\in(\calC\setminus\{0\})\times(\calC\setminus\{0\}) .
     \end{align} holds 
     with probability at least $1-2\exp(-cn)$. By the simple fact \[\frac{\arccos(\bu^\top\bv)}{\pi}\le\frac{\|\bu-\bv\|_2}{2},\qquad \forall\bu,\bv\in\mathbb{S}^{d-1}\] and a standard bound (see, e.g., \citealp[Equation (2.12)]{chen2024robust})  
     \[\bigg\|\frac{\bu}{\|\bu\|_2}-\frac{\bv}{\|\bv\|_2}\bigg\|_2\le \frac{2\|\bu-\bv\|_2}{\max\{\|\bu\|_2,\|\bv\|_2\}},\qquad \forall \bu,\bv\ne 0,\] along with $\|\bu-\bx\|_2\le c_0\|\bx\|_2$ as in (\ref{c0place}),    
     \[
   \frac{1}{\pi}\arccos\Big(\Big\langle\frac{\bu}{\|\bu\|_2},\frac{\bx}{\|\bx\|_2}\Big\rangle\Big) \le \frac{1}{2}\left\|\frac{\bu}{\|\bu\|_2}-\frac{\bx}{\|\bx\|_2}\right\|_2\le \frac{\|\bu-\bx\|_2}{\|\bx\|_2}\le c_0.
     \] 
     Thus, (\ref{c113}) leads to 
     \[
         |\bR_{\bu,\bx}|\le 2c_0n,\qquad \forall(\bu,\bx)\in(\calC\setminus\{0\})\times(\calC\setminus\{0\}).\]
     On this event, $\Xi_6',\,\Xi_6''$ in (\ref{defineT2p}) are controlled by 
     \begin{align}  \nn 
         \Xi_6'+\Xi_6''&\le \sup_{\bw\in\calK}\max_{\substack{I\subset[n]\\|I|\le 2c_0n}}\frac{2}{n}\sum_{i\in I}(\ba_i^\top\bw)^2 + \sup_{\bw\in\calC_{(1)}}\max_{\substack{I\subset[n]\\|I|\le 2c_0n}}\frac{2}{n}\sum_{i\in I}(\ba_i^\top\bw)^2\\\label{C115}
         &\stackrel{(a)}{\le} C_1\bigg( \frac{\omega^2(\calK)+\omega^2(\calC_{(1)})}{n}+c_0\log(e/c_0)\bigg)\stackrel{(b)}{\le} \frac{1}{40},
     \end{align}
     where   $(a)$ holds with probability at least $1-2\exp(-c_1c_0n\log(e/c_0))$ for some universal constants $C_1>c_1>0$ due to Lemma \ref{lem:max_ell_sum}, and under $n\gtrsim \omega^2(\calK)+\omega^2(\calC_{(1)})$, $(b)$ can be guaranteed by setting $c_0>0$ as a small enough universal constant. Substituting (\ref{C115}) into (\ref{116}) yields $\Xi_6\le \frac{1}{40}\|\bu-\bx\|_2$. Combining with (\ref{t1boundpr}) and (\ref{decompr}) completes the proof. 
\end{proof}
In the following, we prove Theorem \ref{traicpr}.
\begin{proof}[Proof of Theorem \ref{traicpr}] By identifying $\mathbb{R}^{d_1\times d_2\times d_3}$ with $\mathbb{R}^{d_1d_2d_3}$, we invoke Theorem \ref{thm:prraic} with $\calC=T^{\bd}_{\br}$, $\calK = T^{\bd}_{2\br}\cap \mathbb{B}_{\rm F}$ and use the bounds in Proposition \ref{pro:tensorbound} to conclude the proof.  
\end{proof}
 
\subsection{ReLU Regression (Proof of Theorem \ref{thm:traicrelu})}\label{app:raicrelu}
Recall that we proposed to work with $\calH_{\tcalX^\star}(\tcalU)$ in (\ref{relutengra}), whose counterpart in $\mathbb{R}^d$ is given by 
\begin{align}\label{reluhxuvec}
    \bh_{\bx}(\bu) = \frac{1}{n}\sum_{i=1}^n\big[\relu(\ba_i^\top\bu)-\relu(\ba_i^\top\bx)\big]\big[1+\sign(\ba_i^\top\bu)\big]\ba_i.
\end{align}
Here, $\relu(a)=\max\{a,0\}=\frac{1}{2}(a+|a|)$. For some cone $\calC$, the following theorem shows that $\bh_{\bx}(\bu)$ satisfies a general RAIC, uniformly for all $\bx\in\calC\setminus\{0\}$.

\begin{theorem}[RAIC for ReLU Regression]\label{thm:raicrelusparse}
   Let $\calC$ be a cone in $\mathbb{R}^d$, and $\calK$ be a symmetric set being contained in $\mathbb{B}_2^d.$ Under Assumption \ref{gaussianaiassump}, we define  
   $\bh_{\bx}(\bu)$ as in (\ref{reluhxuvec}). Then there exist universal constants $C>c>0$ such that, if \[n\ge C\big(\log\calN(\calC\cap \mathbb{S}^{d-1};c)+\omega^2(\calK)+\omega^2(\calC_{(1)})\big),\] then with probability at least $1-C\exp(-cn)$,  \[\bh_{\bx}(\bu)\sim{\rm RAIC}\bigg(\calK;\calC\cap \mathbb{B}_2^d(\bx;c\|\bx\|_2),\bx,\frac{\|\bu-\bx\|_2}{20}\bigg),\qquad \forall\bx\in\calC\setminus\{0\}.\]
\end{theorem}
\begin{proof}
    We only need to prove 
   \begin{align}\label{reluraic}
        \big\|\bu-\bx-\bh_{\bx}(\bu)\big\|_{\calK^\circ}\le \frac{\|\bu-\bx\|_2}{20},\qquad \forall  \bu\in \calC\cap \mathbb{B}_2^d(\bx;c_0\|\bx\|_2),~\bx\in\calC\setminus\{0\},
    \end{align}where $c_0>0$ is some sufficiently small universal constant.  By triangle inequality,  
  \[
       \|\bu-\bx-\bh_{\bx}(\bu)\|_{\calK^\circ}\le   \underbrace{\bigg\|\bigg(\bI_d - \frac{1}{n}\sum_{i=1}^n \ba_i\ba_i^\top\bigg)(\bu-\bx)\bigg\|_{\calK^\circ}}_{:=\Xi_7}+\underbrace{\bigg\|\frac{1}{n}\sum_{i=1}^n\ba_i\ba_i^\top(\bu-\bx) -\bh_{\bx}(\bu)\bigg\|_{\calK^\circ}}_{:=\Xi_8}.
       \] 

       \paragraph{Bounding $\Xi_7$.} By reiterating the argument in the part of ``bounding $\Xi_1$'' in the proof of Theorem \ref{thm:raicsim} (with   changed   constants as needed), we obtain that, under $n\gtrsim \omega^2(\calK)+\omega^2(\calC_{(1)})$, $\Xi_7\le \frac{1}{80}\|\bu-\bx\|_2$ holds uniformly for all $\bx,\bu\in\calC$ with probability at least $1-2\exp(-cn)$.

  \paragraph{Bounding $\Xi_8.$} Substituting $\bh_{\bx}(\bu)$ and writing $\relu(a)$ as $a\mathbbm{1}(a\ge 0)$, 
   \begin{align*}
       &\frac{1}{n}\sum_{i=1}^n\ba_i\ba_i^\top(\bu-\bx) - \bh_{\bx}(\bu)\\
       &= \frac{1}{n}\sum_{i=1}^n
       \left( 
       \ba_i^\top(\bu-\bx) -\big[\ba_i^\top\bu \mathbbm{1}(\ba_i^\top\bu\ge 0)
       -\ba_i^\top\bx\mathbbm{1}(\ba_i^\top\bx\ge 0)\big]
       \big[1+\sign(\ba_i^\top\bu)\big]
      \right)\ba_i\\
       &= \frac{1}{n}\sum_{i=1}^n
       \left( 
       \ba_i^\top\bu-\ba_i^\top\bx
       -2 \ba_i^\top\bu \mathbbm{1}(\ba_i^\top\bu\ge 0) +\ba_i^\top\bx\mathbbm{1}(\ba_i^\top\bx\ge 0)
       +\ba_i^\top\bx\mathbbm{1}(\ba_i^\top\bx\ge 0)\sign(\ba_i^\top\bu)
      \right)\ba_i \\
       &= \frac{1}{n}\sum_{i=1}^n
       \underbrace{\left( 
       -|\ba_i^\top\bu|-\ba_i^\top\bx\mathbbm{1}(\ba_i^\top\bx<0)+\ba_i^\top\bx\mathbbm{1}(\ba_i^\top\bx\ge 0)\sign(\ba_i^\top\bu)
        \right)}_{:=F^{(i)}_{\bu,\bx}}\ba_i  
   \end{align*}
   We now introduce $\bR_{\bu,\bv}=\{i\in[n]:\sign(\ba_i^\top\bu)\ne\sign(\ba_i^\top\bv)\}$. It follows from some algebra that 
   \[F^{(i)}_{\bu,\bx}=
   \begin{cases} 
        -|\ba_i^\top\bu| - \ba_i^\top\bx \mathbbm{1}(\ba_i^\top\bx<0) 
        - \ba_i^\top\bx\mathbbm{1}(\ba_i^\top\bx\ge 0)
       = -|\ba_i^\top\bu| - \ba_i^\top\bx
      ,&\textrm{if}~i\in \bR_{\bu,\bx}, \\ 
       -|\ba_i^\top\bu|- \ba_i^\top\bx \mathbbm{1}(\ba_i^\top\bx<0)
       + \ba_i^\top\bx\mathbbm{1}(\ba_i^\top\bx\ge 0)
       =\sign(\ba_i^\top\bx)\ba_i^\top(\bx-\bu)
   ,&\textrm{if}~i\notin  \bR_{\bu,\bx}
   \end{cases}.\]
   Therefore
   \begin{align*}
       &\frac{1}{n}\sum_{i=1}^n\ba_i\ba_i^\top(\bu-\bx) - \bh_{\bx}(\bu) = \frac{1}{n}\sum_{i=1}^n F^{(i)}_{\bu,\bx}\ba_i \\
       & = \frac{1}{n}\sum_{i\in \bR_{\bu,\bx}}(-|\ba_i^\top\bu|-\ba_i^\top\bx)\ba_i + \frac{1}{n}\sum_{i\notin \bR_{\bu,\bx}} \sign(\ba_i^\top\bx)\ba_i^\top(\bx-\bu)\ba_i \\
       & = \frac{1}{n}\sum_{i\in \bR_{\bu,\bx}}[-|\ba_i^\top\bu|-\ba_i^\top\bx+\sign(\ba_i^\top\bx)\ba_i^\top(\bu-\bx)]\ba_i + \frac{1}{n}\sum_{i=1}^n  \sign(\ba_i^\top\bx)\ba_i^\top(\bx-\bu)\ba_i\\
       & = - \frac{1}{n}\sum_{i\in \bR_{\bu,\bx}}(\ba_i^\top\bx+|\ba_i^\top\bx|+2|\ba_i^\top\bu|)\ba_i + \frac{1}{n}\sum_{i=1}^n  \sign(\ba_i^\top\bx)\ba_i^\top(\bx-\bu)\ba_i,
   \end{align*}
   which implies 
   \begin{align}
       \Xi_8\le \bigg\|\frac{1}{n}\sum_{i\in \bR_{\bu,\bx}}(\ba_i^\top\bx+|\ba_i^\top\bx|+2|\ba_i^\top\bu|)\ba_i\bigg\|_{\calK^\circ}+\bigg\|\frac{1}{n}\sum_{i=1}^n  \sign(\ba_i^\top\bx)\ba_i^\top(\bx-\bu)\ba_i\bigg\|_{\calK^\circ}:=\Xi_{8,1}+\Xi_{8,2}. \label{reluraic3}
   \end{align}

 \paragraph{Bounding $\Xi_{8,1}$.} We start with  
   \begin{align*}
       &\Xi_{8,1} =\sup_{\bw\in\calK}\frac{1}{n}\sum_{i\in \bR_{\bu,\bx}}(\ba_i^\top\bx+|\ba_i^\top\bx|+2|\ba_i^\top\bu|)\ba_i^\top\bw \\
       &\stackrel{(a)}{\le}  \sup_{\bw\in\calK}\frac{1}{n}\sum_{i\in \bR_{\bu,\bx}}2\big(|\ba_i^\top\bx|+|\ba_i^\top\bu|\big)|\ba_i^\top\bw|
       \\&\stackrel{(b)}{\le}  \sup_{\bw\in\calK}\frac{1}{n}\sum_{i\in \bR_{\bu,\bx}}4|\ba_i^\top(\bu-\bx)||\ba_i^\top\bw|\\
       &\stackrel{(c)}{\le} \|\bu-\bx\|_2 \cdot \sup_{\bw\in\calK}\sup_{\bv\in \calC_{(1)}}\frac{4}{n}\sum_{i\in \bR_{\bu,\bx}}|\ba_i^\top\bv||\ba_i^\top\bw|\\
       &\stackrel{(d)}{\le} \|\bu-\bx\|_2 \cdot \bigg(\sup_{\bw \in\calK} \frac{2}{n}\sum_{i\in \bR_{\bu,\bx}}(\ba_i^\top\bw)^2 + \sup_{\bw \in\calC_{(1)}} \frac{2}{n}\sum_{i\in \bR_{\bu,\bx}}(\ba_i^\top\bw)^2 \bigg)\stackrel{(e)}{\le}\frac{\|\bu-\bx\|_2}{40}, 
   \end{align*}
   where $(a)$ is due to triangle inequality, $(b)$ holds because $\max\{|\ba_i^\top\bx|,|\ba_i^\top\bu|\}\le |\ba_i^\top(\bu-\bx)|$ when $i\in \bR_{\bu,\bx}$, in $(c)$ we suppose $\bu\ne \bx$ (otherwise, we are done) and take supremum over $\bv= \frac{\bu-\bx}{\|\bu-\bx\|_2}\in \calC_{(1)}$,  in $(d)$ we use $2|\ba_i^\top\bv||\ba_i^\top\bw|\le |\ba_i^\top\bv|^2+|\ba_i^\top\bw|^2$, in $(e)$ we reuse the bound (\ref{C115}) established in the proof of Theorem \ref{thm:prraic}, up to change of universal constants.

\paragraph{Bounding $\Xi_{8,2}$.}  
  For any $(\bu,\bx)\in \calC\times(\calC\setminus\{0\})$, 
    \begin{align*}       
    \Xi_{8,2}&=\bigg\|\frac{1}{n}\sum_{i=1}^n  \sign(\ba_i^\top\bx)\ba_i^\top(\bx-\bu)\ba_i\bigg\|_{\calK^\circ}  = \sup_{\bw\in\calK}\frac{1}{n}\sum_{i=1}^n  \sign(\ba_i^\top\bx)\ba_i^\top(\bx-\bu)\ba_i^\top\bw 
    \\
    &\stackrel{(a)}{\le} \|\bu-\bx\|_2 \sup_{\bw\in\calK}\sup_{\bv\in\calC_{(1)}}\frac{1}{n}\sum_{i=1}^n  \sign(\ba_i^\top\bx)\bv^\top\ba_i\ba_i^\top\bw 
    \\&\stackrel{(b)}{\le} \|\bu-\bx\|_2 \sup_{\bx\in \calC\cap \mathbb{S}^{d-1}} \underbrace{\sup_{\bw\in\calK}\sup_{\bv\in\calC_{(1)}}\frac{1}{n}\sum_{i=1}^n  \sign(\ba_i^\top\bx)\bv^\top\ba_i\ba_i^\top\bw}_{:=\xi(\bx)}, 
    \end{align*}
    where in $(a)$ we suppose $\bu\ne \bx$ and take supremum over $\bv:=\frac{\bx-\bu}{\|\bx-\bu\|_2}\in \calC_{(1)}$, $(b)$ holds because $\sign(\ba_i^\top\bx)=\sign(\ba_i^\top\bx/\|\bx\|_2)$ and $\bx/\|\bx\|_2\in\calC\cap \mathbb{S}^{d-1}$. 
    We seek to show $\sup_{\bx\in\calC\cap \mathbb{S}^{d-1}}\xi(\bx) \le\frac{1}{40}$. We accomplish this by a covering argument. For some small enough universal constant $\hat{c}$, we let $N_{\hat{c}}$ be a minimal $\hat{c}$-net of $\calC\cap \mathbb{S}^{d-1}$ such that $ |N_{\hat{c}}|= \calN(\calC\cap \mathbb{S}^{d-1};\hat{c})$. For each fixed $\bx\in \calC\cap \mathbb{S}^{d-1}$, notice that $\mathbb{E}[\sign(\ba_i^\top\bx)\bv^\top\ba_i\ba_i^\top\bw]=0$. By viewing $\xi(\bx)$ as a centered product process consisting of two factors $\sign(\ba_i^\top\bx)\ba_i^\top\bv$ and $\ba_i^\top\bw$, Lemma \ref{lem:product_process} with ``$\calU=\calC_{(1)}$, $K_{\calU}=O(1)$, $r_{\calU}=O(1)$,~$\calV=\calK$,~$K_{\calV}=O(1)$,~$r_{\calV}=O(1)$'' establishes 
    \[\mathbb{P}\bigg(\xi(\bx)\lesssim  \frac{\big(\omega(\calC_{(1)})+t\big)\big(\omega(\calK)+t\big)}{n}+\frac{\omega(\calC_{(1)})+\omega(\calK)+t}{\sqrt{n}}\bigg)\ge 1-2 \exp(-c't^2),\qquad\forall  t\ge 1.\]
    Therefore a union bound over $\bx\in N_{\hat{c}}$, along with $t=c_0\sqrt{n}$ for sufficiently small universal constant $c_0$ and $n\gtrsim \omega^2(\calC_{(1)})+\omega^2(\calK)+ \log\calN(\calC\cap \mathbb{S}^{d-1};\hat{c})$, yields
    \begin{align}
        \mathbb{P}\bigg(\sup_{\bx\in N_{\hat{c}}}\xi(\bx)\le \frac{1}{160}\bigg)\ge 1- 2\exp(-c'n/2).\label{netboundrelu}
    \end{align}
    For any $\bx\in\calC\cap \mathbb{S}^{d-1}$, we let $\bx' = \textrm{arg}\min_{\bu\in N_{\hat{c}}}\|\bu-\bx\|_2$ with ties broken arbitrarily. Then 
    \begin{align*}
        \sup_{\bx\in \calC\cap \mathbb{S}^{d-1}}\xi(\bx)\le \sup_{\bx\in \calC\cap \mathbb{S}^{d-1}}\xi(\bx') + \sup_{\bx\in \calC\cap \mathbb{S}^{d-1}}[\xi(\bx)-\xi(\bx')]\le \frac{1}{160}+  \sup_{\bx\in \calC\cap \mathbb{S}^{d-1}}[\xi(\bx)-\xi(\bx')], 
    \end{align*}
    Moreover, using $\bR_{\bx,\bx'}=\{i\in[n]:\sign(\ba_i^\top\bx)\ne\sign(\ba_i^\top\bx')\}$, 
    \begin{align*}
         \sup_{\bx\in \calC\cap \mathbb{S}^{d-1}}[\xi(\bx)-\xi(\bx')]&\stackrel{(a)}{\le} \sup_{\bx\in \calC\cap \mathbb{S}^{d-1}}\sup_{\substack{\bw\in\calK\\ \bv\in\calC_{(1)}}}\frac{1}{n}\sum_{i=1}^n [\sign(\ba_i^\top\bx)-\sign(\ba_i^\top\bx')]\bv^\top\ba_i\ba_i^\top\bw
        \\
        &\stackrel{(b)}{\le} \sup_{\bx\in \calC\cap \mathbb{S}^{d-1}} \sup_{\substack{\bw\in\calK\\ \bv\in \calC_{(1)}}} \frac{2}{n}\sum_{i\in \bR_{\bx,\bx'}}|\ba_i^\top\bv||\ba_i^\top\bw|
        \\&\stackrel{(c)}{\le} \sup_{\bx\in \calC\cap \mathbb{S}^{d-1}}\sup_{\substack{\bw\in \calK\\ \bv\in \calC_{(1)}}} \frac{1}{n}\sum_{i\in \bR_{\bx,\bx'}}\big[(\ba_i^\top\bw)^2+(\ba_i^\top\bv)^2\big]\stackrel{(d)}{\le} \frac{1}{160}, 
    \end{align*}
    where $(a)$ follows from $\sup A-\sup B\le \sup[A-B]$, $(b)$ follows from triangle inequality and the definition of $\bR_{\bx,\bx'}$, in $(c)$ we use $|\ba_i^\top\bv||\ba_i^\top\bw|\le\frac{1}{2}[(\ba_i^\top\bv)^2+(\ba_i^\top\bw)^2]$, and $(d)$ can be ensured, with the promised probability, by repeating the arguments in (\ref{c113})--(\ref{C115}) along with the small enough $\hat{c}$. Combining these developments, we obtain $\Xi_{8,2}\le \frac{1}{80}\|\bu-\bx\|_2$.  Further combining with $\Xi_{8,1}\le \frac{\|\bu-\bx\|_2}{40}$ and $\Xi_7\le \frac{1}{80}\|\bu-\bx\|_2$, we arrive at the desired (\ref{reluraic}). 
\end{proof}

\begin{proof}[Proof of Theorem \ref{thm:traicrelu}] By identifying $\mathbb{R}^{d_1\times d_2\times d_3}$ with $\mathbb{R}^{d_1d_2d_3}$, we invoke Theorem \ref{thm:raicrelusparse} with $\calC = T^{\bd}_{\br},$ $\calK= T^{\bd}_{2\br}\cap \mathbb{B}_{\rm F}$ and use the estimates in Proposition \ref{pro:tensorbound} to conclude the proof.  
\end{proof}

\subsection{One-Bit Compressed Sensing (Proof of Theorem \ref{traic1bcs})}\label{app:raic1bcs}
The counterpart of (\ref{1bcsgrad}) in $\mathbb{R}^d$ is given by
\begin{align} \label{1bcshuvec}
    \bh_{\bx}(\bu):=\frac{1}{n}\sum_{i=1}^n\big[\sign(\ba_i^\top\bu)-\sign(\ba_i^\top\bx)\big]\ba_i.
\end{align}
The following RAIC for a general cone is adapted from \citet[Theorem 5]{chen2024optimal}. 

\begin{theorem}
    [RAIC for 1bCS] \label{thm:1bcsraic}
    Let $\calC$ be a cone in $\mathbb{R}^d$. Under Assumption \ref{gaussianaiassump}, let $\bh_{\bx}(\bu)$ be defined in (\ref{1bcshuvec}). There exist universal constants $C>c>0$ such that the following holds. Given any $\delta\in(0,c)$, 
 if \[
      n\ge C\frac{\omega^2(\calC_{(1)})+\log\calN(\calC\cap \mathbb{S}^{d-1};\delta/2)}{\delta}, 
 \]
  then with probability at least $1-C\exp(-c\log\calN(\calC\cap \mathbb{S}^{d-1};\delta))$,   \[
           \sqrt{\frac{\pi}{2}}\bh_{\bx}(\bu)\sim {\rm RAIC}\bigg(\calC_{(1)};\calC\cap \mathbb{S}^{d-1},\bx,\sqrt{C\delta\|\bu-\bx\|_2}+ C\delta\sqrt{\log(1/\delta)} \bigg),\qquad \forall \bx\in \calC\cap \mathbb{S}^{d-1}.\]
\end{theorem}
\begin{proof}
   This is a consequence of the more general RAIC for a star-shaped set in \citet[Theorem 5]{chen2024optimal}. Using the notation in \citet{chen2024optimal}, in the specific one-bit compressed sensing model, we have $(\Delta,\Lambda)=(2,0)$, $\mu_4=2$, and $\varepsilon^{(i)}=0$ for $i\in \{1,2,3,4\}$; see \citep[Section 4.1]{chen2024optimal}. While $\calK$ in \citet{chen2024optimal} can be a star-shaped set (that is more general than a cone), when specialized to a cone, the sample complexity in their Equation (46) reduces to 
    \[ n \gtrsim  \frac{\omega^2(\calC_{(1)})+\log\calN(\calC\cap \mathbb{S}^{d-1};r/2)}{r}\]
    for some small enough $r>0$. Under such sample complexity,  Theorem 5 therein gives an RAIC that reads (see \citealp[Definition 1]{chen2024optimal} for their definition of RAIC) 
    \begin{align*}
       \|\bu-\bx-\sqrt{\pi/2}\cdot\bh_{\bx}(\bu)\|_{\calC_{(1)}^\circ}\le  \sqrt{O(r)\|\bu-\bx\|_2}+ O\big(r\sqrt{\log(1/r)}\big),\qquad  \forall \bu,\bx\in\calC\cap \mathbb{S}^{d-1} 
     \end{align*}
     with probability at least $1-\exp(-c\log\calN(\calC\cap \mathbb{S}^{d-1};r))$. By replacing $r$ with $\delta$, we arrive at the claim.  
\end{proof}
\begin{proof}[Proof of Theorem \ref{traic1bcs}] By identifying $\mathbb{R}^{d_1\times d_2\times d_3}$ with $\mathbb{R}^{d_1d_2d_3}$, we invoke Theorem \ref{thm:1bcsraic} with $\calC=T^{\bd}_{2\br}$. Due to $T^{\bd}_{2\br}\cap \mathbb{B}_{\rm F} \subset \calC_{(1)}\subset T^{\bd}_{4\br}\cap \mathbb{B}_{\rm F}$, this yields the following: if for some small $\delta$ it holds that
\begin{align}\label{t1bcssamp}
    n\gtrsim \delta^{-1}\big(\omega^2(T^{\bd}_{4\br}\cap \mathbb{B}_{\rm F}) + \log\calN(T^{\bd}_{2\br}\cap \mathbb{S}_{\rm F};\delta/2)\big),
\end{align}
then with probability at least $1-C\exp(-c\log\calN(T^{\bd}_{2\br}\cap \mathbb{S}_{\rm F};\delta))$,
\begin{align}
    \sqrt{\frac{\pi}{2}}\tcalH_{\tcalX^\star}(\tcalU)\sim {\rm RAIC}\bigg(T^{\bd}_{2\br}\cap \mathbb{B}_{\rm F};T^{\bd}_{\br}\cap \mathbb{S}_{\rm F},\tcalX^\star,\sqrt{O(\delta)\|\tcalU-\tcalX^\star\|_{\rm F}}+O(\delta\sqrt{\log(1/\delta)})\bigg).\label{traicwithdelta}
\end{align}
By Proposition \ref{pro:tensorbound}, a sufficient condition for (\ref{t1bcssamp}) is $n\gtrsim \delta^{-1}T_{\rm df}\log(\delta^{-1})$, and so it is sufficient to choose $\delta\asymp \frac{T_{\rm df}}{n}\log\frac{n}{T_{\rm df}}$, which is small enough under $n\gtrsim T_{\rm df}$. Substituting this $\delta$ into (\ref{traicwithdelta}), along with $\sqrt{O(\delta)\|\tcalU-\tcalX^\star\|_{\rm F}}\le \frac{1}{40}\|\tcalU-\tcalX^\star\|_{\rm F}+ O (\delta)$, gives the desired RAIC. It is  also easy to show that $\log\calN(T^{\bd}_{2\br}\cap \mathbb{S}_{\rm F};\delta)\gtrsim T_{\rm df}$. The proof is complete.   
\end{proof}

\section{Proof of Spectral Initialization (Theorem \ref{thm:initsim},  Corollary \ref{corCpoint1})}
Throughout this section, the assumptions and notation are those of Section \ref{sec:speinitial}. By rotational invariance of $\tcalA_i$ and \eqref{Cpoint2}, 
\begin{equation}
\begin{aligned}
 \mathbb E[y_i\tcalA_i]
 &=\mathbb E\left[f_i(\langle\tcalA_i,\tcalX^\star\rangle)\tcalA_i \right]\\
 &=\frac{\mathbb E\left[f_i(\langle\tcalA_i,\tcalX^\star\rangle)
 \langle\tcalA_i,\tcalX^\star\rangle\right]}{\rho^2}\tcalX^\star= \frac{\mathbb{E}_{g\sim N(0,1)}[f_i(\rho g)g]}{\rho}\tcalX^\star
 =\mu\tcalX^\star.
\end{aligned}\label{Cpoint25}
\end{equation}
We first establish two concentration lemmas.

\begin{lem} \label{lemCdot1}
There exist universal constants $C>c>0$ such that, if $n\geq CT_{\mathrm{df}}$, then
\[
 \sup_{\tcalD\in T_{\br}^{\bd}\cap \mathbb{B}_{\rm F}}
 \bigg|
 \frac{1}{n}\sum_{i=1}^n
 \left\{
 y_i\langle\tcalA_i,\tcalD\rangle
 -\mathbb E\bigl[y_i\langle\tcalA_i,\tcalD\rangle\bigr]
 \right\}
 \bigg|
 \leq C\nu\sqrt{\frac{T_{\mathrm{df}}}{n}}   
\] 
with probability at least $1-C\exp(-cT_{\mathrm{df}})$.
\end{lem}

\begin{proof}
By $y_i = f_i(\langle\tcalA_i,\tcalX^\star\rangle)$ as assumed in Section \ref{sec:speinitial},  the process can be written as 
\[\Xi_9:=\sup_{\tcalD\in T_{\br}^{\bd}\cap \mathbb{B}_{\rm F}}\bigg|\frac{1}{n}\sum_{i=1}^nf_i(\langle\tcalA_i,\tcalX^\star\rangle)\langle\tcalA_i,\tcalD\rangle - \mathbb{E}\big[f_i(\langle\tcalA_i,\tcalX^\star\rangle)\langle\tcalA_i,\tcalD\rangle\big]\bigg|.\]
By identifying $\mathbb{R}^{d_1\times d_2\times d_3}$ with $\mathbb{R}^{d_1d_2d_3}$ and using Lemma \ref{lem:product_process} with $\calU = \{\tcalX^\star\}$, $K_{\calU}=0$ (because $\calU$ is a singleton), $r_{\calU} = \nu$ (due to (\ref{Cpoint2})), $\calV=T^{\bd}_{\br}\cap \mathbb{B}_{\rm F}$, $K_{\calV}=O(1)$, $r_{\calV}=O(1)$, we obtain 
\[\mathbb{P}\bigg(\Xi_9 \lesssim \frac{t\nu[\omega(T^{\bd}_{\br}\cap \mathbb{B}_{\rm F})+t]}{n}  + \frac{\nu[\omega(T^{\bd}_{\br}\cap \mathbb{B}_{\rm F})+t]}{\sqrt{n}}\bigg)\ge 1-2\exp(-ct^2),\qquad \forall t\ge 1.\]
Setting $t=\sqrt{T_{\rm df}}$ and using $\omega(T^{\bd}_{\br}\cap \mathbb{B}_{\rm F})=O(\sqrt{T_{\rm df}})$ from Proposition \ref{pro:tensorbound} conclude the proof. 
\end{proof}

In the following, we establish mode-wise concentration bounds for $ \widehat{\bG}_j:=\bY_j\bY_j^\top-\tau_j\bI_{d_j}$, $j=1,2,3.$ In Algorithm \ref{alg:ini}, we invoke the T-HOSVD in Algorithm \ref{alg:thosvd}. Hence, we can write the output of Algorithm \ref{alg:ini} as 
\begin{align}\label{inioutput}
    \tcalX_0 = H^{\bd}_{\br}(\tcalX_0^{(f)}) =\tcalX_0^{(f)}\times_{j=1}^3 (\widehat{\bU}_j\widehat{\bU}_j^\top), \qquad\textrm{where~$\widehat{\bU}_j={\rm SVD}_{r_j}(\calM_j(\tcalX_0^{(f)}))$.}
\end{align}

\begin{lem}\label{lemCdot3}
Fix $j\in\{1,2,3\}$ and let
\[
 \bY_j:=\calM_j(\tcalX_0^{(f)}),\qquad
 \tau_j:=\frac{d_{-j}-r_j}{n^2}\sum_{i=1}^n y_i^2,
 \qquad 
 \widehat{\bG}_j:=\bY_j\bY_j^\top-\tau_j\bI_{d_j}.\]
The leading $r_j$ eigenvectors of $\widehat{\bG}_j$ are the leading $r_j$ left singular vectors collected in $\widehat{\bU}_j$   in  (\ref{inioutput}).  If $n\geq Cd_j$, then, with probability at least $1-C\exp(-cd_j)$,
\begin{equation}
 \left\|\widehat{\bG}_j-
 \mu^2\calM_j(\tcalX^\star)\calM_j(\tcalX^\star)^\top\right\|_{\rm op}
 \leq C\bigg\{
 \nu^2\frac{d_j+\sqrt D}{n}
 +\mu\nu\|\calM_j(\tcalX^\star)\|_{\rm op}
 \sqrt{\frac{d_j}{n}}
 \bigg\}. \label{Cpoint29}
\end{equation}
\end{lem}

\begin{proof}
Fix $j\in\{1,2,3\}$ and abbreviate \(
 p:=d_j,~ q:=d_{-j},~ r:=r_j,
 ~\bX_j:=\calM_j(\tcalX^\star).\)
Let
\(
 \bX_j=\bU_j\boldsymbol{\Sigma}_j\bV_j^\top,~\bU_j\in \mathbb{R}^{p\times r},~\bSigma_j\in \mathbb{R}^{r\times r},~\bV_j\in\mathbb{R}^{q\times r}\)
be its compact singular value decomposition, and complete $\bV_j$ to an orthogonal matrix $[\bV_j,\bV_{j,\perp}]\in\mathbb{O}_{q\times q}$.  Define
\[
 \bB_i:=\calM_j(\tcalA_i)\bV_j\in\mathbb R^{p\times r},\qquad \bC_i:=\calM_j(\tcalA_i)\bV_{j,\perp}\in\mathbb R^{p\times(q-r)}.
\]
By Gaussian rotational invariance, $\bB_i$ and $\bC_i$ have \emph{independent} standard Gaussian entries and are \emph{independent of each other}.  Moreover, we have
\(\langle\tcalA_i,\tcalX^\star\rangle
 =\langle\bB_i,\bU_j\boldsymbol{\Sigma}_j\rangle.\) Consequently, $y_i$ is measurable with respect to $(\bB_i,f_i)$ and is independent of $\bC_i$.

Notice that  
\[
 \bY_{j,1}:=\bY_j\bV_j=\frac{1}{n}\sum_{i=1}^n y_i\bB_i,\qquad
 \bY_{j,2}:=\bY_j\bV_{j,\perp}=\frac{1}{n}\sum_{i=1}^n y_i\bC_i.
\]
Since $[\bV_j,\bV_{j,\perp}]$ is orthogonal,
\begin{equation}
 \bY_j\bY_j^\top = ( \bY_j[\bV_j,\bV_{j,\perp}])( \bY_j[\bV_j,\bV_{j,\perp}])^\top
 =\bY_{j,1}\bY_{j,1}^\top
 +\bY_{j,2}\bY_{j,2}^\top. \label{Cpoint30}
\end{equation}
The population identity \eqref{Cpoint25} gives
\[ \mathbb E \bY_{j,1}= 
 \mathbb E[y_i\bB_i] = \mathbb{E}[y_i\calM_j(\tcalA_i)\bV_j]=\mu  \calM_j(\tcalX^\star)\bV_j=\mu\bU_j\boldsymbol{\Sigma}_j.
\]

Set
\(
 \bE_j:=\bY_{j,1}-\mu\bU_j\boldsymbol{\Sigma}_j.
\) 
Subtracting the scalar matrix $\tau_j\bI_p$ does not change any eigenvector of $\bY_j\bY_j^\top$, which proves the first assertion of the lemma.  By \eqref{Cpoint30},
\[ 
 \widehat{\bG}_j-\mu^2\bX_j\bX_j^\top
 =\mu\bU_j\boldsymbol{\Sigma}_j\bE_j^\top
 +\mu\bE_j\boldsymbol{\Sigma}_j\bU_j^\top
 +\bE_j\bE_j^\top+\bY_{j,2}\bY_{j,2}^\top-\tau_j\bI_p.
\]
Therefore, 
 \begin{align}\label{tria1e3}
 \left\|\widehat{\bG}_j-\mu^2\bX_j\bX_j^\top\right\|_{\rm op}
   \leq 2\mu\|\bX_j\|_{\rm op}\|\bE_j\|_{\rm op}
 +\|\bE_j\|_{\rm op}^2 +\left\|\bY_{j,2}\bY_{j,2}^\top-\tau_j\bI_p\right\|_{\rm op}. 
 \end{align} 
 
\paragraph{Bounding $\|\bE_j\|_{\rm op}$.} We write $\bE_j = \bY_{j,1} - \mu\bU_j\bSigma_j = \frac{1}{n}\sum_{i=1}^nf_i(\langle\tcalA_i,\tcalX^\star\rangle)\bB_i- \mathbb{E}[f_i(\langle\tcalA_i,\tcalX^\star\rangle)\bB_i]$, and hence
\[\|\bE_j\|_{\rm op} = \sup_{\substack{\rank(\bD)=1\\\|\bD\|_{\rm F}=1}}\langle \bE_j, \bD\rangle = \sup_{\substack{\rank(\bD)=1\\\|\bD\|_{\rm F}=1}}  \frac{1}{n}\sum_{i=1}^nf_i(\langle\tcalA_i,\tcalX^\star\rangle)\langle\bB_i,\bD\rangle- \mathbb{E}[f_i(\langle\tcalA_i,\tcalX^\star\rangle)\langle\bB_i,\bD\rangle].\]
We now apply Lemma \ref{lem:product_process} with $\calU=\{\tcalX^\star\}$, $K_{\calU}=0$, $r_{\calU}=\nu$, $\calV=\{\bD\in \mathbb{R}^{p\times r}:\rank(\bD)=1,~\|\bD\|_{\rm F}=1\}$, $K_{\calV}=O(1)$, $r_{\calV}=O(1)$, along with the standard estimate $\omega(\calV)=O(\sqrt{p})$ (see, e.g., \citealp{chandrasekaran2012convex}), to obtain 
\[\mathbb{P}\bigg(\|\bE_j\|_{\rm op}\lesssim \Big[1+\frac{t}{\sqrt{n}}\Big]\frac{\nu(\sqrt{p}+t)}{\sqrt{n}}\bigg)\ge 1-2\exp(-ct^2),\qquad \forall t\ge 1.\]
Under $n\ge p$, setting $t=\sqrt{p}$ yields
\begin{align}
    \label{Ewbound}
    \mathbb{P}\big(\|\bE_j\|_{\rm op}\lesssim \nu\sqrt{p/n}\big) \ge 1-2\exp(-cp).
\end{align}

\paragraph{Bounding $\|\bY_{j,2}\bY_{j,2}^\top-\tau_j\bI_p\|_{\rm op}$.} Next, since $\|y_i\|_{\psi_2}\leq\nu$, we have $\|y_i^2\|_{\psi_1}=O(\nu^2)$, and therefore Bernstein's inequality (e.g., \citealp[Theorem 2.8.1]{vershynin2018high}) gives
\begin{equation}
\mathbb{P}\bigg( S_y:=\sum_{i=1}^n y_i^2\leq C\nu^2n\bigg)\ge 1-2\exp(-cn). \label{Cpoint32}
\end{equation}
 Conditioning on $(\bB_i,f_i)_{i=1}^n$, and hence on $(y_i)_{i=1}^n$, the entries of $\bY_{j,2}$ are independent Gaussian variables with variance $S_y/n^2$.  Thus, conditional on the responses, we have 
\(
 \bY_{j,2}\stackrel{d}{=}\frac{\sqrt{S_y}}{n}\bW,
\) 
where $\bW\in\mathbb R^{p\times(q-r)}$ has independent $N(0,1)$ entries.  Standard covariance estimation bound (e.g., \citealp[Exercise 4.7.3]{vershynin2018high}) implies
\[
 \mathbb{P}\Big(\|\bW\bW^\top-(q-r)\bI_p\|_{\rm op}
 \leq C\bigl(p+\sqrt{p(q-r)}\bigr)\Big)\ge 1-2\exp(-cp).
\]
Combining this with $\frac{S_y}{n^2}\lesssim \frac{\nu^2}{n}$ from \eqref{Cpoint32}, and then removing the conditioning, we obtain 
\begin{equation}
 \mathbb{P}\bigg(\left\|\bY_{j,2}\bY_{j,2}^\top-\tau_j\bI_p\right\|_{\rm op}
 \leq C\nu^2\frac{p+\sqrt{pq}}{n}\bigg)\ge 1-C\exp(-cp).  \label{Cpoint33}
\end{equation}

We are now ready to complete the proof. Substituting \eqref{Ewbound}, \eqref{Cpoint33} into (\ref{tria1e3}), and recalling $p=d_j$ and $pq=D$, we arrive at \eqref{Cpoint29}.
\end{proof}

\subsection{Proof of Theorem \ref{thm:initsim} (Performance of Algorithm \ref{alg:ini})}
\begin{proof}  
By \eqref{Cpoint25}, $\mathbb{E}\tcalX_0^{(f)}=\mu\tcalX^\star$. Let \(\tcalX^\star=\tcalC\times_1\bU_1\times_2\bU_2\times_3\bU_3\)
be a Tucker decomposition. For every $j\in\{1,2,3\}$, let $\widehat{\bU}_j$ contain the leading $r_j$ left singular vectors of $\calM_j(\tcalX_0^{(f)})$, and define
\(\bP_j:=\bU_j\bU_j^\top,
    ~
    \widehat{\bP}_j:=\widehat{\bU}_j\widehat{\bU}_j^\top.\)
By the definition of Algorithm \ref{alg:ini} and (\ref{inioutput}),
we have \(
    \tcalX_0
    =
    \tcalX_0^{(f)}\times_{j=1}^3\widehat{\bP}_j,
\)
whereas $\tcalX^\star=\tcalX^\star\times_{j=1}^3\bP_j$. Therefore, the triangle inequality gives
\begin{align} 
    \|\tcalX_0-\mu\tcalX^\star\|_{\rm F}
    \leq 
    \big\|
        \big(\tcalX_0^{(f)}-\mu\tcalX^\star\big)
        \times_{j=1}^3\widehat{\bP}_j
    \big\|_{\rm F} +
    \mu\big\|
        \tcalX^\star\times_{j=1}^3\widehat{\bP}_j-\tcalX^\star
    \big\|_{\rm F}
    =:\Xi_{10}+\mu\Xi_{11}. 
\label{Cpoint26}
\end{align}

\paragraph{Bounding $\Xi_{10}$.} Since the $\widehat{\bP}_j$'s are orthogonal projectors, the variational representation of the Frobenius norm and the self-adjointness of the mode-wise projections yield
\begin{align*}
    \Xi_{10}
    &=
    \sup_{\|\tcalD\|_{\rm F}\leq 1}
    \left|
        \left\langle
            \big(\tcalX_0^{(f)}-\mu\tcalX^\star\big)
            \times_{j=1}^3\widehat{\bP}_j,
            \tcalD
        \right\rangle
    \right|\\
    &=
    \sup_{\|\tcalD\|_{\rm F}\leq 1}
    \left|
        \left\langle
            \tcalX_0^{(f)}-\mu\tcalX^\star,
            \tcalD\times_{j=1}^3\widehat{\bP}_j
        \right\rangle
    \right|\leq
    \sup_{\substack{\widetilde{\tcalD}\in T_{\br}^{\bd}\\
                    \|\widetilde{\tcalD}\|_{\rm F}\leq 1}}
    \left|
        \left\langle
            \tcalX_0^{(f)}-\mu\tcalX^\star,
            \widetilde{\tcalD}
        \right\rangle
    \right|.
\end{align*}
Indeed, for every $\tcalD$ satisfying $\|\tcalD\|_{\rm F}\leq1$, the tensor
\(
    \widetilde{\tcalD}
    :=
    \tcalD\times_{j=1}^3\widehat{\bP}_j\) 
has Tucker rank at most $\br$ and satisfies
$\|\widetilde{\tcalD}\|_{\rm F}\leq\|\tcalD\|_{\rm F}\leq1$.
Consequently, Lemma \ref{lemCdot1} and \eqref{Cpoint25} imply
\begin{equation}
    \mathbb{P}\big(\Xi_{10}
    \leq
    C\nu\sqrt{T_{\rm df}/n}\big)\ge 1-C\exp(-cT_{\mathrm{df}}). 
    \label{X1bound11}
\end{equation} 

\paragraph{Bounding $\Xi_{11}.$} A telescoping expansion gives
 \begin{align*}
 &\tcalX^\star\times_1\widehat{\bP}_1\times_2\widehat{\bP}_2\times_3\widehat{\bP}_3
 -\tcalX^\star\times_1\bP_1\times_2\bP_2\times_3\bP_3\\
 & =\tcalX^\star\times_1(\widehat{\bP}_1-\bP_1)
 \times_2\widehat{\bP}_2\times_3\widehat{\bP}_3\\
 &\quad+\tcalX^\star\times_1\bP_1\times_2(\widehat{\bP}_2-\bP_2)
 \times_3\widehat{\bP}_3\\
 &\quad+\tcalX^\star\times_1\bP_1\times_2\bP_2
 \times_3(\widehat{\bP}_3-\bP_3).
 \end{align*} 
Since all $\bP_j$ and $\widehat{\bP}_j$ are orthogonal projectors,
\begin{equation}
 \Xi_{11}\leq\rho\sum_{j=1}^3
 \|\widehat{\bP}_j-\bP_j\|_{\rm op}. \label{Cpoint28}
\end{equation}
For each mode $j$, define $\widehat{\bG}_j$ as in Lemma \ref{lemCdot3}, then 
Lemma \ref{lemCdot3} yields
\[ 
 \left\|\widehat{\bG}_j-
 \mu^2\calM_j(\tcalX^\star)\calM_j(\tcalX^\star)^\top\right\|_{\rm op}
 \leq C\bigg\{
 \nu^2\frac{d_j+\sqrt D}{n}
 +\mu\nu s_j\sqrt{\frac{d_j}{n}}
 \bigg\}. 
\]
The signal eigengap of
$\mu^2\calM_j(\tcalX^\star)\calM_j(\tcalX^\star)^\top$ is $\mu^2\lambda_j^2$.  After increasing the universal constant in \eqref{Cpoint6}, that condition ensures
$\left\|\widehat{\bG}_j-
 \mu^2\calM_j(\tcalX^\star)\calM_j(\tcalX^\star)^\top\right\|_{\rm op}\leq\mu^2\lambda_j^2/2$ for every $j\in[3]$.  Since the leading $r_j$ eigenvectors of $\widehat{\bG}_j$ are precisely $\widehat{\bU}_j$, the Davis--Kahan theorem \citep{davis1970rotation} gives
\begin{equation}
 \|\widehat{\bP}_j-\bP_j\|_{\rm op}
 \leq\frac{2\big\|\widehat{\bG}_j-
 \mu^2\calM_j(\tcalX^\star)\calM_j(\tcalX^\star)^\top\big\|_{\rm op}}{\mu^2\lambda_j^2}
 \leq C\bigg\{
 \frac{\nu^2(d_j+\sqrt D)}{\mu^2\lambda_j^2n}
 +\frac{\nu s_j}{\mu\lambda_j^2}\sqrt{\frac{d_j}{n}}
 \bigg\}. \label{Cpoint35}
\end{equation}
Substituting \eqref{Cpoint35} into \eqref{Cpoint28}, and then combining with \eqref{Cpoint26} and \eqref{X1bound11}, proves \eqref{Cpoint7}.  A union bound over the three modes and the first-moment event gives total success probability at least
\(
 1-C\exp(-cT_{\mathrm{df}})-C\sum_{j=1}^3\exp(-cd_j)
 \geq1-C\exp(-c\underline d).\)
Putting pieces together completes the proof.
\end{proof}

\subsection{Proof of Corollary \ref{corCpoint1}}
\begin{proof}
Write $\Lambda:=\Lambda_{\min}(\tcalX^\star)$.  First, let $g\sim N(0,1)$, Cauchy--Schwarz and \eqref{Cpoint2} give 
\begin{equation}
 \mu\rho
 =|\mathbb{E}[gf(\rho g)]|\le \|g\|_{L_2}\|f(\rho g)\|_{L_2}\lesssim \nu. \label{Cpoint40}
\end{equation}
Also, for every mode $j$,
\begin{equation}
 \rho\geq s_j\geq\lambda_j\geq\Lambda,
 \qquad \kappa\geq1. \label{Cpoint41}
\end{equation}
Relations \eqref{Cpoint40}--\eqref{Cpoint41} show that the first term in \eqref{Cpoint8} implies $n\geq CT_{\mathrm{df}}$.  Moreover,
\[
 \frac{d_j}{\lambda_j^2}
 \leq\frac{\overline d}{\Lambda^2}
 \leq\frac{\rho^2\kappa^2\overline d}{\Lambda^4},
 \qquad
 \frac{\sqrt D}{\lambda_j^2}
 \leq\frac{\rho\sqrt D}{\Lambda^3},\qquad  \frac{s_j^2d_j}{\lambda_j^4}
 \leq\frac{\kappa^2\overline d}{\Lambda^2}
 \leq\frac{\rho^2\kappa^2\overline d}{\Lambda^4}.
\]
Thus, after increasing $C_{c_0}$, condition \eqref{Cpoint8} implies \eqref{Cpoint6}. We now bound the terms on the right-hand side of \eqref{Cpoint7}.  The three terms in \eqref{Cpoint8}, respectively, imply
\[
 \nu\sqrt{\frac{T_{\mathrm{df}}}{n}}
 \leq c\mu\Lambda,
 \qquad
 \frac{\rho\nu^2\sqrt D}{\mu\Lambda^2n}
 \leq c\mu\Lambda,\qquad  \frac{\rho\nu\kappa}{\Lambda}
 \sqrt{\frac{\overline d}{n}}
 \leq c\mu\Lambda.
\] 
The remaining contribution
$\rho\nu^2\overline d/(\mu\Lambda^2n)$ is also bounded by $c\mu\Lambda$ using the last term in \eqref{Cpoint8} and the fact that $\Lambda\leq\rho\kappa^2$.  Choosing the constant sufficiently large proves the first result in \eqref{Cpoint9}.

For nonzero tensors $\tcalU$ and $\tcalV$, the standard normalization inequality (see, e.g., \citealp[Equation (2.12)]{chen2024robust}) gives
\(
 \left\|\frac{\tcalU}{\|\tcalU\|_{\rm F}}-
 \frac{\tcalV}{\|\tcalV\|_{\rm F}}\right\|_{\rm F}
 \leq\frac{2\|\tcalU-\tcalV\|_{\rm F}}{\|\tcalV\|_{\rm F}}.\)
Taking $\tcalU=\tcalX_0$ and $\tcalV=\mu\tcalX^\star$, and using $c_0<1$ and $\Lambda\leq\rho$, proves the second result in \eqref{Cpoint9}.

It remains to verify the condition-number form.  For every mode $j$, we have
\(
 \rho^2=\|\calM_j(\tcalX^\star)\|_{\rm F}^2
 \leq r_j\|\calM_j(\tcalX^\star)\|_{\rm op}^2
 \leq\overline r\Lambda_{\max}^2(\tcalX^\star).\)
Consequently,
\begin{equation}
 \frac{1}{\Lambda}
 =\frac{\kappa}{\Lambda_{\max}(\tcalX^\star)}
 \leq\frac{\sqrt{\overline r}\,\kappa}{\rho}. \label{Cpoint42}
\end{equation}
Together with
$T_{\mathrm{df}}\leq\overline r^3+3\overline r\,\overline d$, relation \eqref{Cpoint42} yields
\begin{align*}
    \Big(\frac{\nu}{\mu}\Big)^2
 \frac{T_{\mathrm{df}}}{\Lambda^2}
 &\leq C\Big(\frac{\nu}{\mu\rho}\Big)^2
 \left(\kappa^2\overline r^4+
 \kappa^6\overline r^2\overline d\right),\\
  \Big(\frac{\nu}{\mu}\Big)^2
 \frac{\rho\sqrt D}{\Lambda^3}
 &\leq\Big(\frac{\nu}{\mu\rho}\Big)^2
 \kappa^3\overline r^{3/2}\sqrt D,\\
  \Big(\frac{\nu}{\mu}\Big)^2
 \frac{\rho^2\kappa^2\overline d}{\Lambda^4}
 &\leq\Big(\frac{\nu}{\mu\rho}\Big)^2
 \kappa^6\overline r^2\overline d.
\end{align*} 
Hence \eqref{Cpoint11} implies \eqref{Cpoint8}, completing the proof.
\end{proof}
\section{Technical Lemmas}
\begin{lem}
    \label{pertur} 
    For any $\mathbfcal{A}\in \mathbb{R}^{d_1\times d_2\times d_3}$ and $\mathbfcal{B}\in T^{\bd}_{\br}$, it holds that
    \begin{align*}
        \|H^\bd_\br(\mathbfcal{A})-\mathbfcal{B}\|_{\rm F} \le (\sqrt{3}+1)\|\mathbfcal{A}-\mathbfcal{B}\|_{\rm F}. 
    \end{align*}
\end{lem}
\begin{proof}
    This is an immediate consequence of $ \|H^{\bd}_{\br}(\mathbfcal{A})-\mathbfcal{A}\|_{\rm F}\le \sqrt{3}\|\calP_{T^{\bd}_\br}(\mathbfcal{A})-\mathbfcal{A}\|_{\rm F}$: 
    \begin{align*}
        &\|H^\bd_\br(\mathbfcal{A})-\mathbfcal{B}\|_{\rm F}\le\| H^\bd_\br(\mathbfcal{A})-\mathbfcal{A}\|_{\rm F}+\|\mathbfcal{A}-\mathbfcal{B}\|_{\rm F}\\
        &\le \sqrt{3}\|\mathcal{P}_{T^{\bd}_{\br}}(\mathbfcal{A})-\mathbfcal{A}\|_{\rm F}+\|\mathbfcal{A}-\mathbfcal{B}\|_{\rm F} \le (\sqrt{3}+1)\|\mathbfcal{A}-\mathbfcal{B}\|_{\rm F},
    \end{align*}
    where the last inequality follows from the definition of $\calP_{T^{\bd}_{\br}}(\cdot)$ and $\mathbfcal{B}\in T^{\bd}_{\br}$.
\end{proof}
\begin{lem}[{\citealp[Lemma 9]{luo2023low}} (see also \citealp{cai2020provable})]
    \label{rie} For any $\mathbfcal{A},\mathbfcal{B}\in\mathbb{R}^{d_1\times d_2\times d_3}$ of Tucker-rank $(r_1,r_2,r_3)$, we  have 
    \begin{align*}
        \|\calP_{T(\tcalB)^\bot}(\tcalA)\|_{\rm F}=\|\mathbfcal{A} - \calP_{T(\tcalB)}(\mathbfcal{A})\|_{\rm F}\le \frac{3\|\mathbfcal{A}-\mathbfcal{B}\|_{\rm F}^2}{\Lambda_{\min}(\mathbfcal{A})}~,
    \end{align*}
    where $T(\tcalB)$ denotes the tangent space of $\{\tcalU\in\mathbb{R}^{d_1\times d_2\times d_3}:\tucker(\tcalU)=(r_1,r_2,r_3)\}$ at $\tcalB$.  
\end{lem}
\begin{lem}[Concentration of Product Process; {\citealp[Theorem 1.13]{mendelson2016upper}}]\label{lem:product_process} {(see also \citealp[Theorem 8]{genzel2023unified}}.)
Let $\{g_{\bu}\}_{\bu\in \calU}$ and $\{h_{\bv}\}_{\bv\in\calV}$ be stochastic processes indexed by two sets    $\calU\subset \mathbb{R}^{p}$ and $\calV\subset \mathbb{R}^q$, both defined on a common probability space $(\Omega,A,\mathbbm{P})$. We assume that there exist $K_{\calU},K_{\calV},r_{\calU},r_{\calV}\geq 0$ such that 
    \begin{align*}
       & \|g_{\bu}-g_{\bu'}\|_{\psi_2}\leq K_{\calU}\|\bu-\bu'\|_2,~\|g_{\bu}\|_{\psi_2} \leq r_{\calU},~\forall\bu,\bu'\in \calU;\\
       & \|h_{\bv}-h_{\bv'}\|_{\psi_2} \leq {K}_{\calV}\|\bv-\bv'\|_2,~\|h_{\bv}\|_{\psi_2}\leq r_{\calV},~\forall\bv,\bv'\in \calV.
    \end{align*}
Suppose that $\ba_1,...,\ba_n$ are independent copies of a random variable $\ba\sim \mathbbm{P}$. There exist universal constants $C>c>0$ such that, for every $t\geq 1$,
\begin{equation}
   \begin{aligned}\nonumber
        &\sup_{\bu\in\calU}\sup_{\bv\in\calV} ~
        \Big|\frac{1}{n}\sum_{i=1}^n g_{\bu}(\ba_i)h_{\bv}(\ba_i)
        -\mathbbm{E}[g_{\bu}(\ba_i)h_{\bv}(\ba_i)]\Big|\\
        &\leq C\left(
        \frac{(K_{\calU}\cdot\omega(\calU)+t\cdot r_{\calU})
        \cdot (K_{\calV}\cdot \omega(\calV)+t\cdot r_{\calV})}{n}+\frac{r_{\calU}\cdot K_{\calV}\cdot \omega(\calV)
        +r_{\calV}\cdot K_{\calU}\cdot \omega(\calU)
        +t\cdot r_{\calU}r_{\calV}}{\sqrt{n}}
        \right)
   \end{aligned}
\end{equation}
holds with probability at least $1-2\exp(-ct^2)$. 
\end{lem}

\begin{lem}[Uniform hyperplane tessellation; \citealp{oymak2015near}]  \label{lem:globinary}
Drawing a matrix $\bA=[\ba_1,\cdots,\ba_n]^\top\in \mathbb{R}^{n\times d}$ with i.i.d. $N(0,1)$ entries, for any $\bu,\bv\in\mathbb{S}^{d-1}$ we let $  \bR_{\bu,\bv}:= \{i\in[n]:\sign(\ba_i^\top\bu)\ne\sign(\ba_i^\top\bv)\}$. 
    For any cone $\calC\subset\mathbb{R}^d$, there are universal constants $C>c>0$ such that, for any small enough $\delta>0$, if
    \begin{align*}
        n\ge C \bigg(\frac{\log\calN(\calC\cap \mathbb{S}^{d-1};c\delta/\sqrt{\log(e/\delta)})}{\delta^2}+\frac{\omega^2(\calC_{(1)})}{\delta\log(e/\delta)}\bigg),
    \end{align*}
    then with probability at least $1-2\exp(-cn\delta^2)$, we have
    \begin{align*}
        \sup_{\bu,\bv\in \calC\cap \mathbb{S}^{d-1}}\bigg|\frac{|\bR_{\bu,\bv}|}{n}-\frac{\arccos(\bu^\top\bv)}{\pi}\bigg|\le\delta.
    \end{align*}
\end{lem}
\begin{proof}
Such a global binary embedding result originated from \citet{oymak2015near}.
    The above statement is obtained from \citet[Theorem 1.5]{dirksen2022sharp} by setting parameters therein as $T=\calC\cap\mathbb{S}^{d-1}$, $\theta={c'\delta}/{\sqrt{\log(e/\delta)}}$, along with the observation $T_{(\theta)}:=(T-T)\cap \theta\mathbb{B}_2\subset (\calC-\calC)\cap\theta\mathbb{B}_2 = \theta \calC_{(1)}$.  
\end{proof}

\begin{lem}
\citep[Theorem 2.10]{dirksen2021non} \label{lem:max_ell_sum}
Let $\ba_1,...,\ba_n$ be independent   random vectors in $\mathbb{R}^d$ satisfying $\mathbbm{E}(\ba_i\ba_i^\top)=\bI_d$ and  $\max_i\|\ba_i\|_{\psi_2}\leq L$ (here, $\|\ba_i\|_{\psi_2}=\sup_{\bv\in \mathbb{S}^{d-1}}\|\bv^\top\ba_i\|_{\psi_2}$). For some   given $\calW\subset \mathbb{R}^d$ and  $1\leq \ell\leq n$, there exist  constants $C,c$ depending only on $L$ such that  
    \begin{equation*}     \sup_{\bw\in\calW}\max_{\substack{I\subset [n]\\|I|\leq \ell}}\bigg(\frac{1}{\ell}\sum_{i\in I}|\langle\ba_i,\bw\rangle|^2\bigg)^{1/2}\leq C\left(\frac{\omega(\calW)}{\sqrt{\ell}}+\rad(\calW)\log^{1/2}\Big(\frac{en}{\ell}\Big)\right)
    \end{equation*}
    holds with probability at least $1-2\exp(-c\ell\log(\frac{en}{\ell}))$, where $\rad(\calW)=\sup_{\bw\in\calW}\|\bw\|_2$. 
\end{lem}

\end{document}